\documentclass{article}
\usepackage{preamble}
\usepackage{longtable}

\title{Moment Relaxations for Data-Driven Wasserstein Distributionally Robust Optimization}

\author{Shixuan~Zhang \thanks{Wm Michael Barnes 64 Department of Industrial \& Systems Engineering,
		Texas A\&M University, College Station, TX, USA, 77843. (shixuan.zhang@tamu.edu).
	} \and
	Suhan~Zhong \thanks{
		Department of Mathematics,
		Texas A\&M University, College Station, TX, USA, 77843. (suzhong@tamu.edu). }}
\date{\vspace{-1cm}}

\begin{document}
	
\maketitle

\begin{abstract}
We propose moment relaxations for data-driven $p$-Wasserstein 
distributionally robust optimization ($p$-WDRO) problems that 
are defined by polynomials.
	The proposed moment relaxations admit Benders-type decomposition with 
	parallel evaluation of the subgradients using each sample subproblem, 
	which enables efficient solution for larger training sets.
We then identify conditions on $p$ and the defining polynomial degrees 
such that the proposed $k$-th order moment relaxations preserve 
the asymptotic consistency of the original $p$-WDRO 
(i.e., the relaxation gap is bounded at most linearly by
the Wasserstein radius). 
In particular, these conditions translate to effective bounds on $k$, 
which lead to polynomially sized semidefinite optimization formulations 
that are compatible with existing solvers.
Numerical experiments on
	a box-constrained regression problem and a two-stage production problem 
	are included to demonstrate the scalability and the effectiveness 
	of the proposed moment relaxations.
\end{abstract}

	\section{Introduction}
Distributionally robust optimization (DRO) is a framework for 
decision-making under uncertainty.
Given an \emph{ambiguity set} $\calA$ of probability measures describing 
the distribution of uncertain parameters $\xi$ over a set $\Xi\subseteq\bbR^{n_0}$, 
it aims to minimize the worst-case mean cost:
\begin{equation}\label{eq:DRO}
	\min_{x\in X}\, \Big\{ f(x)+\sup_{\mu\in\calA}\bbE_{\xi\sim\mu}F(x,\xi) \Big\},
\end{equation}
where $x$ is the decision vector constrained in a 
\emph{feasibility set} $X\subseteq\bbR^{n_1}$,
and $f:X\to\bbR$, $F:X\times\Xi\to\bbR$ are the cost functions.
Here, $\bbE_{\xi\sim\mu}$ denotes the expectation assuming that 
$\xi$ is a random vector associated with the probability measure $\mu\in\calA$. 
DRO has witnessed many successful applications, 
with the solution algorithms largely contingent upon the ambiguity set $\calA$ 
(see e.g.,~\cite{rahimian2022frameworks} for a comprehensive review).

Among various choices of ambiguity sets, 
increasing attention has been given to \emph{Wasserstein balls}.
That is, $\calA$ is defined as a ball of radius $r$, 
with respect to a chosen $p$-Wasserstein distance for probability measures on $\Xi$.
The center of $\calA$ is often an \emph{empirical probability measure} $\hat{\nu}$ 
constructed from $N$ independent and identically distributed (iid) samples 
$\hat{\xi}^{(1)},\dots,\hat{\xi}^{(N)}$ of the underlying truth.
DRO with such $p$-Wasserstein ball ambiguity is known as 
\emph{data-driven $p$-Wasserstein DRO} ($p$-WDRO). 
It offers two notable features:
\begin{itemize}
	\item (performance guarantee) 
	when $N$ is large, concentration in terms of Wasserstein distances dictates that 
	the out-of-sample mean cost is bounded by the in-sample one with high probability~\cite{fournier2015rate};
	\item (asymptotic consistency) 
	by decreasing $r\to0$, the optimal value and solution(s) of the DRO 
	respectively converge to those of the \emph{empirical stochastic optimization (ESO)} 
	based on $\hat{\nu}$:
	\begin{equation}\label{eq:ESO}
		\min_{x\in X}\, \Big\{ f(x)+\frac{1}{N}\sum_{i=1}^{N}F(x,\hat{\xi}^{(i)}) \Big\}.
	\end{equation}
\end{itemize}
Thus, by increasing the data size $N$ and decreasing the Wasserstein radius $r$ accordingly, 
 it is known in the literature that 
	one can avoid excessive conservatism while guaranteeing
	out-of-sample performance~\cite{mohajerin2018data,blanchet2021statistical,gao2023finite}.

Alongside the performance guarantee and the asymptotic consistency, 
data-driven Wasserstein DRO brings significant challenges to the solution procedure.
Based on Slater's condition for $r>0$, 
it is known that the $p$-WDRO problem can be reformulated as
\begin{equation}\label{eq:WDRO}
	\min_{x\in X,\lambda\in\bbR_{\ge0}}\, 
	\Big\{ 
	f(x)+r^p\lambda+\frac{1}{N}\sum_{i=1}^{N}\sup_{\xi\in\Xi}
	\big[F(x,\xi)-\lambda\nVert{\xi-\hat{\xi}^{(i)}}^p\big] 
	\Big\},
\end{equation}
where $\nVert{\cdot}$ is a distance function on $\Xi\subseteq\bbR^{n_0}$ and 
$p\in\bbZ_{\ge1}$ is the chosen order of the Wasserstein distance~\cite{gao2023distributionally}. 
In contrast with the ESO~\eqref{eq:ESO} where the evaluation of $F(x,\hat{\xi}^{(i)})$ 
is usually simple for any $x\in X$, 
the supremum in~\eqref{eq:WDRO} can be computationally intractable 
when $F(x,\cdot)$ is not concave.
For instance, this occurs when $F$ is the \emph{value function} 
in a two-stage optimization problem with fixed linear recourse, i.e.,
\begin{equation}\label{eq:Recourse}
	F(x,\xi) = \min_{x'\in\re^{m_2}}
	\Big\{ c(\xi)^\transpose x'+d(\xi): Ax'= B(\xi)x+b(\xi),\,x'\ge0 \Big\},
\end{equation}
where $A\in\bbR^{n_2\times m_2}$, $B:\Xi\mapsto\bbR^{n_2\times n_1}$, 
$b:\Xi\mapsto\bbR^{n_2}$, $c:\Xi\mapsto\bbR^{m_2}$ and 
$d:\Xi\mapsto\bbR$ are explicitly given.
For any fixed $\xi$, evaluation of $F$ in~\eqref{eq:Recourse} reduces to 
solving linear optimization, and thus the ESO~\eqref{eq:ESO} 
typically allows efficient solution methods~\cite{birge2011introduction,shapiro2021lectures}.
However, $F$ is generally not concave in $\xi$, 
due to the dependency of $b(\xi)$ and $B(\xi)$ on $\xi$, 
which leads to great challenges in solving~\eqref{eq:WDRO}.
Indeed, the decision problem associated with~\eqref{eq:WDRO} is \emph{NP-hard} 
for sufficiently large $r>0$, 
even when $c(\xi)=c$ is independent of $\xi$ and $B(\xi),b(\xi)$ 
are affine-linear in $\xi$~\cite{hanasusanto2018conic}.
The NP-hardness result is then shown for any $r>0$, $p>1$, 
$\nVert{\cdot}=\nVert{\cdot}_p$, $c(\xi)=c$, 
$B(\xi),b(\xi)$ dependent affine-linearly on $\xi$, and $\Xi=\bbR^{n_0}$~\cite{xie2020tractable}.
For these reasons, studies of data-driven $p$-WDRO have been mostly focused on special cases, e.g., $p=1$~\cite{luo2019decomposition,mohajerin2018data,zhao2018data,gamboa2021decomposition,duque2022distributionally,zhang2022distributionally,byeon2025two}, 
or $p=\infty$, for which a different reformulation from~(\ref{eq:WDRO}) is used~\cite{xie2020tractable,jiang2022distributionally}.
In particular, tractability in these cases often relies on 
	further assumptions on the objective function, 
	such as $F(x,\xi)$ being a maximum of concave function in $\xi$~\cite{mohajerin2018data}.
To our knowledge, \cite{hanasusanto2018conic} is the only prior work on two-stage $2$-WDRO: 
it derives a copositive optimization reformulation for~(\ref{eq:WDRO}), 
 which then leads to a semidefinite relaxation 
	that couples all training samples through a common decision vector $x$.

\subsection{Outline of the proposed moment relaxations}

In this work, we propose \emph{moment relaxations} for 
a wide class of data-driven $p$-WDRO problems~\eqref{eq:WDRO} 
as surrogates with better tractability, 
and then show that they preserve asymptotic consistency from the original $p$-WDRO.
Moment relaxations are commonly used in polynomial optimization problems \cite{lasserre2001global,parrilo2003semidefinite,blekherman2012semidefinite,NieBook}.
Due to their flexibility, they have been recently introduced to DRO 
with other types of ambiguity sets~\cite{de2020distributionally,nie2023distributionally,nie2025distributionally}.
To apply it in our $p$-WDRO problem, we focus on distance functions 
given by Euclidean norms $\nVert{\cdot}$, 
i.e., $\nVert{\xi}=(\xi^\transpose H\xi)^{1/2}$ 
for some positive definite matrix $H\succ0$ unless specified, 
and an even Wasserstein order $p>1$. 
	Then the $p$-th power of the distance function 
	\begin{equation}\label{eq:Q_p}
		Q_p^{(i)}(\xi)\coloneqq\nVert{\xi-\hat{\xi}^{(i)}}^p
		= \big[(\xi-\hat{\xi}^{(i)})^{\transpose} H (\xi-\hat{\xi}^{(i)})\big]^{p/2}
	\end{equation}
	is a sum-of-squares polynomial of degree $p$.
	Similarly, for any positive integer $d$, 
	the $2d$-th power of the standard $l_{2d}$-norm, defined as 
	$\|\xi\|_{2d}\coloneqq (\xi_1^{2d}+\cdots+\xi_{n_0}^{2d})^{1/2d}$,
	is also a sum-of-squares polynomial.
We next outline the proposed moment relaxations 
while postponing the background review to Section~\ref{sec:prelim}.

\noindent\textbf{Single-stage $p$-WDRO.}
	Suppose  $F(x,\xi) = \max\{F_1(x,\xi),\allowbreak\dots,F_T(x,\xi)\}$ 
	is the pointwise maximum of $T$ polynomials,
and $\Xi \coloneqq \{\xi\in\bbR^{n_0}: h(\xi)\ge0\}$ 
is defined by a polynomial tuple $h \coloneqq (h_1,\dots,h_{m_1})$.
 For each $t\in[T]$, let $F_{t,x}(\cdot) \coloneqq F_t(x,\cdot)$ 
denote an instantiation of $F_t$ at $x\in X$. 
We use $\bar{\scrS}[h]_{2k}$ to denote the set of normalized
\emph{truncated (pseudo-)moment sequence (tms)} 
of degree-$2k$ associated with $h$, 
which is semidefinite representable by $\binom{n_0+2k}{2k}$ variables 
and $(m_1+1)$ \emph{linear matrix inequalities (LMIs)};
see \eqref{eq:S_2k} for the formal definition.
The proposed moment relaxation for single-stage $p$-WDRO is
\begin{equation}\label{eq:1stage}
	\min_{x\in X,\lambda\in\bbR_{\ge0}}
	\Big\{ 
	f(x) + r^p\lambda +\frac{1}{N}\sum_{i=1}^{N} 
	\max_{t\in[T]}\sup_{y\in\bar{\scrS}[h]_{2k}}
	\pAngle{F_{t,x}-\lambda Q^{(i)}_p}{y}
	\Big\}.
\end{equation}

\noindent\textbf{Two-stage $p$-WDRO.}
By linear optimization duality, 
we can write the linear recourse \eqref{eq:Recourse} as
\begin{equation}\label{eq:RecourseDual}
	F(x,\xi) = \max_{u\in\re^{n_2}} 
	\Big\{ u^\transpose B(\xi)x+b(\xi)^\transpose u+d(\xi): 
	c(\xi)-A^\transpose u\ge0 \Big\}.
\end{equation}
For convenience, let $G_x(\xi,u) \coloneqq u^\transpose B(\xi)x+b(\xi)^\transpose u+d(\xi)$ 
and $g(\xi,u) \coloneqq c(\xi)-A^\transpose u$.
	The decision vector $(\xi,u)$ in the inner maximization of \eqref{eq:WDRO} 
	is subject to $g\ge 0$ and $h\ge 0$.
Similar to the single-stage case, 
	the $2k$-th degree monomial vector of each feasible pair $(\xi,u)$ 
	belongs to the semidefinite representable set $\bar{\scrS}[g,h]_{2k}$,
	which is defined by $\binom{n_0+n_2+2k}{2k}$ variables and $(m_1+m_2+1)$ LMIs.
The proposed moment relaxation for two-stage $p$-WDRO is
\begin{equation}\label{eq:2stage}
	\min_{x\in X,\lambda\in\bbR_{\ge0}}
	\Big\{ 
	f(x) + r^p\lambda +\frac{1}{N}\sum_{i=1}^{N} 
	\sup_{y\in\bar{\scrS}[g,h]_{2k}}\pAngle{G_x-\lambda Q^{(i)}_p}{y}
	\Big\}.
\end{equation}

If we fix the relaxation order $k$, 
the inner suprema in~\eqref{eq:1stage} and \eqref{eq:2stage} can be sought 
via polynomially sized semidefinite optimization subproblems independently for each 
	sample subproblem $i\in[N]$.
	Their solutions can then be used for computing 
	subgradients \eqref{eq:1stage-subgrad}-\eqref{eq:2stage-subgrad} below, respectively, 
	in Benders-type solution methods for the outer minimization 
	(see Algorithm~\ref{alg:Subgradient}).
	However, to ensure that the relaxation gap goes to zero, 
	existing analyses of moment relaxations for polynomial optimization 
	either require the relaxation order $k$ scale at least 
	as the inverse square root of the desired gap~\cite{slot2021near,laurent2024overview}, 
	resulting in exponentially sized moment relaxations; 
	or the diameter of the feasibility set (e.g., $\Xi$ in our case) to shrink~\cite{nie2007complexity,nie2013approximation}, 
	which is not applicable in the $p$-WDRO context.
Thus, to our knowledge, it was unknown whether 
the moment relaxations~\eqref{eq:1stage} and \eqref{eq:2stage} 
can be used for finding near-optimal solutions to the 
$p$-WDRO problems while staying polynomially sized.

\subsection{Overview of main results}

	Assume that the support set $\Xi$ is a \emph{regular closed set}, 
	i.e., it is the closure of its interior.
	To consider potentially unbounded support set $\Xi$,
	we make the following standard assumption on the empirical $p$-th order 
	moments of sample data throughout the paper.
	\begin{assumption}\label{assum}
		There exists a constant $\calB>0$ such that for any sample size $N$, 
		the samples $\hat{\xi}^{(1)}, \dots, \hat{\xi}^{(N)}$ satisfy 
		\[
		\frac{1}{N}\sum_{i=1}^{N}(1+\nVert{\hat{\xi}^{(i)}}_p)^p\le\calB.
		\]
\end{assumption} 

As relaxations, \eqref{eq:1stage} and \eqref{eq:2stage} 
automatically provide more conservative estimates of the mean cost 
compared to the original $p$-WDRO~(\ref{eq:WDRO}).
Consequently, whenever~(\ref{eq:WDRO}) has the out-of-sample performance guarantee, 
so do~\eqref{eq:1stage} and \eqref{eq:2stage}, respectively.
The motivating question of this paper is thus the following.

\vspace{0.05 in}
\begin{question}\label{Or_question}
	Under what conditions does any optimal decision $x$ to the $k$-th order
	moment relaxation~\eqref{eq:1stage} or~\eqref{eq:2stage} 
	have an $\calO(r)$-optimality gap as a decision 
	 in the original $p$-WDRO problem~\eqref{eq:WDRO}?
\end{question}
\vspace{0.05 in}

In Question~\ref{Or_question}, we use the notation $\calO(r)$ to indicate that 
the optimality gap can be bounded by $\scrC\cdot r$ for sufficiently small $r>0$, 
where the constant $\scrC>0$ is independent of $x$, $N$, 
	or the samples $\hat{\xi}^{(1)},\dots,\hat{\xi}^{(N)}$, 
	and may depend only polynomially on the problem coefficients,
	$n_0$ (the dimension of $\xi$),
	and the sample moment bound $\calB$ in Assumption~\ref{assum}. 
Our question is equivalent, up to an $\calO(r)$ difference, 
to the notion of {\it asymptotic consistency} such as in~\cite{mohajerin2018data}: 
	under Assumption~\ref{assum}, the difference between 
	the $p$-WDRO problem~\eqref{eq:WDRO} and the ESO problem~\eqref{eq:ESO} 
	is also bounded by $\calO(r)$~\cite{bartl2021sensitivity}, 
	where the latter is an asymptotically consistent estimator 
	with respect to the true distribution of $\xi$. 
We thus refer to the property in the question as 
\emph{$\calO(r)$-consistency} for simplicity.
	Answers to Question~\ref{Or_question} would, 
	for the first time in the literature, 
	extend the strategy of increasing the data size $N$ 
	while decreasing the radius $r$ 
	to the moment relaxations~\eqref{eq:1stage} and~\eqref{eq:2stage}.
	Moreover, the $\mathcal{O}(r)$-consistency is independent of 
	the specific radius-decreasing scheme, 
	making it broadly applicable to various classes of distributions for $\xi$.

For single-stage problems, we answer 
Question~\ref{Or_question} with the following theorem.
\begin{theorem}\label{thm:1stage-intro}
	Under Assumption~\ref{assum}, 
	suppose $X$ is bounded and the relaxation order 
	$k$ satisfies $2k\ge\max\{\deg(F_x)$, $\deg(h),$ $p\}$.
	Any optimal decision $x$ of the moment relaxation~\eqref{eq:1stage} 
	is an $\calO(r)$-optimal solution to the $p$-WDRO problem~\eqref{eq:WDRO}, 
	if either
	\begin{enumerate}
		\item[(i)] $p\ge\deg(F_x)$; or
		\item[(ii)] there exist $R>0$, $k_1\in\bbN$ such that $R-\nVert{\xi}^2\in\qmod{h}_{2k_1+2}$, 
		and $\deg(F_x)\le 2k-2k_1+p$.
	\end{enumerate}
\end{theorem}
Here, $\deg(h)$ refers to the highest degree of polynomials in $h = (h_1,\dots,h_{m_1})$,
	and $\deg(F_x)$ refers to the maximum of the degrees $\deg(F_{1,x}),\dots,\deg(F_{T,x})$ 
	for generic $x\in X$, i.e., on a Zariski open subset of $X$.
We prove a more elaborate version, Theorem~\ref{thm:asym1str0}, 
of~Theorem~\ref{thm:1stage-intro}, in~Section~\ref{sec:1stage}.
The necessity of the first condition $p\ge\deg(F_x)$ for unbounded $\Xi$ 
is illustrated through~Example~\ref{ex:1stage}.
	Moreover, when $F$ is a polynomial, i.e., $T=1$, 
	we refine the above $\calO(r)$-estimate of the optimality gap.
	In particular, if constraint functions $-h_1,\dots,-h_{m_1}$ 
	are certifiably convex with \emph{sum-of-squares (SOS) Hessian matrices}, 
	we show that any optimal decision $x$ in~\eqref{eq:1stage} has 
	an $\calO(r^2)$-optimality gap as a decision to the original 
	$p$-WDRO~\eqref{eq:WDRO} in Theorem~\ref{thm:Or2gap}.
	Examples on the necessity of the SOS Hessian matrices, 
	as well as more specialized cases for zero optimality gaps, 
	are presented in Section~\ref{ssc:tighter_gap}.

For two-stage problems, we answer Question~\ref{Or_question} with the following theorem.
	For convenience, denote $d_g\coloneqq\max\{ 2\deg(B), 2\deg(b),
	\deg(G_x), \deg(g)\}$ with $B,b,d$ in \eqref{eq:RecourseDual},
	where notations $\deg(G_x)$ and $\deg(B), \deg(b),\deg(g)$ are used in the same way 
	as $\deg(F_x)$ and $\deg(h)$, respectively.
\begin{theorem}\label{thm:2stage-intro}
	Under Assumption~\ref{assum}, 
	suppose $X$ is bounded and there exist $R_0>0$, $k_0\in\N$ such that 
	$R_0-\nVert{u}^2\in\qmod{g,h}_{2k_0}$.
	Let $k$ satisfy that $2k\ge\max\{d_g\deg(h),p,2k_0\}$.
	Then any optimal decision $x$ of the moment relaxation~\eqref{eq:2stage} 
	is an $\calO(r)$-optimal solution to the $p$-WDRO problem~\eqref{eq:WDRO}, 
	if either
	\begin{enumerate}
		\item[(i)] $p\ge d_g$; or
		\item[(ii)] there exist $R>1$, $k_1\in\bbN$ such that 
		$R-\nVert{\xi}^2\in\qmod{h}_{2k_1+2}$ and $d_g\le 2k-2k_1+p$.
	\end{enumerate}
\end{theorem}
We prove a more elaborate version, Theorem~\ref{thm:2stage_idc}, 
of~Theorem~\ref{thm:2stage-intro}, in~Section~\ref{sec:2stage-standard}.
The assumption $R_0-\nVert{u}^2\in\qmod{g,h}_{2k_0}$ 
in Theorem~\ref{thm:2stage-intro} implies that the recourse dual variables $u$ are bounded.
In~Example~\ref{ex:2stage}, we show that this assumption cannot be removed in general, 
as the moment relaxation may have an infinite gap when $u$ is not bounded.
For this reason, we propose a strengthened moment relaxation that 
can preserve the $\calO(r)$-consistency for unbounded 
recourse dual variables in~Section~\ref{sec:2stage-strengthened}.
	Then in Section~\ref{sec:num}, we illustrate the scalability and 
	the effectiveness of the proposed moment relaxations on 
	a box-constrained regression problem and a two-stage production problem.
	Concluding remarks and future research directions are presented in Section~\ref{sec:con}.

\section{Preliminaries}\label{sec:prelim}

In this section, we review some preliminaries on data-driven $p$-WDRO reformulations, 
including a brief derivation of~\eqref{eq:WDRO},
	and some background on the moment relaxations~\eqref{eq:1stage} and \eqref{eq:2stage}. 
	We then outline Benders-type solution methods for these formulations.

\subsection{Wasserstein DRO}\label{sec:WDRO}
Let $\calM(\Xi)$ denote the set of Borel probability measures on $\Xi$.
Given any norm $\nVert{\cdot}$ on $\bbR^{n_0}$ and an integer $p\ge1$, 
define 
\begin{equation*}
	\calW_p(\Xi) \coloneqq 
	\Big\{
	\mu\in\calM(\Xi): \int_{\Xi}\nVert{\xi-\bar{\xi}}^p\diff\mu(\xi)<\infty
	\text{ for some }\bar{\xi}\in\Xi \Big\}
\end{equation*}
as the set of all probability measures that have finite $p$-th order moments on $\Xi$.
Then the \emph{$p$-Wasserstein distance} on $\calW_p(\Xi)$ can be defined as
\[
W_p(\alpha,\beta) \coloneqq \inf_{\gamma\in\calM(\Xi\times\Xi)}
\Big\{ \int_{\Xi\times\Xi} \nVert{\xi-\eta}^p\diff\gamma(\xi,\eta):
\begin{array}{c}
	\gamma\text{ is a joint probability}\\
	\text{measure of }\alpha\text{ and }\beta
\end{array}
\Big\}.
\]
Using these notations, the Wasserstein ball in \eqref{eq:DRO} is
\begin{equation}\label{eq:WassBall}
	\calA = \calB_{p,r}(\hat{\nu}) \coloneqq 
	\{\mu\in\calW_p(\Xi):W_p(\mu,\hat{\nu})\le r\},
\end{equation}
where $r>0$ is a predetermined radius and 
$\hat{\nu}\coloneqq \frac{1}{N}\sum_{i=1}^{N}\delta_{\hat{\xi}^{(i)}}$ 
is an empirical probability measure based on the samples 
$\hat{\xi}^{(1)},\ldots, \hat{\xi}^{(N)}$ drawn from the true probability measure $\nu$.
The term $\delta_{\hat{\xi}^{(i)}}$ represents the Dirac atomic measure 
supported at $\hat{\xi}^{(i)}$ for each $i\in[N]\coloneqq\{1,\dots,N\}$.
Given a decision $x\in X$, we consider
\[
\begin{array}{rl}
	\mbox{\emph{in-sample} mean cost:} & \quad f(x)+\sup_{\mu\in\calB_{p,r}(\hat{\nu})}\bbE_{\xi\sim\mu}F(x,\xi),\\
	\mbox{\emph{out-of-sample} mean cost:} & \quad 
	f(x)+\bbE_{\xi\sim\nu}F(x,\xi),
\end{array}
\]
and estimate the latter through resampling of the truth $\nu$ in experiments.
The concentration result~\cite{fournier2015rate} ensures that
$\nu\in\calB_{p,r}(\hat{\nu})$ with high probability
when $Nr^{n_0}$ is sufficiently large, 
which implies the out-of-sample mean cost is bounded by the in-sample one.

When $\Xi$ is an infinite set, $\calM(\Xi)$ is infinite-dimensional, 
so it is generally difficult to seek the supremum over $\calA$ numerically. 
Instead, people use the finite-dimensional reformulation~\eqref{eq:WDRO}, 
under the assumption that $F_x$ has a finite order-$p$ growth rate 
$\Gamma_p(x)$ on $\Xi$ for each $x\in X$.
Precisely, if $\Xi$ is bounded, then $\Gamma_p(x)=0$ by convention; otherwise,
\begin{equation*}    
	\Gamma_p(x) \coloneqq \limsup_{\xi\in\Xi:\nVert{\xi-\bar{\xi}}\to\infty}
	\frac{F(x,\xi)-F(x,\bar{\xi})}{\nVert{\xi-\bar{\xi}}^p}
\end{equation*}
for any fixed $\bar{\xi}\in\Xi$.
One can check that the growth rate is independent of the choice of $\bar{\xi}\in\Xi$.
Next, we outline the derivation of~\eqref{eq:WDRO}, 
which is a simpler version of the arguments in~\cite{gao2023distributionally}.
For any joint probability measure $\gamma\in\calM(\Xi\times\Xi)$ of $\mu$ and $\hat{\nu}$, 
where $\hat{\nu}$ is an empirical measure,
using conditional probability measures 
(or more generally \emph{disintegration}~\cite{chang1997conditioning}), 
there exist $\mu^{(1)},\dots,\mu^{(N)}$ such that 
$\gamma=\frac{1}{N}\sum_{i=1}^{N}\mu^{(i)}\otimes\delta_{\hat{\xi}^{(i)}}$, 
and consequently $\mu=\frac{1}{N}\sum_{i=1}^{N}\mu^{(i)}$.
For any $\rho>0$, let
\begin{equation}\label{eq:MeasConstr}
	v_p(x;\rho) \coloneqq \Big\{
	\sup_{\substack{\mu^{(i)}\in\mc{M}(\Xi)\\i\in[N]}}
	\frac{1}{N}\sum_{i=1}^{N}\int_{\Xi}F(x,\xi)\diff\mu^{(i)}(\xi):
	\frac{1}{N}\sum\limits_{i=1}^N Q_p^{(i)}(\xi) 
	\le\rho\Big\}.
\end{equation}
For notational consistency, we set
\[
v(x;0)=v_p(x;0)=\bbE_{\xi\sim\hat{\nu}}F(x,\xi).
\]
Since $(\delta_{\hat{\xi}^{(1)}}, \ldots, \delta_{\hat{\xi}^{(N)}})$ 
is always feasible to \eqref{eq:MeasConstr},
$v_p(x;\rho)\ge v(x;0) >-\infty$ for all $x\in X$ and $\rho\ge 0$.
Suppose $\rho = r^p>0$. For all $0<\epsilon<r^p$,
\[
v_p(x;r^p-\epsilon)
\le \sup_{\mu\in\mc{A}=\calB_{p,r}(\hat{\nu})}\bbE_{\xi\sim \mu}F_x(\xi)
\le v_p(x;r^p).
\]
The (generalized) Slater's condition implies the strong Lagrangian duality, i.e.,
\[
\begin{aligned}
	v_p(x;\rho) = & \min_{\lambda\in\re_{\ge0}} 
	\Big\{ \lambda\rho+\frac{1}{N}\sum_{i=1}^{N}\sup_{\mu^{(i)}\in\calM(\Xi)} 
	\int_{\Xi} [F(x,\xi)- \lambda Q_p^{(i)}(\xi)]\diff\mu^{(i)}(\xi) \Big\} \\
	=& \min_{\lambda\in\re_{\ge0}} 
	\Big\{ \lambda\rho+\frac{1}{N}\sum_{i=1}^{N}\sup_{\xi^{(i)}\in\Xi}
	\big[F(x,\xi^{(i)})-\lambda Q_p^{(i)}(\xi^{(i)})\big] \Big\},
\end{aligned}
\]
where the second equality holds by the choices of Dirac atomic measures 
$\mu^{(i)} = \delta_{\xi^{(i)}}$.
Thus, $v_p(x;\rho)$ is concave in $\rho$, and has finite values 
since it is bounded from above by setting $\lambda=\Gamma_p(x)$.
Consequently, $v_p(x;\cdot)$ is continuous on the open interval $(0,+\infty)$, 
which implies that $v_p(x;r^p)=\sup_{\mu\in\calA}\bbE_{\xi\sim \mu}F(x,\xi)$ 
for any $r>0$, and hence the reformulation~\eqref{eq:WDRO}.

\subsection{Moment relaxations}\label{sec:MomRelax}

The tractability of moment relaxations is closely related to 
certifying nonnegativity of polynomials through sums of squares. 
Let $z$ be the vector of $l$ indeterminates, 
which may represent $\xi$ or $u$ in our context.
We use $\bbR[z]$ to denote the real polynomial ring in $z$, 
and $\bbR[z]_t$ its subset of polynomials with degrees up to $t$.
Let $\N$ denote the set of all nonnegative integers.
Given $\alpha=(\alpha_1,\ldots, \alpha_l)\in\N^l$, 
a monomial with exponent $\alpha$ is written as 
$z^{\alpha} \coloneqq z_1^{\alpha_1}\cdots z_l^{\alpha_l}$ 
with degree $\aVert{\alpha} \coloneqq \alpha_1+\cdots+\alpha_l$.
For $t\in\N$, denote $\N_t^l \coloneqq \{(\alpha_1,\ldots,\alpha_l)\in\N^l: \aVert{\alpha}\le t\}$.
In the context, we will use
$[z]_t \coloneqq [ 1, z_1, \cdots, z_l, z_1^2, z_1z_2, \cdots, z_l^t]$
to denote the vector of all monomials of $z$ up to degree $t$ 
(with the degree-lexicographic order).
We denote the set of all \emph{sum-of-square (SOS) polynomials} by $\Sigma[z]$ 
and its degree-$t$ truncation $\Sigma[z]_t \coloneqq \Sigma[z]\cap\bbR[z]_t$.
For a set $Z\subseteq\bbR^l$, the set of polynomials that are nonnegative on $Z$ 
is denoted as 
$\scrP(Z) \coloneqq \{q\in\bbR[z]:q(z)\ge0,\,\forall z\in Z\}$, 
with its degree-$t$ truncation defined as 
$\scrP_t(Z) \coloneqq \scrP(Z)\cap\bbR[z]_t$.
If $Z = \{z\in\bbR^l: g(z)\ge0\}$ for a polynomial tuple $g = (g_1,\ldots, g_m)$, 
then any polynomial that can be written as $\sigma_0+\sigma_1g_1+\cdots+\sigma_mg_m$ 
for some $\sigma_0,\dots,\sigma_m\in\Sigma[z]$ is manifestly nonnegative on $Z$.
We thus define the degree-$t$ truncated \emph{quadratic module} generated by $g$ as
\begin{equation}\label{eq:QMg}
	\qmod{g}_{t} \coloneqq 
	\Sig[z]_{t} + g_1\cdot\Sig[z]_{t-\deg(g_1)} + \cdots +
	g_{m}\cdot \Sig[z]_{t-\deg(g_{m})}.
\end{equation}
Both $\qmod{g}_t$ and $\scrP_t(Z)$ are closed under addition and nonnegative scaling, 
and thus form convex cones in the $\bbR$-vector space $\bbR[z]_t$.
When $Z$ has a nonempty interior, the cone $\qmod{g}_t$ is closed~\cite[Theorem~2.5.2]{NieBook}. 

A probability measure $\mu\in\calM(Z)$ leads to a sequence of moments
$(y_\alpha)_{\alpha\in\bbN^l_t}$ 
where $y_\alpha=\int_Z z^\alpha\diff\mu$ for each $\alpha\in\bbN^l_t$.
Such a sequence determines the expectation of any polynomial 
$q = \sum_{\alpha\in\bbN^l_t}q_\alpha z^\alpha\in\bbR[z]_t$ by
\begin{equation}\label{eq:bilinear}
	\pAngle{q}{y} \coloneqq \sum_{\alpha\in\N_t^l} q_{\alpha}y_{\alpha}
	= \int_Z q\diff\mu. 
\end{equation}	
By selecting a degree bound $t$, an optimization over measures, 
e.g.,~\eqref{eq:MeasConstr}, can be approximated by a finite sequence.
In fact, when $Z$ is compact, $y\in \re^{\alpha\in\bbN^l_t}$ admits a 
(not necessarily probability) measure supported on $Z$ if and only if 
$\pAngle{q}{y}\ge0$ for all $q\in\scrP_t(Z)$~\cite[Theorem 17.3]{schmudgen2017moment}.
However, it is usually difficult to check $\pAngle{q}{y}\ge0$ 
for all nonnegative polynomials $q\in\scrP_t(Z)$.
Instead, it is relatively convenient to check the condition only for 
manifestly nonnegative polynomials, such as those in $\scrQ[g]_t$.
We say that $y=(y_\alpha)_{\alpha\in\bbN^l_t}$ is a 
\emph{truncated (pseudo-)moment sequence (tms)} 
if it lies in the dual cone $\scrS[g]_t\coloneqq\qmod{g}_t^*$ of $\qmod{g}_t$ 
under the bilinear relation~\eqref{eq:bilinear},
	In particular, if $g\ge 0$ determines a set with nonempty interior, 
	then $\qmod{g}_{t}$ is also a proper cone, 
	thus $\qmod{g}_t = \mathscr{S}[g]_t^*$.
	This duality a well-established concept in polynomial optimization \cite{lasserre2001global,laurent2009sums,NieBook}.
	For convenience,
we further denote the set of all tms with the 
\emph{normalizing condition} $\bar{\scrS}[g]\coloneqq\scrS[g]\cap\{y_0=1\}$,
	of which the $2k$-th order truncation is explicitly 
	displayed in \eqref{eq:S_2k}.

In terms of computation, conditions for a real sequence $(y_\alpha)_{\alpha\in\bbN^l_t}$ 
to be a tms can be represented by LMIs using \emph{localizing matrices}. 
Given any $q\in\bbR[z]_{2t}$ and $k\ge t$, 
the $k$-th \emph{localizing matrix} of $q$ is the symmetric matrix, 
denoted as $L_{q}^{(k)}[y]$, such that
\begin{equation} \label{locmat:gi}
	\vectorize(a_1)^\transpose L_{q}^{(k)}[y]\vectorize(a_2)= \pAngle{qa_1a_2}{y},
\end{equation}
for all polynomials $a_1,a_2\in\bbR[z]_{k-t}$.
Here, $\vectorize(a_1),\vectorize(a_2)$ represent the coefficient vectors of $a_1,a_2$, 
in terms of the monomial basis $[z]_{k-t}$, respectively.
For the special case that $q=1$ is the constant polynomial, 
the localizing matrix is called $k$-th order \emph{moment matrix} of $y$:
\begin{equation}\label{eq:mommat}
	M_k[y]\,\coloneqq\, L_1^{(k)}[y].
\end{equation} 
For each relaxation order $k$ such that $2k\ge \deg(g)$,
the set of normalized tms can be written as
\begin{equation}\label{eq:S_2k}
	\bar{\scrS}[g]_{2k} \coloneqq \Big\{y\in\re^{\N_{2k}^l}: 
	y_0 = 1, M_k[y]\succeq 0, 
	L_{g_1}^{(k)}[y]\succeq 0, \ldots, L_{g_m}^{(k)}[y]\succeq 0\Big\}.    
\end{equation}
	Suppose $F$ is a polynomial.
	For a given pair $(x,\lambda)\in X\times \re_{\ge 0}$, 
	the $t$-th inner maximization problem in \eqref{eq:WDRO} 
	is a polynomial optimization problem.
	It can be equivalently reformulated as 
	\[
	\sup\limits_{\xi\in \Xi}\, [F_{x}-\lambda Q_p^{(i)}](\xi)
	= \sup\limits_{\mu\in\mc{M}(\Xi)}\, \int_{\Xi} 
	[F_{x}-\lambda Q_p^{(i)}](\xi){\tt d}\mu(\xi).
	\]
	Suppose $d=\deg(F_x-\lambda Q_p^{(i)})$.
	For any $\mu\in\mc{M}(\Xi)$, if $y\in\re^{\N_{d}^{n_0}}$ equals 
	the $d$-th order truncated moment of $\mu$, then 
	\[
	\int_{\Xi} [F_x-\lambda Q_p^{(i)}](\xi){\tt d}\mu(\xi) 
	= \langle F_x-\lambda Q_p^{(i)}, y\rangle,
	\]
	where $\langle \cdot, \cdot\rangle$ is the bilinear form be defined in \eqref{eq:bilinear}.
	For all integers $2k\ge d$, the $i$-th inner maximization problem in \eqref{eq:WDRO} 
	is thus equivalent to
	\[
	\sup\limits_{y}\,\langle F_x-\lambda Q_p^{(i)},y\rangle\quad \st\quad 
	\text{$y\in\re^{\N_{2k}^{n_0}}$ admits a representing measure in $\mc{M}(\Xi)$}.
	\] 
We then provide an explicit example of the moment relaxation~\eqref{eq:1stage}.

\begin{example}\label{ex:separable}
	Consider $X=\{0\}$, $f(x)\equiv 0$, $\Xi=[-1,1]^2$, 
	with $N=1$, $\hat{\xi}^{(1)}=(0,0)$, and $F_x(\xi)=\xi_1^3+\xi_2$.
	Since $\Xi$ is bounded, $\Gamma_p(x)=0$ for any $p\ge1$. 
	Set $p=2$ and suppose $r<\frac{1}{2}$.
	Since $X$ is a singleton, the data-driven $p$-WDRO~\eqref{eq:1stage} 
	reduces to an evaluation at $x=0$:
	\begin{equation*}
		\min_{\lambda\in\re_{\ge 0}}\Big\{\lambda r^2+\sup_{\xi\in\Xi}\big[\xi_1^3+\xi_2-\lambda(\xi_1^2+\xi_2^2)\big]\Big\}
		= \min\limits_{\lambda\in\re_{\ge 0}} \Big\{\lambda r^2+\frac{1}{4\lambda}\Big\} = r,
	\end{equation*}
	where the minimum is attained by $\lambda=\frac{1}{2r}$.
	Since $F_x(\hat{\xi}^{(1)}) = 0$, the difference between the objective values 
	of~\eqref{eq:ESO} and~\eqref{eq:WDRO} is exactly $r$.
	For $h=(1-\xi_1^2, 1-\xi_2^2)$ and relaxation order $k=2$,
	the decision vector $y=(1,y_{10},\ldots, y_{04})\in\bar{\mathscr{S}}[h]_4$ 
	must satisfy $M_2[y]\succeq 0, L_{1-\xi_1^2}^{(2)}[y]\succeq 0, L_{1-\xi_2^2}^{(2)}[y]\succeq 0$. 
	where
	\begin{equation*}
		\begin{aligned}
			M_2[y] &=\begin{bmatrix}
				1      & y_{10} & y_{01} & y_{20} & y_{11} & y_{02} \\
				y_{10} & y_{20} & y_{11} & y_{30} & y_{21} & y_{12} \\
				y_{01} & y_{11} & y_{02} & y_{21} & y_{12} & y_{03} \\
				y_{20} & y_{30} & y_{21} & y_{40} & y_{31} & y_{22} \\
				y_{11} & y_{21} & y_{12} & y_{31} & y_{22} & y_{13} \\
				y_{02} & y_{12} & y_{03} & y_{22} & y_{13} & y_{04}
			\end{bmatrix},
		\end{aligned}
		\begin{aligned}
			L_{1-\xi_1^2}^{(2)}[y] &= \begin{bmatrix}
				1-y_{20}      & y_{10}-y_{30} & y_{01}-y_{21} \\
				y_{10}-y_{30} & y_{20}-y_{40} & y_{11}-y_{31} \\
				y_{01}-y_{21} & y_{11}-y_{31} & y_{02}-y_{22}
			\end{bmatrix},\\[2mm]
			L_{1-\xi_2^2}^{(2)}[y] &= \begin{bmatrix}
				1-y_{02}      & y_{10}-y_{12} & y_{01}-y_{03} \\
				y_{10}-y_{12} & y_{20}-y_{22} & y_{11}-y_{13} \\
				y_{01}-y_{03} & y_{11}-y_{13} & y_{02}-y_{04}
			\end{bmatrix}.
		\end{aligned}
	\end{equation*}
	Since $M_2[y]$ is psd, we get $y_{30}^2\le y_{20}y_{40}\le y_{20}^2$, 
	and $y_{01}^2\le y_{02}\le 1$.
	Similarly, by $L_{1-\xi_1^2}^{(2)}[y]\succeq0$, we deduce that $y_{40}\le y_{20}$, 
	and by $L_{1-\xi_2^2}^{(2)}[y]\succeq0$ that $y_{02}\le1$.
	Thus, the problem~\eqref{eq:1stage} is bounded by
	\begin{equation*}\begin{aligned}
			&\min_{\lambda\in\re_{\ge 0}}\Big\{
			\lambda r^2+\sup_{y\in \bar{\mathscr{S}}[g]_4}\big[y_{30}+y_{01}-\lambda(y_{20}+y_{02})\big]\Big\}\\
			\le & \min_{\lambda\in\re_{\ge 0}}\Big\{ 
			\lambda r^2+\sup_{0\le y_{20}\le 1,y_{01}^2\le 1}
			[(1-\lambda)y_{20}+y_{01}-\lambda y_{01}^2]
			\Big\}.
		\end{aligned}
	\end{equation*} 
	By setting $\lambda=\frac{1}{2r}>1$, we see that this bound becomes 
	$y_{01}-\frac{y_{01}^2}{2r}\le \frac{r}{2}$, 
	where the equality is achieved at $y_{01} = r$.
	So the objective of~\eqref{eq:1stage} in this case is at most $r$, 
	which means that the moment relaxation is exact.
\end{example}

	\subsection{Solution methods}
	
	A key feature of the proposed moment relaxations~\eqref{eq:1stage} and~\eqref{eq:2stage} 
	is that they admit Benders-type solution methods,
	which can be easily distributed over each sample subproblem $i\in[N]$.
	To be more specific, let
	\begin{equation*}
		\calV_{p,2k}^{(i)}(x,\lambda)\coloneqq
		\begin{cases}
			\displaystyle\max_{t\in[T]}\sup_{y\in\bar{\scrS}[h]_{2k}}
			\pAngle{F_{t,x}-\lambda Q_p^{(i)}}{y}, \ 
			& \text{for single-stage relaxation~\eqref{eq:1stage};} \\[5mm]
			\displaystyle\sup_{y\in\bar{\scrS}[g,h]_{2k}}
			\pAngle{G_x-\lambda Q_p^{(i)}}{y}, \
			& \text{for two-stage relaxation~\eqref{eq:2stage}}. 
		\end{cases}
	\end{equation*}
	The minimization in both moment relaxations can then be unified as
	\begin{equation}\label{eq:MainMin}
		\min_{x\in X,\lambda\in\re_{\ge0}}
		f(x) + \lambda r^p + \frac{1}{N}\sum_{i=1}^{N} \calV_{p,2k}^{(i)}(x,\lambda).
	\end{equation}
	As commonly assumed in the literature (see e.g.,~\cite{mohajerin2018data}), 
	let us focus on the case where the set $X$ and the cost function $f(x)$ are convex.
	If $\calV_{p,2k}^{(1)},\dots,\calV_{p,2k}^{(N)}$ are also convex, 
	then the minimization problem~\eqref{eq:MainMin} 
	is amenable to \emph{first-order methods} such as the level method~\cite{lemarechal1995new},
	which only use the subgradients of $f$ and $\calV_{p,2k}^{(i)}$.
	In particular, through standard convex analysis, 
	a subgradient of $\calV_{p,2k}^{(i)}$ at a given point 
	$(\bar{x},\bar{\lambda})\in X\times\bbR_{\ge0}$ can be constructed as follows.
	In the single-stage case~\eqref{eq:1stage},
	for any maximizer $t^*\in[T]$ and its corresponding tms 
	$y^{(1)},\dots,y^{(N)}\in\bar{\scrS}[h]_{2k}$,
	\begin{equation}\label{eq:1stage-subgrad}
		\varsigma^{(i)} \coloneqq 
		\Big(\pAngle{\nabla_xF_{t^*}(\bar{x},\xi)}{y^{(i)}},\,
		-\pAngle{Q_p^{(i)}(\xi)}{y^{(i)}}\Big)
		\in \partial\calV_{p,2k}^{(i)}(\bar{x},\bar{\lambda}).
	\end{equation}
	In the two-stage case~\eqref{eq:2stage},
	for any optimal tms $y^{(1)},\dots,y^{(N)}\in\bar{\scrS}[g,h]_{2k}$,
	\begin{equation}\label{eq:2stage-subgrad}
		\varsigma^{(i)} \coloneqq
		\Big(\pAngle{B(\xi)^\transpose u}{y^{(i)}},\, 
		-\pAngle{Q_p^{(i)}(\xi)}{y^{(i)}}\Big)
		\in \partial\calV_{p,2k}^{(i)}(\bar{x},\bar{\lambda}).
	\end{equation}
	Then through aggregation over $i\in[N]$, 
	we can construct a subgradient of the sum 
	$f(x)+\lambda r^p+\frac{1}{N}\sum_{i=1}^{N}\calV_{p,2k}^{(i)}(x,\lambda)$ 
	for solving~\eqref{eq:MainMin}, as summarized in Algorithm~\ref{alg:Subgradient}.
	\begin{algorithm}[h]
		\caption{Evaluation of subgradients for the minimization~\eqref{eq:MainMin}}
		\label{alg:Subgradient}
		\begin{algorithmic}[1]
			\Require{a feasible point $(\bar{x},\bar{\lambda})\in X\times\bbR_{\ge0}$ 
				and a subgradient $\nabla{f}(\bar{x})\in\partial{f}(\bar{x})$}
			\Ensure{a subgradient vector of the function 
				$f(x)+\lambda r^p+\frac{1}{N}\sum_{i=1}^{N}\calV_{p,2k}^{(i)}(x,\lambda)$}
			\For{$i\in[N]$}\Comment{distributable over $i$}
			\State{solve the maximization 
				$$y^{(i)}\in\argmax_{y\in\bar{\scrS}[h]_{2k}}\pAngle{F_{t,\bar{x}}-\bar{\lambda} Q_p^{(i)}}{y}\text{ for each $t\in[T]$ or } 
				y^{(i)}\in\argmax_{y\in\bar{\scrS}[g,h]_{2k}}\pAngle{G_{\bar{x}}-\bar{\lambda} Q_p^{(i)}}{y}$$}
			\State{construct $\varsigma^{(i)}$ using~\eqref{eq:1stage-subgrad} or~\eqref{eq:2stage-subgrad}}
			\EndFor
			\State{return $$(\nabla{f}(\bar{x}),r^p)+\frac{1}{N}\sum_{i=1}^{N}\varsigma^{(i)}$$}
		\end{algorithmic}
	\end{algorithm}
	
	Now we reexamine the premise of Algorithm~\ref{alg:Subgradient} 
	that each $\calV_{p,2k}^{(i)}$ is convex. 
	This is always true for the two-stage case~\eqref{eq:2stage} 
	since each $\calV_{p,2k}^{(i)}$ is a maximum of linear functions $G_x-\lambda Q_p^{(i)}$ 
	in terms of $(x,\lambda)$.
	However, due to the nonlinearity of $F_x$ in $x$, 
	$\calV_{p,2k}^{(i)}$ may fail to be convex in the single-stage case~\eqref{eq:1stage}, 
	even when $T=1$ and  $F_1(\cdot,\xi)$ is convex for every $\xi\in\Xi$.
	Here is an example.
	\begin{example}\label{ex:NonconvexVx}
		Consider $X=[-1,1]$, $f(x)=0$, $\Xi=\{\xi\in\bbR^n:h_1(\xi)=1-\|\xi\|_2^2\ge0\}$, 
		$N=1$, and $\xi^{(1)}=0$.
		Let $\Phi(\xi)$ be any nonnegative form of degree $p\ge4$ 
		that is \emph{not} a sum of squares 
		(e.g., the Motzkin polynomial with degree $p=6$~\cite{motzkin1967arithmetic}),
		and let $F(x,\xi) = F_1(x,\xi) = \Phi(\xi)(x^2-1)$, 
		This $F$ is convex for every $\xi$ since $\Phi(\xi)\ge0$.
		However, $\calV_{p,2k}^{(1)}(x,\lambda)$ is not convex in $(x,\lambda)$
		for $2k=p$.
		This is because $\calV_{p,2k}^{(1)}(x,0)>0$ 
		and $\calV_{p,2k}^{(1)}(-1,0)=\calV_{p,2k}^{(1)}(1,0)=0$ for every $x\in(-1,1)$.
		To see this, let $x\in(-1,1)$ be arbitrary. 
		By polynomial optimization duality introduced in Subsection~\ref{sec:MomRelax}, 
		the optimal value 
		\[
		\calV_{p,2k}^{(i)}(x,0)
		= \min\{\gamma\in\bbR:(1-x^2)\Phi(\xi)+\gamma\in\scrQ[h_1]_p\}.
		\]
		In particular, the minimum is attained 
		since $\scrQ[h_1]_p$ is a closed cone in $\bbR[\xi]_{p}$, 
		due to the nonempty interior of $\Xi$.
		Then it suffices to show that $\gamma=0$ is not feasible in the above minimization: 
		by scaling this would have implied that 
		\[
		\Phi(\xi)=\sigma_0(\xi)+\sigma_1(\xi)h_1(\xi)\quad
		\mbox{for some}\quad 
		\sigma_0\in\Sigma[\xi]_p,\,\sigma_1\in\Sigma[\xi]_{p-2}. 
		\]
		Since $\Phi$ is homogeneous of degree $p$, 
		it vanishes at $\xi=0$ to the order of $p$, 
		while $h_1(0)=1$ and $\deg(\sigma)\le p-2$ force $\sigma_1(\xi)=0$.
		This implies that $\Phi(\xi)=\sigma_0(\xi)$, 
		which is a contradiction with $\Phi$ not being a sum of squares.
	\end{example}
	
	Nevertheless, there remains a large class of convex polynomial cost functions $F(x,\xi)$, 
	including those arising from least-squares and least-absolute-deviations regression, 
	for which each $\calV_{p,2k}^{(i)}$ is still convex.
	\begin{proposition}\label{prop:ConvexVx}
		Suppose $F_t(x,\xi) = \psi_t(a_t(\xi)^\transpose x+b_t(\xi))$ 
		for some polynomial maps $a_t:\Xi\to\bbR^{n_1}$, $b_t:\Xi\to\bbR$, 
		and some univariate convex polynomial $\psi_t\in\bbR[t]$ for each $t\in[T]$.
		Then for any $2k\ge\max\{p,\deg(F_x)\}$, $\calV_{p,2k}^{(i)}$
		for single-stage relaxation \eqref{eq:1stage} 
		is convex in $x$ and $\lambda$ for each $i\in[N]$.
	\end{proposition}
	\begin{proof}
		We claim that for each feasible tms $y\in\bar{\scrS}[h]_{2k}$, 
		the Hessian matrix 
		\[
		\nabla_x^2\pAngle{F_{t,x}-\lambda Q_p^{(i)}}{y} 
		= \pAngle{\nabla_x^2 F_x}{y} \succeq 0
		\]
		is positive semidefinite for any $x$.
		Then the convexity of $\calV_{p,2k}^{(i)}$ follows 
		since it is a maximum of convex functions.
		To see the claim, it suffices to show 
		$\pAngle{v^\transpose\nabla_x^2F_t(x,\xi)v}{y}\ge 0$ for any $v\in\bbR^{n_1}$. 
		By definition, 
		$\nabla_x^2 F_t(x,\xi) = \nabla_x^2\psi_t(a_t(\xi)^\transpose x+b_t(\xi) )$, so
		\[
		\pAngle{v^\transpose\nabla_x^2F_t(x,\xi)v}{y}
		=\pAngle{\psi_t''(a_t(\xi)^\transpose x+b_t(\xi)(a_t(\xi)^\transpose v)^2}{y}.
		\]
		Since $\psi_t(t)$ is a univariate convex polynomial, 
		its second derivative $\psi_t''(t)\ge0$ for all $t\in\bbR$.
		Consequently, $\psi_t''(t)$ is a sum of squares.
		Then $\psi_t''(a_t(\xi^\transpose x+b_t(\xi)))(a_t(\xi)^2v)\in\Sigma[\xi]_{\deg(F_{x})}$, 
		which by the definition of tms implies the desired nonnegativity.
	\end{proof}

\section{Moment Relaxations for Single-Stage Wasserstein DRO}
\label{sec:1stage}

In this section, we study the question regarding $\mc{O}(r)$-consistency 
of moment relaxations for single-stage data-driven $p$-WDRO~\eqref{eq:1stage}. 
	Recall that 
	\begin{equation}\label{eq:def_F_max}
		F(x,\xi) = \max\{ F_1(x,\xi), \ldots, F_T(x,\xi)\},
	\end{equation}
	where each $F_t$ is polynomial in $\xi$. 
	Denote $F_{t,x}(\cdot)\coloneqq F_t(x,\cdot)$ and 
	use $\deg(F_x)$ to denote the highest degree of 
	$F_{t,x}$ over $t\in[T]=\{1,\ldots, T\}$.
	For the special case that $T=1$, 
	$F(x,\xi)$ reduces to a single polynomial in $\xi$ and $F_x(\cdot)=F_{1,x}(\cdot)$.
Recall that $\Xi$ is determined by $h(\xi)\ge 0$ for $h = (h_1,\ldots, h_{m_1})$. 
Let $p_1 \coloneqq p/2$ and
\begin{equation}\label{eq:d1}
	d_1\coloneqq \max\{\lceil\deg(F_{x})/2\rceil,\,\lceil \deg(h)/2\rceil,\, 
	p_1\},
\end{equation}
where $\lceil a\rceil\coloneqq \min\{t\in \N: t\ge a\}$.
This $d_1$ is the lowest moment relaxation order that admits 
all information about $F$, $h$, and $Q_p^{(i)}$.
The value function $v_p(x;r^p)$ from~\eqref{eq:MeasConstr} can be written as
\begin{equation}\label{eq:ValueFunc}
	v_p(x;r^p) = \min_{\lambda\in \re_{\ge 0}} 
	\Big\{ \lambda r^p + \frac{1}{N}\sum\limits_{i=1}^N
	\max\limits_{t\in T} \max\limits_{\xi^{(i)}\in \Xi} 
	\big[ F_{t}(x,\xi^{(i)})-\lambda Q_p^{(i)}(\xi^{(i)}) 
	\big] \Big\}.
\end{equation}  
Its $k$-th order moment relaxation (with respect to the inner maximization) 
becomes
\begin{equation}\label{eq:vf_k}
	v_{p,k}(x;r^p) \coloneqq \min_{\lambda\in \re_{\ge 0}} 
	\Big\{ \lambda r^p + \frac{1}{N}\sum\limits_{i=1}^N
	\max\limits_{t\in[T]}
	\max\limits_{y^{(i)}\in \bar{\mathscr{S}}[h]_{2k}} 
	\langle F_{t,x}-\lambda Q_p^{(i)}(\xi),y^{(i)} \rangle \Big\}.
\end{equation}
	When $r=0$, $v_p(x;0)=v(x;0) = \mathbb{E}_{\xi\sim\hat{\nu}} F(x,\xi)$ by convention.
	For every $r>0$, it holds that $v_{p,k}(x;r^p)\ge v_{p}(x;r^p)\ge v_{p}(x;0)$. 
	This is because $\xi^{(i)} = \hat{\xi}^{(i)}$ is always feasible 
	to the inner maximization of \eqref{eq:ValueFunc}.

\subsection{Proof of asymptotic consistency}
First, we focus on bounding the difference $v_{p,k}(x;r^p)-v_p(x;r^p)=\calO(r)$
for single-stage $p$-WDRO with 
$F$ being a maximum of generic polynomials.
\begin{remark}\label{rem:1stage_opt}
	Such a difference will imply the $\calO(r)$-optimality gap 
	in~Theorem~\ref{thm:1stage-intro} since for any optimal solution 
	$\bar{x}\in\operatorname{argmin}\{f(x)+v_{p,k}(x;r^p):x\in X\}$, 
	we have
	\[
	f(\bar{x})+v_p(\bar{x};r^p)
	\le f(\bar{x})+v_{p,k}(\bar{x};r^p)
	\le \min_{x\in X}\{f(x)+v_p(x;r^p)\}+\calO(r).
	\]
\end{remark}
To obtain the $\mc{O}(r)$-optimality gap for $r>0$,
we begin with the following lemmata. 

\begin{lemma}\label{lem:uni_mom}
	If $y = (1,y_1,\ldots, y_{2k})\in\re^{\N_{2k}^1}$ satisfies $M_k[y]\succeq 0$, 
	then $(y_l)^{2k}\le (y_{2k})^l$ for every $l\in [2k]=\{1,\ldots, 2k\}$.
\end{lemma}
\begin{proof}
		Let $\mathscr{P}_{1,2k}$ denote the set of all nonnegative univariate polynomials 
		and $\mathscr{R}_{1,2k} \coloneqq \operatorname{cone}\{[z_i]_{2k}: z\in\re\}$. 
		As stated in \cite[Theorem 2.4.7]{NieBook}, 
		the closure of $\mathscr{R}_{1,2k}$ is the dual cone of $\mathscr{P}_{1,2k}$, i.e., $cl(\mathscr{R}_{1,2k}) = \mathscr{P}_{1,2k}^*$.
		Hilbert's theorem \cite[Theorem~3.4]{laurent2009sums} ensures that 
		$ \mathscr{P}_{1,2k} = \Sigma[u]_{2k} 
		= \{[u]_k^{\transpose}P[u]_k: P\succeq 0\}$ for $u\in\re^1$.
		Since the psd cone is self-dual, 
		a tms $y = (1, y_1,\ldots,y_{2k})$ satisfies $M_k[y]\succeq 0$ 
		if and only if $y\in cl(\mathscr{R}_{1,2k})$.
		If $y\in \mathscr{R}_{1,2k}$, then there exist $N\in \N$ and 
		$\theta_i\ge 0$, $u_i\in\re^1$ for each $i\in [N]$ such that
		\[
		y = \theta_1[u_1]_{2k}+\cdots+\theta_N[u_N]_{2k}.
		\]
		Note that $y_0 = \sum_{i=1}^N\theta_i = 1$.
		By Jensen's inequality, for all $l\in [2k]$, we have
		\[
		(y_l)^{2k/l} 
		\le \Big(\sum\limits_{i=1}^N \theta_i \big\vert u_i^l\big\vert\Big)^{2k/l}
		\le \sum\limits_{i=1}^N \theta_i \big\vert u_i^l\big\vert^{2k/l}
		= \sum\limits_{i=1}^N \theta_i u_i^{2k}
		= y_{2k}.
		\]
		If $y\not\in \mathscr{R}_{1,2k}$, then there exists a sequence 
		$\{y^{(j)}\}_{j=1}^{\infty}$ in $\mathscr{R}_{1,2k}$ such that 
		$y^{(j)}\to y$ as $j\to \infty$. 
		By previous arguments, $y_{2k}^{(j)}\ge (y_{2k}^{(j)})^{2k/l}$ for all $l\in [2k]$, 
		which implies
		$
		y_{2k} = \lim_{j\to \infty}y_{2k}^{(j)} 
		\ge \lim_{j\to \infty}(y_l^{(j)})^{2k/l} = (y_l)^{2k/l}.
		$
		So the conclusion holds.
\end{proof}

\begin{lemma}\label{lem:multi_mom}
	Suppose $y\in \re^{\N_{2k}^n}$ satisfies $y_0 = 1$ and $M_{k}[y]\succeq 0$.
	Then for every nonzero monomial vector 
	$\alpha = (\alpha_1,\ldots, \alpha_n)\in\N_{2k}^n$, 
	it holds that
	\begin{equation}\label{eq:multi_mom}
		|y_{\alpha}|^{2k}\le (y_{2ke_1})^{\alpha_1}\cdots (y_{2ke_n})^{\alpha_n},
	\end{equation}
	where $e_j\in\re^n$ denotes the $j$-th column vector of the $n\times n$ identity matrix. 
\end{lemma}
\begin{proof}
	Let $s(\alpha)$ denote the number of nonzero entries in $\alpha$.
	We prove this result by inducting on $s(\alpha)$.
	When $s(\alpha) = 1$, there exists some $j\in[n]$ such that 
	$\alpha = |\alpha|e_j$, then the result is implied by $y_{0} = 1$ and Lemma~\ref{lem:uni_mom}.
	Assume \eqref{eq:multi_mom} holds for all $\alpha\in\N_{2k}^n$ with $s(\alpha)\le s$.
	For all $\alpha\in \N_{2k}^n$ with $s(\alpha) = s+1$, 
	we can decompose it as $\alpha = \beta+\gamma$, 
	where $\beta, \gamma\in \N_k^n$ and $s(\beta), s(\gamma) \le s$.
	Without loss of generality, assume $\gamma\ge \beta$ with the degree-lexicographic order.
	Since $M_k[y]\succeq 0$, the principal submatrix
	\[
	\bbm y_{2\beta} & y_{\alpha}\\ y_{\alpha} & y_{2\gamma}\ebm \succeq 0
	\,\Rightarrow\, \left\vert\begin{matrix}
		y_{2\beta} & y_{\alpha}\\ y_{\alpha} & y_{2\gamma}
	\end{matrix}\right\vert = y_{2\beta}y_{2\gamma}-y_{\alpha}^2\ge 0.
	\]
	Then $|y_{\alpha}|^{2k}\le (y_{2\beta}y_{2\gamma})^k
	\le \prod_{l=1}^n (y_{2ke_l})^{\beta_l}(y_{2ke_l})^{\gamma_l}
	= (y_{2ke_1})^{\alpha_1}\cdots (y_{2ke_n})^{\alpha_n}.$
\end{proof}

\begin{lemma}\label{lem:norm}
	Assume $\|\xi\| = (\xi^{\transpose}H\xi)^{1/2}$ 
	for some positive definite matrix $H\succ 0$.
	For all $2k\ge p$, if $y\in\bar{\mathscr{S}}[h]_{2k}$ satisfies 
	$\langle \|\xi\|^p, y\rangle \le \rho\in\re_{\ge 0}$, 
	then 
	\[
	(\theta_{\min})^{p_1}\langle \|\xi\|_p^p, y\rangle = 
	(\theta_{\min})^{p_1}\sum\limits_{l=1}^{n_0} y_{pe_l}\le \rho,
	\]
	where $\theta_{\min}$ is the smallest eigenvalue of $H$.
\end{lemma}
\begin{proof}
	Without loss of generality, we can assume $H = H^{\transpose}$.
	By singular value decomposition, we can write $H = Q\Theta Q^{\transpose}$,
	where $Q$ is an orthogonal matrix and $\Theta = \operatorname{diag}(\theta)\succ 0$ 
	with the vector of eigenvalues $\theta = (\theta_1, \ldots, \theta_{n_0})$ for $H$.
	Then $\|\xi\|^2 = (Q\xi)^{\transpose} \Theta(Q\xi)$.
	Let $ \theta_{\min} $ denote the smallest eigenvalue of $ H $.
		For each $y\in\bar{\mathscr{S}}[h]_{2k}$, it holds that
		\[
		\begin{aligned}
			\langle \|\xi\|^p, y\rangle 
			&= \langle (\xi^{\transpose}Q^{\transpose}\Theta Q\xi)^{p_1}, y\rangle
			= \sum\limits_{\alpha\in\N^{n_0}}^{|\alpha| = p_1}\frac{p_1!}{\alpha!}\cdot \theta^{\alpha}
			\langle (Q^{\transpose}\xi)^{2\alpha}, y\rangle \\
			&= \sum\limits_{\alpha\in \N^{n_0}}^{|\alpha| = p_1}\frac{p_1!}{\alpha!}\cdot \theta^{\alpha}\cdot
			\mbox{vec}[(Q^{\transpose}\xi)^{\alpha}]^{\transpose} M_k[y]\mbox{vec}[(Q^{\transpose}\xi)^{\alpha}],
		\end{aligned}
		\]
		where $\alpha!\coloneqq \alpha_1!\cdots \alpha_{n_0}!$ for 
		each $\alpha = (\alpha_1,\ldots, \alpha_{n_0})\in \N^{n_0}$.
	Since $M_k[y]$ is psd and $p_1!\ge \alpha!$ for all $\alpha\in\N^{n_0}$ with $|\alpha|=p_1$, 
	we further have
	\[  \begin{aligned}
		\langle \|\xi\|^p, y\rangle 
		& \ge (\theta_{\min})^{p_1}\langle (\xi^{\transpose}Q^{\transpose}Q\xi)^{p_1}, y\rangle
		= (\theta_{\min})^{p_1}\langle (\|\xi\|_2^{p}, y\rangle\\
		& \ge (\theta_{\min})^{p_1}\langle \|\xi\|_p^p, y\rangle = (\theta_{\min})^{p_1}\sum\limits_{l=1}^{n_0} y_{pe_l}.
	\end{aligned}\]
	By this inequality, if $\langle \|\xi\|^p, y\rangle \le \rho$, then 
	we must have $(\theta_{\min})^{p_1}\langle \|\xi\|_p^p, y\rangle\le \rho$
\end{proof}
	\begin{remark}
		In computational practice, $H$ is often selected as the identity matrix, 
		which makes $\theta_{\min} = 1$ and reduces the norm $\|\xi\|$ 
		to the $l_2$-norm $\|\xi\|_2$. 
		This leads to a simplified corollary of Lemma~\ref{lem:norm}:
		for all $2k\ge p$, if $y\in \bar{\mathscr{S}}[h]_{2k}$ satisfies 
		$\langle \|\xi\|_2^p, y\rangle \le \rho$, 
		then $\langle \|\xi\|_p^p, y\rangle\le \rho$.
	\end{remark}
	
	For convenience, we set $ \theta_{\min} = 1 $ 
	in the following context without loss of generality.
We are now ready to prove the following more elaborate version of~Theorem~\ref{thm:1stage-intro}.

\begin{theorem}\label{thm:asym1str0}
	Under Assumption~\ref{assum}, 
	suppose $X$ is bounded and the relaxation order $k\ge d_1$. 
	For all $r>0$ that is sufficiently small, 
	the difference $v_{p,k}(x;r^p)-v_p(x;r^p) = \calO(r)$
	if either
	\begin{enumerate}
		\item[(i)] $p\ge \deg(F_x)$; or
		\item[(ii)] there exists $R>0$, $k_1\in \N$ such that 
		$R-\|\xi\|_2^2\in \qmod{h}_{2k_1+2}$, 
		and $\deg(F_x)\le 2k-2k_1+p$.       
	\end{enumerate}
	In either case, any optimal solution $x$ to the moment relaxation~\eqref{eq:1stage} 
	is an $\calO(r)$-optimal solution to the $p$-WDRO~\eqref{eq:WDRO} 
	by~Remark~\ref{rem:1stage_opt}.
\end{theorem}
\begin{proof}
	
		Since $v_{p}(x;r^p)\ge v_{p}(x,0) = \frac{1}{N}\sum_{i=1}^N F(x,\hat{\xi}^{(i)})$ 
		for all $r>0$.
		It suffices to show $v_{p,k}(x;r^p)-v_p(x;0)=\mc{O}(r)$ for all $r>0$ 
		that is sufficiently small.
		By Taylor's expansion, $v_{p,k}(x;r^p)-v_p(x;0) $ equals the optimal value of
		\begin{equation}\label{eq:Taylor_ref}
			\begin{aligned}
				\min_{\lambda\in \re_{\ge 0}}  
				\lambda r^p + \frac{1}{N}\sum\limits_{i=1}^N
				\max\limits_{t\in T} \max\limits_{y^{(i)}\in \bar{\mathscr{S}}[h]_{2k}} 
				& \Big\{\big\langle \sum\limits_{d=1}^{\deg(F_x)} 
				\sum\limits_{\alpha\in\N^{n_0}}^{|\alpha| = d}
				\frac{\partial^{\alpha} F_{t,x}(\hat{\xi}^{(i)})}{\alpha!}
				(\xi-\hat{\xi}^{(i)})^{\alpha}, y^{(i)}
				\big\rangle\\
				& \qquad  -\langle\lambda Q_p^{(i)},y^{(i)}\rangle\Big\}
			\end{aligned}
		\end{equation}
		
		(i) Suppose $\deg(F_x)\le p$.
		We first construct an upper bound of \eqref{eq:Taylor_ref}.
		Let $\hat{\xi}\in \Xi$ and $\alpha = (\alpha_1,\ldots, \alpha_{n_0})\in\N^{n_0}$ 
		with $|\alpha|=d\le p$.
		We bound $|\langle (\xi-\hat{\xi})^{\alpha}, y\rangle|$ according to the value of $d$.
		If $d=1$, then $\alpha=e_l$ for some $l\in[n_0]$,
		\[
		|\langle \xi_l-\hat{\xi}_l, y\rangle|\le \max\limits_{z_l=(z_{l,1},z_{l,2})}
		\left\{
		z_{l,1}: z_{l,1}\ge 0,\, 
		\bbm 1 & z_{l,1}\\z_{l,1} & z_{l,2}\ebm \succeq 0,\,
		z_{l,2} = \langle (\xi_l-\hat{\xi}_l)^2,y\rangle
		\right\}.
		\]
		If $d>1$ is even, then Lemma~\ref{lem:multi_mom} and the AM-GM inequality imply
		\begin{equation}\label{eq:tms_bd}
			0\le\langle (\xi-\hat{\xi})^{\alpha}, y\rangle
			\le \prod\limits_{l=1}^{n_0}[\langle(\xi_l-\hat{\xi}_l)^{d}, y\rangle]^{\frac{\alpha_l}{d}}\le 
			\sum\limits_{l=1}^{n_0}\frac{\alpha_l}{d}\langle(\xi_l-\hat{\xi}_l)^{d}, y\rangle.
		\end{equation}
		Since $d\le p$ and $\binom{n_0+d-1}{d}$ counts the number of monomials 
		$\alpha\in \N^{n_0}$ with $|\alpha|=d$,
		\[
		\sum\limits_{\alpha\in\N^{n_0}}^{|\alpha|=d} \langle (\xi-\hat{\xi})^{\alpha}, y\rangle
		\le \sum\limits_{\alpha\in\N^{n_0}}^{|\alpha|=d}
		\sum\limits_{l=1}^{n_0}\frac{\alpha_l}{d}\langle(\xi_l-\hat{\xi}_l)^{d}, y\rangle
		\le \frac{\binom{n_0+p-1}{p}}{n_0}\langle\|\xi-\hat{\xi}\|_d^d, y\rangle.
		\]    
		If $d>1$ is odd, then there exists some $\beta,\eta\in\N^{n_0}$ 
		with $2|\beta|=d+1$ and $2|\eta|=d-1$ such that $\alpha = \beta+\eta$.
		By the psd property of $M_k[y]$,
		\[
		|\langle (\xi-\hat{\xi})^{\alpha}, y\rangle|
		\le \sqrt{\langle (\xi-\hat{\xi})^{2\beta}, y\rangle}\cdot\sqrt{\langle (\xi-\hat{\xi})^{2\eta}, y\rangle}
		\le \frac{1}{2}\big\langle (\xi-\hat{\xi})^{2\beta}+ (\xi-\hat{\xi})^{2\eta}, y\big\rangle.
		\]
		Based on the bound \eqref{eq:tms_bd} for the even case, we obtain
		\[\begin{aligned}
			\sum\limits_{\alpha\in\N^{n_0}}^{|\alpha|=d}|\langle (\xi-\hat{\xi})^{\alpha}, y\rangle|
			\le \frac{\binom{n_0+p-1}{p}}{2n_0}\sum\limits_{l=1}^{n_0}
			\big\langle\|\xi-\hat{\xi}\|_{d+1}^{d+1}+\|\xi-\hat{\xi}\|_{d-1}^{d-1}, y\big\rangle.
		\end{aligned}\]
		Since $X$ is bounded, there exists a constant $\mc{F}^*>0$ such that
		\[
		\mc{F}^*(1+\|\hat{\xi}\|_p)^{p-|\alpha|}\ge \frac{\binom{n_0+p-1}{p}}{n_0}\left\vert
		\frac{\partial^{\alpha} F_{t,x}(\hat{\xi})}{\alpha!}
		\right\vert
		\]
		for all $t\in[T]$, $\hat{\xi}\in\re^{n_0}$ and $0<\alpha\in \N_p^{n_0}$. 
		Then \eqref{eq:Taylor_ref} is bounded from above by 
		\[
		\min\limits_{\lambda\in \re_{\ge 0}} 
		\left\{ 
		\begin{aligned}
			\lambda r^p +\frac{1}{N}\sum\limits_{i=1}^N 
			\max\limits_{y^{(i)}, z^{(i)}} &  
			\Big\{-\lambda\langle Q_p^{(i)},y^{(i)}\rangle +\mc{F}^*(1+\|\hat{\xi}^{(i)}\|_p)^{p-1}\Big(\sum\limits_{l=1}^{n_0} z_{l,1}^{(i)}\Big)\\
			& +2\mc{F}^*\sum\limits_{d=1}^{p/2}(1+\|\hat{\xi}^{(i)}\|_p)^{p-2d}\langle\|\xi-\hat{\xi}^{(i)}\|_{2d}^{2d}, y^{(i)}\rangle\Big\}\\
			\st &\quad  y^{(i)}\in \hat{\mathscr{S}}[h]_{2k},\, z^{(i)}=(z^{(i)}_{l,1}, z^{(i)}_{l,2}),\, l\in[n_0],\\
			& z_{l,1}^{(i)}\ge 0,\,
			\bbm 1 & z_{l,1}^{(i)}\\
			z_{l,1}^{(i)} & z_{l,2}^{(i)}\ebm\succeq 0,\,
			z_{l,2}^{(i)} = \langle (\xi_l-\hat{\xi}_l^{(i)})^2, y^{(i)}\rangle
		\end{aligned}\right\}.
		\]
		Its Lagrangian dual problem is the semidefinite program
		\begin{equation}\label{eq:def_phi}
			\begin{aligned}
				\phi(x;r^p) \coloneqq \max\limits_{y^{(i)}, z^{(i)}}\quad  & 
				\frac{\mc{F}^*}{N}\sum\limits_{i=1}^N
				\Big\{(1+\|\hat{\xi}^{(i)}\|_p)^{p-1}
				\Big(\sum\limits_{l=1}^{n_0} z_{l,1}^{(i)}\Big) \\
				&\qquad \qquad  +2\sum\limits_{d=1}^{p/2}(1+\|\hat{\xi}^{(i)}\|_p)^{p-2d}\langle\|\xi-\hat{\xi}^{(i)}\|_{2d}^{2d}, y^{(i)}\rangle\Big\}\\
				\st\quad  &  y^{(i)}\in \hat{\mathscr{S}}[h]_{2k},\, z^{(i)}=(z^{(i)}_{l,1}, z^{(i)}_{l,2}),\,l\in[n_0],\,\\
				& z_{l,1}^{(i)}\ge 0,\,
				\bbm 1 & z_{l,1}^{(i)}\\
				z_{l,1}^{(i)} & z_{l,2}^{(i)}\ebm\succeq 0,\,
				z_{l,2}^{(i)} = \langle (\xi_l-\hat{\xi}_l^{(i)})^2, y^{(i)}\rangle,\\
				& \frac{1}{N}\sum\limits_{i=1}^N \langle Q_p^{(i)}, y^{(i)}\rangle \le r^p.
			\end{aligned}
		\end{equation}
		Since $\Xi$ is regular closed, 
		problem \eqref{eq:def_phi} is strictly feasible for all $r>0$.
		Thus, the strong duality holds by the Slater's condition 
		and $v_{p,k}(x;r^p)-v_p(x;0)\le \phi(x;r^p)$.
		
		We then show $\phi(x;r^p)=\mc{O}(r)$ when $r>0$ is sufficiently small.
		Suppose $(y^{(i)}, z^{(i)})_{i\in[N]}$ is an optimizer of \eqref{eq:def_phi}. 
		For convenience, denote 
		\begin{equation}\label{eq:s1sd}
			s_1^{(i)} = \sum\limits_{l=1}^{n_0} z_{l,1}^{(i)},\quad
			s_{2d}^{(i)} = \langle\|\xi-\hat{\xi}^{(i)}\|_{2d}^{2d}, y^{(i)}\rangle,\, \forall d\in[p_1].
		\end{equation}
		Then $s_1^{(i)}\le \sum_{l=1}^{n_0}
		\langle (\xi_l-\hat{\xi}_l^{(i)})^2, y^{(i)}\rangle^{\frac{1}{2}}$.
		By Lemma~\ref{lem:multi_mom} and H\"{o}lder's inequality,
		\begin{align}
			\label{eq:s1bd}
			s_1^{(i)}& 
			\le \sum\limits_{l=1}^{n_0} 
			\langle (\xi_l-\hat{\xi}_l^{(i)})^p, y^{(i)}\rangle^{\frac{1}{p}}
			\le n_0^{\frac{p-1}{p}}
			\langle \|\xi-\hat{\xi}^{(i)}\|_p^p,y^{(i)}\rangle^{\frac{1}{p}},\\
			\label{eq:sdbd}
			s_{2d}^{(i)}& \le \sum\limits_{l=1}^{n_0}
			\langle (\xi_l-\hat{\xi}_l^{(i)})^p, y^{(i)}\rangle^{\frac{2d}{p}}
			\le n_0^{\frac{p-2d}{p}} \langle \|\xi-\hat{\xi}^{(i)}\|_p^p,y^{(i)}\rangle^{\frac{2d}{p}},
			\quad \forall d\in[p_1].
		\end{align}
		By H\"{o}lder's inequality (applied to the sum on $i$) and
		Assumption~\ref{assum}, we obtain
		\[\begin{aligned}
			\phi(x;r^p) & 
			\le \mc{F}^*\Big(\frac{1}{N} \sum\limits_{i=1}^N 
			(1+\|\hat{\xi}^{(i)}\|_p)^p\Big)^{\frac{p-1}{p}}
			\Big(\frac{1}{N}\sum\limits_{i=1}^N (s_1^{(i)})^p\Big)^{\frac{1}{p}}\\
			&\quad +2\mc{F}^*\sum\limits_{d=1}^{p_1}
			\Big(\frac{1}{N} \sum\limits_{i=1}^N (1+\|\hat{\xi}^{(i)}\|_p)^p\Big)^{\frac{p-2d}{p}}
			\Big(\frac{1}{N}\sum\limits_{i=1}^N (s_{2d}^{(i)})^{\frac{p}{2d}}\Big)^{\frac{2d}{p}}\\
			& \le \mc{F}^*\mc{B}^{\frac{p-1}{p}}
			\Big(\frac{1}{N}\sum\limits_{i=1}^N 
			n_0^{p-1}\langle \|\xi-\hat{\xi}^{(i)}\|_p^p, y^{(i)}\rangle\Big)^{\frac{1}{p}}\\
			& \quad
			+2\mc{F}^*\sum\limits_{d=1}^{p_1}\mc{B}^{\frac{p-2d}{p}}
			\Big(\frac{1}{N}\sum\limits_{i=1}^N n_0^{\frac{p-2d}{2d}}
			\langle \|\xi-\hat{\xi}^{(i)}\|_p^p, y^{(i)}\Big)^{\frac{2d}{p}}.
		\end{aligned}
		\]
		By Lemma~\ref{lem:norm} (with $\theta_{\min}=1$),
		the constraint $\sum_{i=1}^N\langle Q_p^{(i)}, y^{(i)}\rangle \le Nr^p$ 
		implies $\sum_{i=1}^N\langle\|\xi-\hat{\xi}^{(i)}\|_p^p, y^{(i)}\rangle\le Nr^p$.
		Then 
		\begin{equation}\label{eq:estimateVcasei}
			\phi(x;r^p)\le \mc{F}^*(\sqrt[p]{\mc{B}n_0})^{p-1}r+2\mc{F}^*\sum\limits_{d=1}^{p_1}(\sqrt[p]{\mc{B}n_0})^{p-2d}r^{2d},
		\end{equation}
		which is dominated by the linear term when $r>0$ is sufficiently small.

		(ii) Suppose  $R-\|\xi\|_2^2\in \qmod{h}_{2k_1+2}$ for 
		some $R>0, k_1\in\N$ and $\deg(F_x)\le 2k-2k_1+p$.
		Then $\Xi$ is compact and $\mc{B} = (1+\sqrt{R})^p$ is the constant 
		that satisfies Assumption~\ref{assum}.
		Similar to (i), 
		there exists a constant $\hat{\mc{F}}^*>0$ such that 
		$v_{p,k}(x;r^p)-v_p(x;0)\le \hat{\phi}(x;r^p)$  for
		\begin{equation}\label{eq:def_phi_hat}
			\begin{aligned}
				\hat{\phi}(x;r^p) \coloneqq \max\limits_{y^{(i)}, z^{(i)}}
				\quad  & \frac{\hat{\mc{F}}^*}{N}\sum\limits_{i=1}^N
				\Big\{\Big(\sum\limits_{l=1}^{n_0} 
				z_{l,1}^{(i)}\Big) +2\sum\limits_{d=1}^{\lceil \deg(F_x)/2\rceil}
				\langle\|\xi-\hat{\xi}^{(i)}\|_{2d}^{2d}, y^{(i)}\rangle\Big\}\\
				\st\quad  &  y^{(i)}\in \hat{\mathscr{S}}[h]_{2k},\, z^{(i)}=(z^{(i)}_{l,1}, z^{(i)}_{l,2}),\,l\in[n_0],\,\\
				& z_{l,1}^{(i)}\ge 0,\,
				\bbm 1 & z_{l,1}^{(i)}\\
				z_{l,1}^{(i)} & z_{l,2}^{(i)}\ebm\succeq 0,\,
				z_{l,2}^{(i)} = \langle (\xi_l-\hat{\xi}_l^{(i)})^2, y^{(i)}\rangle,\\
				& \frac{1}{N}\sum\limits_{i=1}^N \langle Q_p^{(i)}, y^{(i)}\rangle \le r^p.
		\end{aligned}\end{equation}
		Suppose $(y^{(i)}, z^{(i)})_{i\in[N]}$ is an optimizer of \eqref{eq:def_phi}. 
		Let $s_1^{(i)}$, $s_{2d}^{(i)}$ be defined in \eqref{eq:s1sd} for all $d\in[t]$.
		With the conclusion of (i), it suffices to show 
		$\frac{1}{N}\sum_{i=1}^N s_{2d}^{(i)}\le\mc{O}(r)$ 
		for $d=p_1+1,\ldots, \lceil \frac{\deg(F_x)}{2}\rceil$ 
		when $r>0$ is sufficiently small.
		Since $y^{(i)}\in\hat{\mathscr{S}}[h]_{2k}= \mathscr{S}[h]_{2k}\cap \{y_0=1\}$,
		for all $k\ge k_1+1$,
		the dual relation $\mathscr{S}[h]_{2k} = \qmod{h}_{2k}^*\subseteq \qmod{h}_{2k_1}^*$
		\cite[Theorem~2.5.2]{NieBook} and Lemma~\ref{lem:norm} (with $\theta_{\min} = 1$) ensure
		\[
		\langle R-\|\xi\|_2^2, y^{(i)}\rangle\ge
		\langle R-\|\xi\|^2, y^{(i)}\rangle \ge 0.
		\]
		For each $d \in[k-k_1]$, since $\|\xi\|_2^{2d-2}(R-\|\xi\|^2)\in\qmod{h}_{2k}$, 
		we also have
		\[
		R\langle \|\xi\|^{2d-2}, y^{(i)}\rangle \ge \langle \|\xi\|^{2d}, y^{(i)}\rangle
		\quad\Rightarrow\quad 
		\langle \|\xi\|^{2d}, y^{(i)}\rangle \le R^d,\,\forall d\in[k-k_1].
		\]
		Noting that each $\hat{\xi}^{(i)}\in \Xi$ ensures 
		$\|\hat{\xi}^{(i)}\|^2\le R$, the above result yields
		\[
		\langle (\xi_l-\hat{\xi}_l^{(i)})^{2d},y^{(i)}\rangle\le 4^{d}R^d,\quad 
		\forall d\in[k-k_1].
		\]
		Since $\lceil\frac{\deg(F_x)}{2}\rceil\le k-k_1+p_1$, 
		for each $d=p_1+1,\ldots, \lceil\frac{\deg(F_x)}{2}\rceil$, 
		\[\begin{aligned}
			\langle (\xi_l-\hat{\xi}_l^{(i)})^{2d}, y^{(i)}\rangle& \le \sqrt{\langle (\xi_l-\hat{\xi}_l^{(i)})^{2d-p},y^{(i)}\rangle}
			\sqrt{\langle (\xi_l-\hat{\xi}_l^{(i)})^{p},y^{(i)}\rangle}\\
			& \le (4R)^{\lceil\frac{\deg(F_x)}{2}\rceil-p_1}
			\sqrt{\langle (\xi_l-\hat{\xi}_l^{(i)})^p, y^{(i)}\rangle}.
		\end{aligned}\]
		Therefore, we get
		\[\begin{aligned}
			\frac{1}{N}\sum\limits_{i=1}^N
			s_{2d}^{(i)} & \le 
			\frac{(4R)^{\lceil\frac{\deg(F_x)}{2}\rceil-p_1}\sqrt{n_0}}{N}\sum\limits_{i=1}^N
			\sqrt{\langle \|\xi-\hat{\xi}^{(i)}\|_p^p,y^{(i)}\rangle}\\
			& \le (4R)^{\lceil\frac{\deg(F_x)}{2}\rceil-p_1}\sqrt{n_0}r^{p_1}.
		\end{aligned}\]
		Since $p_1=p/2\ge 1$, $\frac{1}{N}\sum_{i=1}^N s_{2d}^{(i)}\le \mc{O}(r)$ 
		when $r>0$ is sufficiently small.
		The conclusion is thus implied. In particular,
		if $p=2$, then with analysis in part (i),
		\[
		\hat{\phi}(x;r^p)\le \Big(\hat{\mc{F}}^*(1+\sqrt{R})\sqrt{n_0}+
		\mc{F}^*(4R)^{\lceil\frac{\deg(F_x)}{2}\rceil-1}\sqrt{n_0}\Big)r+o(r).
		\]
		If $p>2$, then 
		$\hat{\phi}(x;r^p)\le \hat{\mc{F}}^*(1+\sqrt{R})^{p-1}(n_0)^{\frac{p-1}{p}}r+o(r).$ 
\end{proof}

The following example shows that when $\Xi$ is unbounded, 
generally we need $p\ge\deg(F_x)$ for the moment relaxations 
to preserve the $\mc{O}(r)$-consistency.
\begin{example}\label{ex:1stage}
	Consider $X=\{0\}$, $f(x)=0$, $\Xi=\{\xi = (\xi_1,\xi_2,\xi_3)\in\bbR^3:\xi\ge0\}$ 
	is the nonnegative orthant, with $N=1$, $\hat{\xi}^{(1)}=(1,1,1)$, 
	and $F_x(\xi)=3\xi_1\xi_2\xi_3-\xi_1^2\xi_2-\xi_1\xi_2^2-\xi_3^3$.
	By the inequality of arithmetic and geometric means, $F_x(\xi)\le0$ on $\Xi$, 
	and thus $\Gamma_p(x)=0$ for any $p\ge1$.
	Thus we can use $p=2$ despite that $\deg(F_x)=4$ for 
	the data-driven $p$-WDRO~\eqref{eq:WDRO}. 
	For any $r>0$, we can rewrite the problem as 
	\begin{equation*}
		\min_{\lambda\in\re_{\ge0}}\Big\{\lambda r^2+
		\sup_{\xi\in\Xi}\Big[
		3\xi_1\xi_2\xi_3-\xi_1^2\xi_2-\xi_1\xi_2^2-\xi_3^2
		-\lambda\sum_{j=1}^{3}(\xi_j-1)^2\Big]\Big\}.
	\end{equation*}
	We claim that the inner supremum here is $+\infty$ for any 
	$\lambda\ge0$ and $2k\ge4$, 
	which implies that moment relaxation~\eqref{eq:1stage} is unbounded 
	and thus does not have $\calO(r)$-consistency.
	Assume for contradiction that the moment relaxation is bounded from above by $v\in\bbR$.
	Since $\Xi$ has a nonempty interior,
	$v-F_x+\lambda Q_p^{(1)}\in\qmod{h}_{2k}$ for $h=(\xi_1,\xi_2,\xi_3)$.
	Namely, there exist $\sigma_0\in\Sigma[\xi]_{2k}$, $\sigma_1,\sigma_2,\sigma_3\in\Sigma[\xi]_{2k-2}$ such that
	\begin{equation*}
		v+\xi_1^2\xi_2+\xi_1\xi_2^2+\xi_3^3-3\xi_1\xi_2\xi_3+\lambda\sum_{j=1}^{3}(\xi_j-1)^2=\sigma_0+\sigma_1\xi_1+\sigma_2\xi_2+\sigma_3\xi_3.
	\end{equation*}
	First we observe that $\deg(\sigma_i)\le2$ for $i=0,1,\ldots,3$ because otherwise either of the following situations gives a contradiction.
	\begin{itemize}
		\item[(i)] If $\deg(\sigma_0)=2k\ge4$, then its homogeneous part of highest degree of the right-hand side is given by $\sigma_0$ and thus must be nonzero.
		\item[(ii)] If $\deg(\sigma_0)<2k$ and $\max\{\deg(\sigma_1),\deg(\sigma_2),\deg(\sigma_3)\}=2k-2\ge4$, then the homogeneous part of highest degree of the right-hand side is given by $\sigma_1\xi_1+\sigma_2\xi_2+\sigma_3\xi_3$.
		By taking the homogeneous part $\bar{\sigma}_1,\bar{\sigma}_2,\bar{\sigma}_3$ of degree $k-1$ of $\sigma_1,\sigma_2,\sigma_3$, respectively, we see that $\bar{\sigma}_1\xi_1+\bar{\sigma}_2\xi_2+\bar{\sigma}_3\xi_3=0$.
		Since $\bar{\sigma}_1,\bar{\sigma}_2, \bar{\sigma}_3$ are also sums of squares, evaluation on the Zariski open subset $\bbR_{>0}^3$ shows that $\bar{\sigma}_1=\bar{\sigma}_2=\bar{\sigma}_3=0$ in $\bbR[\xi]$.
	\end{itemize}
	We then take the degree-3 part of both sides, which gives
	\begin{equation*}
		\xi_1^2\xi_2+\xi_1\xi_2^2+\xi_3^3-3\xi_1\xi_2\xi_3=\bar{\sigma}_1\xi_1+\bar{\sigma}_2\xi_2+\bar{\sigma}_3\xi_3, 
	\end{equation*}
	where each of $\bar{\sigma}_1,\bar{\sigma}_2,\bar{\sigma}_3$ 
	is a homogeneous sum of squares (possibly 0) of degree 2.
	This is not possible because by substituting $\xi_j=\zeta_j^2$
	for some indeterminate $\zeta_j$, $j=1,2,3$, 
	the left-hand side becomes the Motzkin polynomial~\cite{motzkin1967arithmetic}, 
	which cannot be written as a sum of squares, such as the right-hand side, 
	in $\bbR[\zeta_1,\zeta_2,\zeta_3]$.
\end{example}

\subsection{Refined estimates of optimality gaps}\label{ssc:tighter_gap}
	Consider the special case that $F$ is a single polynomial in $\xi$, i.e., 
	$F(x,\xi) = F_1(x,\xi)$.
	In computational practice, the moment relaxations \eqref{eq:vf_k} 
	may attain the same optimal value as the $p$-WDRO problem \eqref{eq:ValueFunc} 
	for some finite $k$.
	While the specific relaxation order required for such exactness 
	is generally unknown a priori, 
	there exist certain classes of problems whose moment relaxations 
	are exact at the initial order. 
	We summarize below some sufficient conditions for this initial-order exactness.
	\begin{proposition}
		For all $k\ge d_1$, the moment relaxation \eqref{eq:vf_k}
		is exact if one of the following conditions holds.
		\begin{enumerate}
			\item[(i)]  $\|\cdot\|=\|\cdot\|_p$ and $F(x,\xi)$ is 
			a separate polynomial in $\xi$ \footnote{
				A polynomial $p\in\re[\xi]$ is said to be separable 
				if it can be decomposed as a sum of univariate polynomials, 
				i.e., $p(\xi) = p_1(\xi_1)+\cdots+p_{n_0}(\xi_{n_0})$} 
			and $\Xi = \mc{I}_1\times\cdots\times \mc{I}_{n_0}$ 
			is a Cartesian product of intervals (possibly unbounded)
			determined by linear or quadratic polynomials.
			\item[(ii)] $F(x,\xi)-\lambda\|\xi-\hat{\xi}^{(i)}\|^p$ 
			is a quadratic polynomial in $\xi$ and 
			$h$ is a single quadratic polynomial with some $u$ satisfying $h(u)>0$.
		\end{enumerate}
	\end{proposition}
	\begin{proof}
		(i) Fix $x\in X$, $\lambda\in\re_{\ge 0}$ and $i\in[N]$. 
		Under given conditions, there exist univariate scalar polynomials $q_l\in\re[\xi_l]$ and vector-valued polynomials $\bar{h}_l\in\re[\xi_l]$
		for $l\in [n_0]$ such that 
		$F(x,\xi)-\lambda\|\xi-\hat{\xi}^{(i)}\|^p = \sum_{l=1}^{n_0} q_l(\xi_l)$
		and $h = (\bar{h}_1,\ldots, \bar{h}_{n_0})$, $\mc{I}_l = \{\xi_l\in\re^1: \bar{h}_l(\xi_l)\ge 0\}$.
		Then the $i$-th inner maximization problem in \eqref{eq:WDRO} is equivalent to
		a summation of univariate polynomial optimization problems over intervals. 
		We can solve it via moment relaxations for each individual variate, i.e., 
		\begin{equation}\label{eq:ind_mom_rel}
			\sum\limits_{l=1}^{n_0} \max\limits_{\xi_l\in\mc{I}_l} q_l(\xi_l)
			\,\xrightarrow[\text{relaxation}]{\text{moment}}\,
			\sum\limits_{l=1}^{n_0} \max\limits_{y\in \bar{\mathscr{S}}[h_l]_{2k}^{u}} \langle q_l, y\rangle,
		\end{equation}
		where each $\bar{\mathscr{S}}[h_l]_{2k}^u$ denotes the tms set associated with the
		polynomial $h_l$ and the univariate variable $\xi_l$.
		When $\mathcal{I}_l$ is an interval and $\bar{h}_l$ is a linear 
		or quadratic polynomial vector, 
		the set of polynomials nonnegative on $\mathcal{I}_l$ is 
		entirely characterized by the quadratic module generated by $\bar{h}_l$,
		i.e., $\mathscr{P}_{2k}(\mc{I}_l) = \qmod{\bar{h}_l}_{2k}$;
		see \cite[Proposition~3.1.3, Theorem~3.2.1 and Theorem~3.2.3]{NieBook}.
		Since $\mathscr{S}[\bar{h}_l]_{2k} = \qmod{\bar{h}_l}_{2k}^* = \mathscr{P}_{2k}(\mc{I}_l)^*$,
		the moment relaxations for all univariate subproblems are exact 
		if $k\ge d_1$.
		Since $h = (\bar{h}_1,\ldots, \bar{h}_{n_0})$, 
		for each $\bar{y}\in\bar{\mathscr{S}}[h]$ and $l\in[m]$,
		there always exists $y\in \bar{\mathscr{S}}[h_l]_{2k}^u$ 
		such that $y$ is a truncation of $\bar{y}$.
		This implies \eqref{eq:1stage} is a tighter moment relaxation than 
		that in \eqref{eq:ind_mom_rel}, thus exact.

		(ii) This is implied by the well-known S-Lemma \cite[Theorem~5.2.6]{NieBook}:
		for any $\gamma_1\in \re$, if the quadratic polynomial 
		$F(x,\xi)-\lambda\|\xi-\hat{\xi}^{(i)}\|^p-\gamma_1\ge 0$ 
		over $\Xi = \{\xi: h(\xi)\ge 0\}$,
		then there exists a scalar $\gamma_2\ge 0$ such that 
		the polynomial
		$F(x,\xi)-\lambda\|\xi-\hat{\xi}^{(i)}\|^p-\gamma_1-\gamma_2h(\xi)$ is 
		a sum-of-squares.   
	\end{proof}
	
	Another interesting type of single-stage $p$-WDRO problems consists of those whose initial
	order moment relaxations exhibit an optimality gap strictly smaller than $\mc{O}(r)$.
	We characterize a subclass of these problems with a SOS-convexity assumption.
	A polynomial $q\in\re[\xi]$ is {\it SOS-convex} if its Hessian can be decomposed as
	$\nabla^2 q(\xi) = H(\xi)H(\xi)^T$ 
	for some matrix polynomials $H(\xi)$ \cite[Definition~7.1.1]{NieBook}. 
	While not all convex polynomials are SOS-convex \cite{AAAPP2013}, 
	several widely used ones such as linear and convex quadratic functions are SOS-convex. 
	In particular, $\|\xi\|^p = (\xi H\xi)^{p/2}$ is SOS-convex for any 
	positive even integer $p$. This can be verified by direct computation:
	if $p=2$, $\nabla^p \|\xi\|^2 = 2H\succ 0$; otherwise, 
	\[ 
	\nabla^2 \|\xi\|^p = p \, (\xi^\transpose H \xi)^{\frac{p}{2} - 2} 
	\left[ (p - 2) H \xi \xi^\transpose H + (\xi^\transpose H \xi) H \right]. 
	\]
	To prove the result, we need to use the Lagrangian dual of problem \eqref{eq:ValueFunc}:
	\begin{equation}\label{eq:UnderApprox}
		V_p(x;r^p) \coloneqq \max_{\xi^{(i)}\in \Xi, i\in[N]} 
		\Big\{
		\frac{1}{N}\sum_{i=1}^{N} F(x,\xi^{(i)}):  \frac{1}{N}\sum_{i=1}^{N}Q_p^{(i)}(\xi^{(i)})\le r^p \Big\},
	\end{equation}
	and its $k$-th order moment relaxation (with respect to the inner maximization)
	\begin{equation}\label{eq:Vlowk}
		V_{p,k}(x;r^p) \coloneqq \max_{y^{(i)}\in \bar{\mathscr{S}}[h]_{2k},i\in[N]} 
		\Big\{
		\frac{1}{N}\sum\limits_{i=1}^{N}\langle F_x, y^{(i)}\rangle 
		:\frac{1}{N}\sum_{i=1}^{N}\langle Q_p^{(i)}, y^{(i)}\rangle\le  r^p
		\Big\}.
	\end{equation}
	\begin{lemma}\label{lem:gap}
		Suppose $F$ is a polynomial in $\xi$.
		For each $r\ge 0$ and $k\ge d_1$,
		\[\begin{array}{ccc}
			v_p(x;r^p) & \le & v_{p,k}(x;r^p) \\
			\rotatebox{90}{$\le$} &  & \rotatebox{90}{$\le$} \\
			V_p(x;r^p) & \le & V_{p,k}(x;r^p). \\
		\end{array}
		\]
		If $r>0$, then $0\le v_{p,k}(x;r^p)-v_p(x;r^p)\le V_{p,k}(x;r^p)-V_p(x;r^p)$.
	\end{lemma}
	\begin{proof}
		As \eqref{eq:vf_k} and \eqref{eq:Vlowk} are moment relaxations of 
		\eqref{eq:ValueFunc} and \eqref{eq:UnderApprox}, respectively, 
		we have $v_p(x;r^p)\le v_{p,k}(x;r^p)$ and $V_p(x;r^p)\le V_{p,k}(x;r^p)$.
		The weak Lagrangian duality implies
		$V_p(x;r^p)\le v_p(x;r^p)$ and $V_{p,k}(x;r^p)\le v_{p,k}(x;r^p)$.
		When $r>0$, $y^{(i)} = [\hat{\xi}^{(i)}]_{2k}$ for $i\in[N]$ 
		is strictly feasible to \eqref{eq:Vlowk}.
		By Slater's condition, the strong duality holds between \eqref{eq:vf_k} and \eqref{eq:Vlowk},
		making $v_{p,k}(x;r^p)-v_p(x;r^p)\le V_{p,k}(x;r^p)-V_p(x;r^p)$.
		
	\end{proof}
	\begin{theorem}\label{thm:Or2gap}
		Under conditions of Theorem~\ref{thm:asym1str0}, 
		suppose $F(x,\xi)$ is a polynomial in $\xi$ 
		and $-h_1,\ldots, -h_{m_1}$ are all SOS-convex.
		Then $v_{p,d_1}(x;r^p)-v_p(x;r^p)=\mc{O}(r^2)$ for all $r>0$ that is sufficiently small.
	\end{theorem}
	\begin{proof}
		By Lemma~\ref{lem:gap}, 
		it suffices to show $V_{p,d_1}(x;r^p)-V_p(x;r^p)=\mc{O}(r^2)$ 
		when $r>0$ is sufficiently small.
		Suppose $(y^{(1)},\ldots, y^{(N)})$ is an optimizer of \eqref{eq:Vlowk} for $k=d_1$.
		Since $y^{(i)}\in\bar{\mathscr{S}}[h]_{2d_1}$, we must have $\langle h_j,y^{(i)}\rangle \ge 0$
		for each $j\in[m_1]$,
		as $\langle h_j,y^{(i)}\rangle$ is the (1,1)-entry of the psd matrix $L_{h_j}^{(d_1)}[y^{(i)}]$; 
		check \eqref{eq:S_2k} for the defining equation of $\bar{\mathscr{S}}[h]_{2d_1}$.
		Let $\pi(y^{(i)}) \coloneqq (y_{e_1}^{(i)},\ldots, y_{e_{n_0}}^{(i)})$ for each $i\in[N]$. 
		Since each $-h_j$ is SOS-convex, 
		the Jensen's inequality (\cite[Theorem~7.1.6]{NieBook}, 
		\cite[Theorem~2.6]{Las09SOSconvex}) ensures that
		\[
		-h_j(\pi(y^{(i)}))\le \langle -h_j,y^{(i)}\rangle \le 0,
		\quad \forall j\in[m_1],\, i\in[N].
		\]
		Since $\|\xi\|^p$ is SOS-convex,  
		the sum of its translation $\sum_{i=1}^N  Q_p^{(i)}(\xi)$
		is also SOS-convex. By using Jensen's inequality again, we obtain
		\[
		\sum\limits_{i=1}^N Q_p^{(i)}(\pi(y^{(i)})) 
		\le \sum\limits_{i=1}^N \langle Q_p^{(i)}, y^{(i)}\rangle \le Nr^p.
		\]
		This implies that $(\pi(y^{(1)}),\ldots, \pi(y^{(N)}))$ is a feasible point of \eqref{eq:UnderApprox}.
		Denote $w^{(i)} = [\pi(y^{(i)})]_{2d_1}$ for $i\in[N]$. 
		It is clear that $(w^{(1)}, \ldots, w^{(N)})$ is feasible to \eqref{eq:Vlowk}
		and $\pi(y^{(i)}) = \pi(w^{(i)}) = (w_{e_1}^{(1)},\ldots, w_{e_{n_0}}^{(i)})$ for each $i$.
		Hence, 
		\[\begin{aligned}
			& V_{p,d_1}(x;r^p)-V_p(x;r^p)
			\le \frac{1}{N}\sum\limits_{i=1}^N \big[
			\langle F_x, y^{(i)}\rangle - F_x(\pi(y^{(i)}))
			\big] \\
			& = 
			\frac{1}{N}\sum\limits_{i=1}^N\langle F_x, y^{(i)}-w^{(i)}\rangle
			= \frac{1}{N}\sum\limits_{i=1}^N\sum\limits_{|\alpha|=2}^{\deg(F_x)}
			(F_x)_{\alpha}( y_{\alpha}^{(i)}-w_{\alpha}^{(i)})\\
			& \le \Big\vert\frac{1}{N}\sum\limits_{i=1}^N\sum\limits_{|\alpha|=2}^{\deg(F_x)}
			(F_x)_{\alpha}( y_{\alpha}^{(i)}-(\hat{\xi}^{(i)})^{\alpha})\Big\vert
			+ \Big\vert\frac{1}{N}\sum\limits_{i=1}^N\sum\limits_{|\alpha|=2}^{\deg(F_x)}
			(F_x)_{\alpha}( w_{\alpha}^{(i)} -(\hat{\xi}^{(i)})^{\alpha})\Big\vert.
		\end{aligned}\]
		Since $(y^{(1)}, \ldots, y^{(N)})$ and $(w^{(1)},\ldots, w^{(N)})$ 
		are both feasible to \eqref{eq:Vlowk},
		with the estimation \eqref{eq:sdbd},
		the conclusion is therefore implied by Theorem~\ref{thm:asym1str0}.
\end{proof} 

	We remark that the SOS-convexity assumption on $-h_1,\dots,-h_{m_1}$ 
	in Theorem~\ref{thm:Or2gap} cannot be relaxed to $\Xi$ being a convex set, 
	as shown by the next example.
	\begin{example}
		Consider $X=\{0\}$, $f(x)=0$, $F_x(\xi)=\xi_1-\xi_3$,
		\[
		\Xi=\{\xi=(\xi_1,\xi_2,\xi_3)\in\bbR^3: 
		h_1(\xi)=\xi_3^2-\xi_1^2-\xi_2^2\ge0,\,h_2(\xi)=\xi_3\ge0\}
		\]
		is the second-order cone, 
		with $p=2$, $N=1$, $\hat{\xi}^{(1)}=(0,0,0)$.
		Since $\xi_3^2\ge\xi_1^2$ on $\Xi$, 
		we see that $F_x(\xi)\le0$.
		Thus for any $r\ge0$, the value function~\eqref{eq:ValueFunc} is
		\[
		v_2(0;r^2)=\min_{\lambda\in\bbR_{\ge0}}\Bigl\{
		\lambda r^2+\sup_{\xi\in\Xi}\bigl[\xi_1-\xi_3-\lambda(\xi_1^2+\xi_2^2+\xi_3^2)\bigr]
		\Bigr\}=0.
		\]
		Note that $p_1=d_1=1$, 
		so the first order ($k=1$) moment relaxation~\eqref{eq:vf_k} 
		can be used here 
		\[
		v_{2,1}(0;r^2)=\min_{\lambda\in\bbR_{\ge0}}
		\Bigl\{\lambda r^2+\max_{y\in\bar{\scrS}[h]_2}
		\pAngle{\xi_1-\xi_3-\lambda(\xi_1^2+\xi_2^2+\xi_3^2)}{y}
		\Bigr\}.
		\]
		We claim that $v_{2,1}(0,r^2)\ge\frac{r}{\sqrt{2}}$,
		so we cannot have the $\calO(r^2)$ bound as in Theorem~\ref{thm:Or2gap}.
		To see this, we explicitly write the inner maximization problem as
		\[
		\begin{aligned}
			\max_{y}\quad &y_{100}-y_{001}-\lambda(y_{200}+y_{020}+y_{002})\\
			\suchthat\quad & y_{001}\ge0,\\
			& y_{002}-y_{200}-y_{020}\ge0,\\
			& M_1[y]=\begin{bmatrix}
				1 & y_{100} & y_{010} & y_{001} \\
				y_{100} & y_{200} & y_{110} & y_{101} \\
				y_{010} & y_{110} & y_{020} & y_{011} \\
				y_{001} & y_{101} & y_{011} & y_{002}
			\end{bmatrix}\succeq0.
		\end{aligned}
		\]
		Take $y_{100}=r/\sqrt{2}$, $r_{200}=r_{002}=r^2/2$, 
		and all others to be zero.
		This gives a feasible solution to the inner maximization, 
		since the first two constraints are clearly satisfied, 
		while the moment matrix
		\[
		M_1[y]=\begin{bmatrix}
			1 & \frac{r}{\sqrt{2}} & 0 & 0 \\
			\frac{r}{\sqrt{2}} & \frac{r^2}{2} & 0 & 0 \\
			0 & 0 & 0 & 0 \\
			0 & 0 & 0 & \frac{r^2}{2}
		\end{bmatrix}\succeq0,
		\]
		since its Schur complement of the top-left entry is
		\[
		\begin{bmatrix}
			\frac{r^2}{2} & 0 & 0 \\
			0 & 0 & 0 \\
			0 & 0 & \frac{r^2}{2} 
		\end{bmatrix}-\begin{bmatrix}
			\frac{r}{\sqrt{2}} \\ 0 \\ 0
		\end{bmatrix}\begin{bmatrix}
			\frac{r}{\sqrt{2}} & 0 & 0
		\end{bmatrix} = \begin{bmatrix}
			0 & 0 & 0 \\
			0 & 0 & 0 \\
			0 & 0 & \frac{r^2}{2} 
		\end{bmatrix}\succeq0.
		\]
		Therefore, we obtain
		\[
		v_{2,1}(0;r^2)\ge\min_{\lambda\in\re_{\ge 0}}
		\left\{\lambda r^2+(\frac{r}{\sqrt{2}}-\lambda r^2)\right\} = \frac{r}{\sqrt{2}}.
		\]
\end{example} 

\section{Moment Relaxations for Two-Stage Wasserstein DRO}
\label{sec:2stage}

In this section, we focus on \eqref{eq:WDRO} 
where $F$ is a recourse function given by \eqref{eq:RecourseDual}.
For convenience, we repeat the explicit expressions of the recourse objective 
$G_x$ and the constraint tuple $g$ as
\[
G_x(\xi,u) = u^{\transpose}B(\xi)x+b(\xi)^{\transpose}u+d(\xi),\quad
g(\xi,u) = c(\xi)-A^{\transpose}u.
\]
Denote the joint feasible region 
$U \coloneqq \{(\xi,u): h(\xi)\ge 0,\, g(\xi,u)\ge 0\}$, 
which we assume is a regular closed set.
Let $p_1=p/2$ and 
\begin{equation}\label{eq:d_g}
		d_g\coloneqq \max\{ 2\deg(B), 2\deg(b),
		\deg(G_x), \deg(g)\},
\end{equation}
\begin{equation}\label{eq:d2}
	d_2 \coloneqq \max\{\lceil d_g/2\rceil, \, \lceil \deg(h)/2\rceil,\, p_1 \}.
\end{equation}
This $d_2$ is the lowest relaxation order that admits 
all information about $G_x$, $g$, $h$, and $Q_p^{(i)}$ for each $i\in[N]$.

\subsection{Standard moment relaxations}\label{sec:2stage-standard}

To distinguish it from the single-stage case, 
we denote the two-stage value function 
\begin{equation}\label{eq:2stageVF}
	\bar{v}_p(x;r^p )\coloneqq \min\limits_{\lambda\in\re_{\ge 0}}
	\Big\{ \lambda r^p+\frac{1}{N}\sum\limits_{i=1}^N\max\limits_{(\xi^{(i)},u^{(i)})\in U} 
	\big[ G_x(\xi^{(i)},u^{(i)})-\lambda Q_p^{(i)}(\xi^{(i)})\big]\Big\}.
\end{equation}
For the special case that $r = 0$, we set 
\[
\bar{v}(x;0) \coloneqq \frac{1}{N}\sum\limits_{i=1}^N 
\max_{u\in\re^{n_2}} \big\{G(\hat{\xi}^{(i)}, u): g(\hat{\xi}^{(i)},u)\ge 0 \big\}.
\]
For $k\ge d_2$, the $k$-th order moment relaxation of \eqref{eq:2stageVF} 
with respect to the inner maximization  
\begin{equation}\label{eq:2stage_vf_k}
	\bar{v}_{p,k}(x;r^p)\coloneqq \Big\{
	\lambda r^p+\frac{1}{N}\sum\limits_{i=1}^N \max\limits_{y^{(i)}\in\bar{\mathscr{S}}[g,h]_{2k}}
	\langle G_x-\lambda Q_p^{(i)},y^{(i)}\rangle\Big\}.
\end{equation}
\begin{remark}\label{rem:2stage_opt}
	In the following, we focus on bounding the discrepancy 
	$\bar{v}_{p,k}(x;r^p)-\bar{v}_p(x;r^p) = \mc{O}(r)$. 
	Such a bound will imply the $\mc{O}(r)$-optimality gap 
	in Theorem~\ref{thm:2stage-intro} similar to~Remark~\ref{rem:1stage_opt} 
	in the single-stage case.
\end{remark}
Consider the Lagrangian dual of \eqref{eq:2stageVF}: 
\begin{equation}\label{eq:2stage_Vstar}
	\bar{V}_p(x;r^p) \coloneqq \max\limits_{\substack{(\xi^{(i)}, u^{(i)})\in U\\
			i\in[N]}}
	\Big\{
	\frac{1}{N}\sum\limits_{i=1}^N G_x(\xi^{(i)}, u^{(i)}): 
	\frac{1}{N}\sum\limits_{i=1}^N Q_p^{(i)}(\xi^{(i)})\le r^p
	\Big\}.
\end{equation}
This is a polynomial optimization problem. 
Its $k$-th order moment relaxation is
\begin{equation}\label{eq:2stage_Vk}
	\bar{V}_{p,k}(x;r^p)\coloneqq 
	\max\limits_{\substack{y^{(i)}\in\bar{\mathscr{S}}[g,h]_{2k}\\
			i\in[N]}}
	\Big\{\frac{1}{N}\sum\limits_{i=1}^N 
	\langle G_x, y^{(i)}\rangle: \frac{1}{N}
	\sum\limits_{i=1}^N \langle Q_p^{(i)},y^{(i)}\rangle \le r^p
	\Big\},
\end{equation}
which is also the dual problem of \eqref{eq:2stage_vf_k}. 
Since $U$ has a nonempty interior, 
the strong duality holds between \eqref{eq:2stage_vf_k} 
and \eqref{eq:2stage_Vk} for all $r>0$.
The proof of the following lemma is almost a word-by-word repetition of 
that of~Lemma~\ref{lem:gap}, and thus omitted.

\begin{lemma}\label{lem:2stage_gap}
	For each $r\ge 0$ and $k\ge d_2$, it holds that
	\[
	\begin{array}{ccc}
		\bar{v}_p(x;r^p) & \le & \bar{v}_{p,k}(x;r^p) \\
		\rotatebox{90}{$\le$} &  & \rotatebox{90}{$\le$} \\
		\bar{V}_p(x;r^p) & \le & \bar{V}_{p,k}(x;r^p). \\
	\end{array}
	\]
	If $r>0$, then $0\le \bar{v}_{p,k}(x;r^p)-\bar{v}_p(x;r^p)
	\le \bar{V}_{p,k}(x;r^p)-\bar{V}_p(x;r^p).$
\end{lemma}

\begin{theorem}\label{thm:2stage_gap}
	For every $x\in X$ and $k\ge \max\{d_2, \deg(G_x)-p_1\}$, 
	\[
	\bar{v}(x;0) = \bar{V}_p(x;0) = \bar{V}_{p,k}(x;0).
	\]
\end{theorem}
\begin{proof}
	Suppose $r = 0$ and $k\ge \max\{d_2,  \deg(G_x)-p_1\}$. 
	Since $(\hat{\xi}^{(1)},\ldots, \hat{\xi}^{(N)})$ 
	is the only feasible option of \eqref{eq:2stage_Vstar}, 
	we have $\bar{V}_p(x; 0) = \bar{v}(x,0)$.
	Suppose $(y^{(1)},\ldots, y^{(N)})$ is a feasible point of \eqref{eq:2stage_Vk}.
	For each $i\in[N]$ and $\alpha\in \N^{n_0}, \beta\in \N^{n_2}$, 
	denote the translated tms $\tilde{y}^{(i)}$ and the corresponding translation matrix
	\begin{align}\label{eq:2stage_trans}
		\tilde{y}^{(i)} & = (\langle (\xi-\hat{\xi}^{(i)})^{\alpha}u^{\beta}, y^{(i)}\rangle)_{(\alpha,\beta)\in \N_{2k}^{n_0+n_2}},\\
		\label{eq:2stage_transmat}
		T_i &= [\mbox{vec}((\xi-\hat{\xi}^{(i)})^{\alpha}u^{\beta})]_{(\alpha,\beta)\in\N_{2k}^{n_0+n_2}}.
	\end{align}
	Since $M_k[\tilde{y}^{(i)}] = T_i^{\transpose}M_k[y^{(i)}]T_i\succeq 0$, 
	all its diagonal entries are nonnegative. On the other hand, 
	by Lemma~\ref{lem:norm} (with $\theta_{\min} = 1$),
	\[
	\sum\limits_{i=1}^N\sum\limits_{l=1}^{n_0}\tilde{y}^{(i)}_{2p_1e_l,0}
	= \sum\limits_{i=1}^N\langle \|\xi\|_p^p,\tilde{y}^{(i)}\rangle
	\le \sum\limits_{i=1}^N\langle Q_p^{(i)},y^{(i)}\rangle = 0.
	\]
		This means each involved $\tilde{y}^{(i)}_{2p_1e_l,0}= 0$.
		Then $\tilde{y}^{(i)}_{te_l} = 0$ for all $t\in[p]$ and $l\in[n_0]$ 
		by Lemma~\ref{lem:multi_mom}.
		In addition, if $p_1<t\le p_1+k$, then $\tilde{y}_{te_l}^{(i)} = 0$
		since the 2-by-2 minor 
		$\tilde{y}^{(i)}_{2(t-k)e_l}\tilde{y}^{(i)}_{2ke_l}-(\tilde{y}_{te_l}^{(i)})^2\ge 0$.
		For any $\alpha\in\N^{n_0}$ with $0<|\alpha|\le p_1+k$, 
		$\tilde{y}_{\alpha,0}^{(i)} = 0$
		by the nonnegativity of $2$-by-$2$ minors.
		Similarly, we can get $\tilde{y}^{(i)}_{\alpha,e_t} = 0$ for all nonzero 
		$(\alpha,e_t)\in \N_{p_1+k}^{n_0,n_2}$ by using the nonnegativity of associated 2-by-2 minors 
		of $M_k[\tilde{y}^{(i)}]$.
		In addition, for each $t\in [n_2]$, 
		$\tilde{y}^{(i)}_{\alpha, e_t} = 0$ for all nonzero $(\alpha, e_t)\in \N_{p_1+k}^{n_0}$,
		since there always exists
		some $\beta\in \N_{p_1+k}^{n_0}$ that makes 
		\[\mbox{either}\quad
		\left\vert\begin{array}{ll} 
			\tilde{y}_{2(\alpha-\beta),2e_t}^{(i)} & \tilde{y}_{\alpha,e_t}^{(i)}\\
			\tilde{y}_{\alpha,e_t}^{(i)} & 0 \end{array}\right\vert\ge 0,\quad\mbox{or}\quad
		\left\vert\begin{array}{ll} 
			0 & \tilde{y}_{\alpha,e_t}^{(i)}\\
			\tilde{y}_{\alpha,e_t}^{(i)} &  \tilde{y}_{2(\alpha-\beta),2e_t}^{(i)}  
		\end{array}\right\vert\ge 0. 
		\]
	Let $u^{(i)} = \pi_u(y^{(i)})\coloneqq (\tilde{y}^{(i)}_{0,e_1},\ldots, \tilde{y}^{(i)}_{0,e_{n_2}})$ 
	be the projection of $y^{(i)}$ onto the $u$-space and let $\mathbf{0}\in\re^{n_0}$ be the constant zero vector.
	It holds that
	\[
	\forall (\alpha,\beta)\in\N_{p_1+k}^{n_0,n_2},\, |\beta|\le 1:\quad
	\tilde{y}^{(i)}_{\alpha,\beta} = \mathbf{0}^{\alpha}(u^{(i)})^{\beta}\,\Rightarrow\,
	y^{(i)}_{\alpha, e_t} = (\hat{\xi}^{(i)})^{\alpha}(u^{(i)})^{\beta}.
	\]
	Since $\deg(G_x)\le k+p_1$ and $G_x$ is linear in $u$, 
	we get $G_x(\hat{\xi}^{(i)}, u^{(i)}) = \langle G_x, y^{(i)}\rangle$. 
	Note that $u^{(i)}$ has a one-to-one relation to $y^{(i)}$ and
	$(\hat{\xi}^{(i)}, u^{(i)})\in U$ by $y^{(i)}\in \bar{\mathscr{S}}[g,h]_{2k}$. 
	For each $i\in [N]$, we have
	\[
	\max_{u\in\re^{n_2}} \big\{G_x(\hat{\xi}^{(i)},u): g(\hat{\xi}^{(i)}, u)\ge 0\big\} 
	= \max_{y\in \bar{\mathscr{S}}[g,h]_{2k}}
	\big\{ \langle G_x, y\rangle : \langle Q_p^{(i)}, y\rangle = 0 \big\}.
	\]
	The conclusion follows by averaging above both sides over $i\in [N]$.
	
\end{proof}

We then study the $\calO(r)$-consistency property for the moment relaxation \eqref{eq:2stage_vf_k} when $r>0$.
Let $y\in\bar{\mathscr{S}}[g,h]_{2k}$. 
The projections of $y$ onto the $\xi$-space and $u$-space are respectively denoted by
\[
\pi_{\xi}(y) \coloneqq (y_{e_1,0},\ldots, y_{e_{n_0},0}),\quad
\pi_u(y) \coloneqq (y_{0,e_{1}},\ldots, y_{0,e_{n_2}}).
\]

\begin{theorem}\label{thm:2stage_idc}
	Under Assumption~\ref{assum}, 
	suppose $X$ is bounded and there exist $R_0>0, k_0\ge 1$ 
	such that $R_0-\|u\|^2\in \qmod{g,h}_{2k_0}$. 
	Assume the relaxation order $k\ge \max\{d_2,k_0\}$.
	For all $r>0$ that is sufficiently small, 
	the difference $\bar{v}_{p,k}(x;r^p)-\bar{v}_p(x;r^p) = \mathcal{O}(r)$ if either
	\begin{enumerate}
		\item[(i)] $p\ge d_g$
			for $d_g$ in \eqref{eq:d_g}; or
		\item[(ii)] there exists $R>1, k_1\in \N$ such that 
		$R-\|\xi\|^2\in \qmod{h}_{2k_1+2}$, 
		and $d_g\le 2k-2k_1+p$.
	\end{enumerate}
	In either case, any optimal solution $x$ to the moment relaxation \eqref{eq:2stage} is an $\mc{O}(r)$-optimal 
	solution to the ESO \eqref{eq:ESO} by Remark~\ref{rem:2stage_opt}.
\end{theorem}
\begin{proof}
		By Lemma~\ref{lem:2stage_gap} and Theorem~\ref{thm:2stage_gap},
		it suffices to show $\bar{V}_{p,k}(x;r^p)-\bar{V}_p(x;0)=\mc{O}(r)$
		for all sufficiently small $r>0$.
		Let $(y^{(1)},\ldots,y^{(N)})$ be an optimizer of \eqref{eq:2stage_Vk}.
		For each $i\in [N]$, 
		let $\tilde{y}^{(i)}$ be the translated tms of $y^{(i)}$ as in \eqref{eq:2stage_trans}, i.e.,
		\[
		\tilde{y}^{(i)} = (\langle (\xi-\hat{\xi}^{(i)})^{\alpha}u^{\beta}, y^{(i)}\rangle)_{(\alpha,\beta)\in \N_{2k}^{n_0+n_2}}.
		\]
		This translation is independent with $u$. The membership $R_0-\|u\|^2\in\qmod{g,h}_{2k}$ implies that 
		$\langle R_0-\|u\|^2,\tilde{y}^{(i)}\rangle \ge 0$.
		For each $i\in [N]$, define
		\begin{align}\label{eq:inner_i}\tag{$LT_i$}
			F_x^{(i)} & \coloneqq \max\limits_{u}\, 
			\{G_x(\hat{\xi}^{(i)}, u): c(\hat{\xi}^{(i)})-A^{\transpose}u\ge 0\},\\
			\label{eq:mom_inner_i}\tag{$LT_i^k$}
			F_x^{(i,k)} & \coloneqq\max\limits_{u} \,
			\{ G_x(\hat{\xi}^{(i)}, u): \langle c, y^{(i)}\rangle-A^{\transpose}u\ge 0\},
		\end{align}
		where $\langle c, y^{(i)}\rangle$ is the vector obtained by applying the bilinear operation \eqref{eq:bilinear} entrywise between $c$ and $y^{(i)}$.
		Since \eqref{eq:inner_i} and \eqref{eq:mom_inner_i} are linear programs with the same objective function and constraint matrix,
		there exists a positive constant $\mc{G}>0$ such that 
		\begin{equation}\label{eq:Fx_diff}
			|F_x^{(i)} - F_x^{(i,k)}| \le \mc{G}\|c(\hat{\xi}^{(i)})- \langle c, y^{(i)}\rangle\|_{\infty}.
		\end{equation}
		By Theorem~\ref{thm:2stage_gap}, $\bar{V}_p(x;0) = \bar{v}(x;0) = \frac{1}{N}\sum_{i=1}^N F_x^{(i)}$. 
		Observe that $\pi_u(y^{(i)})$ is a feasible point of \eqref{eq:mom_inner_i}
		since $y^{(i)}\in \bar{\mathscr{S}}[g,h]_{2k}$.
		Then
		\[\begin{aligned}
			\bar{V}_{p,k}(x;r^p)-\bar{V}_p(x;0) 
			= \frac{1}{N}\sum\limits_{i=1}^N \big( \langle G_x, y^{(i)}\rangle -F_x^{(i,k)} \big) + \frac{1}{N} \sum\limits_{i=1}^N \big(F_x^{(i,k)}-F_x^{(i)}\big)\\
			\le \frac{1}{N}\sum\limits_{i=1}^N \big| \langle G_x, y^{(i)}\rangle-G_x(\hat{\xi}^{(i)}, \pi_u(y^{(i)})) \big|+\frac{1}{N}\sum\limits_{i=1}^N |F_x^{(i,k)}-F_x^{(i)}|.
		\end{aligned}\]
		Since $G_x$ is linear in $u$, $\langle G_x, y^{(i)}\rangle-G_x(\hat{\xi}^{(i)}, \pi_u(y^{(i)}))$ is a weighted sum of $\tilde{y}^{(i)}_{\alpha,0}$ and $\tilde{y}^{(i)}_{\alpha,e_t}$
		for $0\neq\alpha\in\N^{n_0}, 2|\alpha|\le d_g$ and $t\in[n_2]$.
		Since $\tilde{y}^{(i)}_{0,2e_t}\le \langle \|u\|_2^2, \tilde{y}^{(i)}\rangle
		\le \langle \|u\|^2, \tilde{y}^{(i)}\rangle \le R_0$,
		Using Lemma~\ref{lem:multi_mom} and arguments for Theorem~\ref{thm:asym1str0}~(i), with $|\alpha|=d$,
		\[
		\sum\limits_{t=1}^{n_2}|\tilde{y}_{\alpha,e_t}^{(i)}|\le\left\{
		\begin{array}{lr}
			\sqrt{n_2R_0}\sqrt{\tilde{y}_{2\alpha,0}^{(i)}}, 
			&\quad \mbox{if $d=1$};\\
			\sqrt{n_2R_0}\sum\limits_{l=1}^{n_0}\frac{2\alpha_l}{2d}
			\tilde{y}^{(i)}_{2de_l,0},
			&\quad \mbox{otherwise}.
		\end{array}
		\right.\]
		
		(i) Suppose $d_g\le p$. 
		Since $X$ is bounded, by Taylor's expansion and arguments of Theorem~\ref{thm:asym1str0} (i), there exists $\mc{F}^*>0$
		such that for arbitrarily chosen $x\in X$,
		\[\begin{aligned}
			& {\begin{array}{rl}\|c(\hat{\xi}^{(i)})-\langle c,y^{(i)}\rangle\|_{\infty}
					& \le \mc{F}^*(1+\|\hat{\xi}^{(i)}\|_p)^{p-1}
					\sum\limits_{l=1}^{n_0} |\tilde{y}_{e_l,0}^{(i)}| \\
					& \quad +
					2\mc{F}^*\sum\limits_{d=1}^{p_1} (1+\|\hat{\xi}^{(i)}\|_p)^{p-2d}
					\langle \|\xi\|_{2d}^{2d}, \tilde{y}^{(i)}\rangle,
			\end{array}}\\
			& {\begin{array}{r} 
					|\langle G_x, y^{(i)}\rangle-G_x(\hat{\xi}^{(i)}, \pi_u(y^{(i)}))|
					\le (1+\sqrt{n_2R_0})\mc{F}^*(1+\|\hat{\xi}^{(i)}\|_p)^{p-1}
					\sum\limits_{l=1}^{n_0} |\tilde{y}_{e_l,0}^{(i)}|\\
					+
					2(1+\sqrt{n_2 R_0})\mc{F}^*\sum\limits_{d=1}^{p_1} (1+\|\hat{\xi}^{(i)}\|_p)^{p-2d}
					\langle \|\xi\|_{2d}^{2d}, \tilde{y}^{(i)}\rangle. 
			\end{array}}
		\end{aligned}\]
		Then the conclusion is then implied by Theorem~\ref{thm:asym1str0}~(i).

		(ii) Suppose $d_g\le 2k-2k_1+p$ and that there exists $R>1, k_1\in \N$ such that $R-\|\xi\|^2\in \qmod{h}_{2k_1+2}$.
		Then $U$ is bounded and there exists a constant $\hat{\mc{F}}^*$ such that 
		for arbitrarily chosen $x\in X$,
		\[\begin{aligned}
			& 
			\|c(\hat{\xi}^{(i)})-\langle c,y^{(i)}\rangle\|_{\infty}
			\le \hat{\mc{F}}^*
			\sum\limits_{l=1}^{n_0} |\tilde{y}_{e_l,0}^{(i)}|+
			2\hat{\mc{F}}^*\sum\limits_{d=1}^{\lceil \frac{d_g}{2}\rceil} 
			\langle \|\xi\|_{2d}^{2d}, \tilde{y}^{(i)}\rangle, \\
			& {\begin{array}{rcl} 
					|\langle G_x, y^{(i)}\rangle-G_x(\hat{\xi}^{(i)}, \pi_u(y^{(i)}))|
					& \le & (1+\sqrt{n_2R_0})\hat{\mc{F}}^*
					\sum\limits_{l=1}^{n_0} |\tilde{y}_{e_l,0}^{(i)}|\\
					& & \quad +
					2(1+\sqrt{n_2 R_0})\hat{\mc{F}}^*\sum\limits_{d=1}^{\lceil \frac{d_g}{2}\rceil}
					\langle \|\xi\|_{2d}^{2d}, \tilde{y}^{(i)}\rangle. 
			\end{array}}
		\end{aligned}\]
		Then the conclusion is implied by 
		Theorem~\ref{thm:asym1str0}~(ii).
\end{proof}

	Unlike the special cases in Section~\ref{ssc:tighter_gap}, 
	the following example shows that the $\calO(r)$ 
	estimate of the relaxation gap cannot be refined in the two-stage case 
	even when the defining polynomials $-h_1,\dots,-h_{m_1}$ are SOS-convex. 
	\begin{example}\label{ex:2stage-Or}
		Consider $X=\{0\}, f(x)=0$, $\Xi=\{\xi\in\bbR:h(\xi):=2\xi-\xi^2\ge0\}=[0,2]$, with $N=1$, $\hat{\xi}^{(1)}=0$, and the recourse function $F_x(\xi)=\max_{u}\{-\xi u:g(u):=u\ge0\}=0$ where the maximum is attained at $u=0$.
		Thus, we can use $p=\deg(F_x)=2$, in which case the corresponding~\eqref{eq:WDRO} is
		\begin{equation*}
			\min_{\lambda\in\bbR_{\ge0}}\Big\{\lambda r^2+\sup_{\xi\in[0,2],u\ge0}[-\xi u-\lambda\xi^2]\Big\}=0.
		\end{equation*}
		Note that $-h(\xi)$ is SOS-convex since it is quadratic and convex.
		We next show that the optimal value of the degree-2 moment 
		relaxation~\ref{eq:2stage} is at least $r/2$.
		To see this, we use normalized tms $y\in\bar{\scrS}[g,h]_2$ 
		with $y_{ij}$ denoting the pseudomoment of $\xi^i u^j$ 
		(e.g., $y_{10}$ corresponding to $\xi$ and $y_{02}$ corresponding to $u^2$).
		Then the moment relaxation~\eqref{eq:2stage} can be written as
		\begin{equation*}
			\min_{\lambda\in\bbR_{\ge0}}\Big\{\lambda r^2+\sup_{y\in\bar{\scrS}[2\xi-\xi^2,u]_2}[-y_{11}-\lambda y_{20}]\Big\}\ge\frac{r}{2}.
		\end{equation*}
		Note that the moment matrix
		\begin{equation*}
			M_1[y]=\begin{bmatrix}
				1 & y_{10} & y_{01} \\
				y_{10} & y_{20} & y_{11} \\
				y_{01} & y_{11} & y_{02}
			\end{bmatrix}=\begin{bmatrix}
				1 & \frac{1}{2} & \frac{r}{2} \\
				\frac{1}{2} & 1 & -\frac{r}{2} \\
				\frac{r}{2} & -\frac{r}{2} & r^2
			\end{bmatrix}
		\end{equation*}
		is feasible, because $\pAngle{g}{y}=y_{01}>0$ and $\pAngle{h}{y}=2y_{10}-y_{20}=0$.
		Moreover, $M_1[y]\succeq0$ since Schur complement of the top-left entry is
		\begin{equation*}
			\begin{bmatrix}
				1 & -\frac{r}{2} \\ -\frac{r}{2} & r^2
			\end{bmatrix}-\begin{bmatrix}
				\frac{1}{2} \\ \frac{r}{2}
			\end{bmatrix}\begin{bmatrix}
				\frac{1}{2} & \frac{r}{2}
			\end{bmatrix}=\frac{3}{4}\begin{bmatrix}
				1 & -r \\
				-r & r^2
			\end{bmatrix}\succeq0.
		\end{equation*}
		Thus the inner supremum satisfies 
		$\sup_{y\in\bar{\scrS}[2\xi-\xi^2,u]_2}[-y_{11}-\lambda y_{20}]\ge \frac{r}{2}-\lambda r^2$, 
		which shows that the relaxation~\eqref{eq:2stage} is at least $r/2$.
	\end{example}

	In Theorem~\ref{thm:2stage_idc}, the condition that 
	$R_0-\|u\|^2\in \mathscr{S}[g,h]_{2k_0}$ ensures 
	the second-stage decision vector $u$ is bounded, 
	which is necessary for the conclusion to hold.
The following example shows that when the second-stage problem has 
an unbounded set of decisions $u$, the $\calO(r)$-consistency may fail.

\begin{example}\label{ex:2stage}
	Consider $X=\{0\}$, $f(x)=0$, $\Xi=\{\xi\in\bbR:1-\xi^2\ge0\}=[-1,1]$, 
	with $N=1$, $\hat{\xi}^{(1)}=0$, 
	and the recourse function $F_x(\xi):=\max_u\{-(1-\xi)u:u\ge0\}=0$.
	The growth rate $\Gamma_p(x)=0$ for any $p\ge1$, 
	and the corresponding~\eqref{eq:WDRO} is
	\begin{equation*}
		\min_{\lambda\in\re_{\ge0}}\Big\{
		\lambda r^2+\sup_{\xi\in\Xi,u\in\bbR_{\ge0}}\big[-(1-\xi)u-\lambda\xi^2 \big]
		\Big\}= 0.
	\end{equation*}
	However, its moment relaxation~\eqref{eq:2stage} using a normalized tms $y$ leads to
	\begin{equation*}
		\min_{\lambda\in\re_{\ge0}}\Big\{
		\lambda r^2+\sup_{y\in\bar{\scrS}[1-\xi^2,u]}\big[-y_{01}+y_{11}-\lambda y_{20}\big]
		\Big\}=+\infty,
	\end{equation*}
	as the inner supremum is unbounded for any $\lambda\ge0$.
	To see this, suppose the supremum is bounded from above by $v\in\bbR$ for some $\lambda\ge0$.
	Since $\Xi\times \re_{\ge 0}$ has a nonempty interior, 
	we know as before that
	\begin{equation*}
		v+(1-\xi)u+\lambda\xi^2=\sigma_0+\sigma_1(1-\xi^2)+\sigma_2 u,
	\end{equation*}
	for some $\sigma_0,\sigma_1,\sigma_2\in\Sigma[\xi,u]_{2k}$.
	However, this is not possible for any relaxation order $k$ 
	by~\cite[Theorem 2]{powers2005polynomials}.
\end{example}

In~Example~\ref{ex:2stage}, one can include more manifestly nonnegative 
polynomials on the half-strip $[-1,1]\times[0,+\infty)$ 
such as $(1-\xi)u\in\bbR[\xi,u]$ to resolve the unboundedness issue 
(see e.g.,~\cite{nguyen2012polynomials} for a discussion on the nonnegativity over half-strips). 
In general, this leads to the \emph{preorder}-based moment relaxations, 
which impose positive semidefinite constraints on 
localizing matrices associated with products of the defining polynomials $g$ and $h$.
However, such preorder-based moment relaxations are less favored by optimizers 
as they typically have many more constraints than the standard moment relaxations, 
e.g.,~\eqref{eq:2stage}, 
and it remains unclear to us whether there is any effective bound 
on the relaxation order $k$ of these moment relaxations to preserve the $\calO(r)$-consistency.
As an alternative, we exploit the trait of our data-driven $p$-WDRO~\eqref{eq:MeasConstr} 
as an optimization over measures, 
and derive strengthened two-stage moment relaxations 
that can guarantee the $\calO(r)$-consistency in low-degree situations.

\subsection{Strengthened two-stage moment relaxations}
\label{sec:2stage-strengthened}

Our main idea on strengthening the two-stage moment relaxations 
is to add constraints that are redundant in the optimization over measures~\eqref{eq:MeasConstr}, 
but conducive in ensuring the consistency of the resulting moment relaxations.
We begin with some definitions. Let
\begin{equation}\label{eq:g_eps}
	g_{\epsilon} = (g_{\epsilon}^{(i)})_{i=1}^N\quad \mbox{with}\quad g_{\epsilon}^{(i)}(\xi,u) \coloneqq (\epsilon^2-\|\xi-
	\hat{\xi}^{(i)}\|_2^2)\cdot g(\xi,u),
\end{equation}
where $g(\xi,u) = c(\xi)-A^{\transpose}u$ with $A\in\re^{n_2\times m_2}$. 
When $g(\xi,u)\ge 0$ does not give a compact set, 
we may assume $n_2\ge m_2$ and $A$ has full column rank without loss of generality.
For $j\in[n_2]$ and $J\subseteq [n_2]$, 
denote by $a_j^{\transpose}$ the $j$-th row of $A$, $A_J$ the submatrix of $A$ with rows in $J$. 
Recall that the linear recourse is defined by the dual pair
\begin{align}
	\label{eq:rep_min_linrec}\tag{$P$}
	\min_{x'\in\re^{m_2}}\Big\{ c(\xi)^\transpose x'+d(\xi): Ax'= B(\xi)x+b(\xi),\,x'\ge0\Big\},\\
	\label{eq:rep_max_linrec}\tag{$D$}
	\max_{u\in\re^{n_2}}\Big\{ 
	u^\transpose B(\xi)x+b(\xi)^\transpose u+d(\xi):c(\xi)-A^\transpose u\ge0
	\Big\}.
\end{align}
If $A_J$ is invertible, 
then $x'_J(x;\xi) \coloneqq A_{J}^{-1}[B(\xi)x+b(\xi)]_J$ 
is said to be a {\it basic solution} of \eqref{eq:rep_min_linrec}. 
For the special case that such an $x'_J(x;\xi)$ is an optimizer of \eqref{eq:rep_min_linrec}, 
the corresponding index set $J$ is said to be an \emph{optimal basis}.

Let $d_3\coloneqq \max\{d_2,\lceil \deg(g_{\epsilon})/2\rceil\}$.
Define the $k$-th order strengthened moment relaxation
\begin{equation}\label{eq:2stage_Vk_eps}
	\bar{V}_{p,k}^{\epsilon}(x;r^p)\coloneqq 
	\left\{\max\limits_{\substack{y^{(i)}\in\bar{\mathscr{S}}[h,g,g_{\epsilon}^{(i)}]_{2k}\\ i\in[N]}}\frac{1}{N}\sum\limits_{i=1}^N \langle G_x, y^{(i)}\rangle
	:\frac{1}{N}\sum\limits_{i=1}^N \langle Q_p^{(i)}, y^{(i)}\rangle \le r^p,
	\right\}.
\end{equation}
and its Lagrangian dual
\begin{equation}\label{eq:vbar_eps}
	\bar{v}_{p,k}^\epsilon(x;r^p)\coloneqq \min_{\lambda\in\re_{\ge 0}}\Big\{
	\lambda r^p+\frac{1}{N}\sum\limits_{i=1}^N \max\limits_{y^{(i)}\in \bar{\scrS}[h,g,g_{\epsilon}^{(i)}]_{2k}}
	\langle G_x-\lambda Q_p^{(i)},y^{(i)}\rangle\Big\}.
\end{equation}
When $\epsilon>0$ and $r>0$, the Slater's condition implies strong duality 
between~\eqref{eq:2stage_Vk_eps}--\eqref{eq:vbar_eps}.
Note that the polynomial constraints $g_\epsilon^{(i)}(\xi^{(i)},u^{(i)})\ge0$ 
are not necessarily redundant in~\eqref{eq:2stage_Vstar}. 
To justify the name of moment relaxations, 
we compare $\bar{v}_{p,k}^\epsilon(x;r^p)$ with $\bar{v}_p(x;r^p)$, 
which from the original formulation~\eqref{eq:MeasConstr} is
\begin{equation*}
	\begin{aligned}
		\bar{v}_p(x;r^p)=
		\sup_{\mu^{(i)}\in\calM(\Xi),i\in[N]}\quad &
		\frac{1}{N}\sum_{i=1}^{N}\int_{\Xi}
		\max_u\{G_x(\xi,u):g(\xi,u)\ge0\}\diff\mu^{(i)}(\xi),\\
		\suchthat\quad\quad\,\,\quad
		&\frac{1}{N}\sum_{i=1}^{N}\int_{\Xi} Q_p^{(i)}(\xi)\diff\mu^{(i)}(\xi)\le r^p.
	\end{aligned}
\end{equation*}
The inner maximization over $u$ is a linear optimization problem.
Its optimality can be attained at one of basic solutions, 
which correspond to a subset of columns of $A$, for any $\xi\in \Xi$.
Thus, by indexing these basic solutions and 
taking the smallest index whenever there are multiple optimal basic solutions, 
we can get a Borel measurable solution map $u^*:\Xi\to\bbR^{n_2}$ 
such that $\max_u\{G_x(\xi,u):g(\xi,u)\ge0\}=G_x(\xi,u^*(\xi))$.
Then the pushforward measure of $\mu^{(i)}$ induced by $(1,u^*):\Xi\to\Xi\times\bbR^{n_2}$ 
is a probability measure supported on $U\subset\Xi\times\bbR^{n_2}$.
Hence for any $\epsilon>0$, $r>0$, 
we have $\bar{v}_p(x;r^p)\le\bar{V}_{p,k}^{\epsilon}(x;r^p)=\bar{v}_{p,k}^\epsilon(x;r^p)$ since
\begin{align}
	& \bar{v}_p(x;r^p) \nonumber\\
	& \le \sup_{\substack{\gamma^{(i)}\in\calM(U)\\i\in[N]}}\Big\{
	\frac{1}{N}\sum_{i=1}^{N}\int_{U}G_x(\xi,u)\diff\gamma^{(i)}:
	\frac{1}{N}\sum_{i=1}^{N}\int_{U}Q_p^{(i)}(\xi)\diff\gamma^{(i)}\le r^p\Big\}\label{eq:stren_meas}\\
	&\le \sup_{\substack{y^{(i)}\in\bar{\scrS}[h,g,g_\epsilon^{(i)}]_{2k}\\i\in[N]}}
	\Big\{\frac{1}{N}\sum_{i=1}^{N}\pAngle{G_x}{y^{(i)}}:
	\frac{1}{N}\sum_{i=1}^{N}\pAngle{Q_p^{(i)}}{y^{(i)}}\le r^p\Big\}. \label{eq:stren_mom}
\end{align}
In the above, the second inequality follows as 
any feasible probability measure $\gamma^{(i)}$ in~\eqref{eq:stren_meas} 
determines a degree-$2k$ normalized tms $y^{(i)}$ 
that is feasible to~\eqref{eq:stren_mom}.
Moreover, we also have $\bar{v}_{p,k}^{\epsilon}(x;r^p)\le\bar{v}_{p,k}(x;r^p)$ 
since $\bar{\scrS}[h,g,g_\epsilon^{(i)}]_{2k}\subseteq\bar{\scrS}[h,g]_{2k}$.
These relations can be summarized as follows.
\begin{proposition}\label{prop:stren_valid}
	Suppose $\epsilon\ge Nr^p>0$.
	Then including $g_\epsilon$ from~\eqref{eq:g_eps} gives a valid moment relaxation, i.e., 
	\begin{equation}
		\bar{v}_p(x;r^p)\le \bar{v}_{p,k}^{\epsilon}(x;r^p)\le\bar{v}_{p,k}(x;r^p).
	\end{equation}
\end{proposition}
	We remark that the strengthened moment relaxation~\eqref{eq:vbar_eps} 
	has $m_1+2m_2+1$ LMI constraints, cp.\ $m_1+m_2+1$ LMIs in~\eqref{eq:2stage}, 
	while the number of tms variables remains $\binom{n_0+n_2+2k}{2k}$.

Now we proceed to study the $\calO(r)$-consistency of the strengthened relaxation.
We make the following assumption that 
basic optimal solutions of the second-stage linear optimization remain \emph{feasible} 
with respect to $\epsilon$-perturbations of the uncertain parameters $\xi$, 
which holds for sufficiently small $\epsilon>0$ in all \emph{nondegenerate} cases, 
i.e., all basic variables in the basic optimal solutions are strictly positive 
in the second-stage linear optimization~\eqref{eq:rep_min_linrec}.
\begin{assumption}\label{as:epsilon}
	For given $\epsilon>0$, there exists an optimal basis $J_i$ for each $i\in[N]$, 
	such that $x_{J_i}'(x;\xi) = A_{J_i}^{-1}[B(\xi)x+b(\xi)]$ 
	is feasible to \eqref{eq:rep_min_linrec}
	for every $\xi\in B_{\epsilon}(\xi^{(i)})$.
\end{assumption}
We then show that the $\calO(r)$-consistency is preserved for 
these strengthened moment relaxations 
without requiring $u$ to be bounded.
\begin{theorem}\label{thm:strengthned}
	Assume $X$ is bounded, $\deg(G_x)\le 2$ and $\deg(g)\le p$.
	Under Assumption~\ref{assum} and Assumption~\ref{as:epsilon}, if $k\ge d_3$, 
	then for all $r>0$ that is sufficiently small, we have
	\[  \bar{v}_{p,k}^{\epsilon}(x;r^p)-\bar{v}_p(x;r^p) = \mathcal{O} (r). \]
\end{theorem}
\begin{proof}
	For any fixed $N$, Proposition~\ref{prop:stren_valid} 
	ensures the validness of our strengthened relaxation when $r>0$ is sufficiently small, 
	making $\bar{v}_{p,k}^\epsilon(x;r^p)-\bar{v}_p(x;r^p)\ge0$.
	Let $F_x^{(i)}$ denote the optimal value of \eqref{eq:rep_min_linrec} at $\xi = \xi^{(i)}$.
	By Theorem~\ref{thm:2stage_gap}, we get $\bar{V}_p(x;0) = \frac{1}{N}\sum_{i=1}^N F_x^{(i)}$.
	Under Assumption~\ref{as:epsilon}, for every $i\in[N]$, 
	there is an optimal basis $J_i\subseteq [n_2]$ such that $F_x^{(i)} = c(\hat{\xi}^{(i)})^{\transpose}x'_{J_i}(\hat{\xi}^{(i)})+d(\hat{\xi}^{(i)})$ and $x'_{J_i}(\xi)\ge 0$ for all $\xi\in B_{\epsilon}(\hat{\xi}^{(i)})$.
	Note that $\max\{\deg(B), \deg(d)\}\le\deg(G_x)\le 2$, 
	and the ball $B_{\epsilon}(\hat{\xi}^{(i)})=\{\xi:\epsilon^2-\|\xi-\hat{\xi}^{(i)}\|_2^2\ge 0\}$ 
	is quadratically determined.
	By S-lemma \cite[Theorem~5.2.6]{NieBook}, there exists $\theta^{(i)}\in\re_{\ge0}^{m_2}$ such that 
	\[\begin{aligned}
		& x_{J_i}'(x;\xi)-(\xi^2-\|\xi-\hat{\xi}^{(i)}\|_2^2)\theta^{(i)}\\ 
		& = A_{J_i}^{-1}[B(\xi)x+b(\xi)]_{J_i}-(\xi^2-\|\xi-\hat{\xi}^{(i)}\|_2^2)\theta^{(i)}\in\Sig[\xi]^{m_2},
	\end{aligned}\]
	where $\Sig[\xi]$ is the SOS cone in $\xi$ and $\Sig[\xi]^m$
	denotes its $m$-fold Cartesian product. 
	For each $i\in [N]$, denote
	\[
	\phi_x^{(i)}(\xi,u) \coloneqq\, \big[ x_{J_i}'(x;\xi)-(\xi^2-\|\xi-\hat{\xi}^{(i)}\|_2^2)\theta^{(i)} \big]^{\transpose}
	( c(\xi)-A^{\transpose} u ).
	\]
	We have $\pAngle{\phi_x^{(i)}}{y^{(i)}}\ge0$ for every normalized tms $y^{(i)}$ 
	that is feasible to \eqref{eq:2stage_Vk_eps}.
	After simplification, we get
	\[
	\phi_x^{(i)}(\xi,u) = [c(\xi)^{\transpose}x_{J_i}'(x;\xi)+d(\xi)] - G_x(\xi,u)-(\theta^{(i)})^{\transpose} g_{x,\epsilon}^{(i)}(\xi,u),
	\]
	where $G_x(\xi,u) = u^\transpose B(\xi)x+b(\xi)^\transpose u+d(\xi)$.
	Suppose $(y^{(1)},\ldots, y^{(N)})$ is an optimizer of \eqref{eq:2stage_Vk_eps}. 
	By a proper expansion of 
	$\langle G_x, y^{(i)}\rangle\le 
	\langle G_x, y^{(i)}\rangle +\langle \phi_x^{(i)}, y^{(i)}\rangle$, 
	we obtain 
	\[
	\langle G_x, y^{(i)}\rangle \le  \langle c^{\transpose}x_{J_i}'+d, y^{(i)}\rangle -\langle (\theta^{(i)})^{\transpose} g_{x,\epsilon}^{(i)}, y^{(i)}\rangle  \le \langle c^{\transpose}x_{J_i}'+d, y^{(i)}\rangle.
	\] 
	Since $\bar{V}_p^{\epsilon}(x;r^p)\ge 
	\bar{V}_p(x;0)$, $\bar{V}_p(x;0) = \frac{1}{N}\sum_{i=1}^N F_x^{(i)}$
	and $\deg(g)\le p$, by Theorem~\ref{thm:2stage_idc}, 
	when $r>0$ is sufficiently small,
	\begin{align}\nonumber
		\bar{V}_{p,k}^{\epsilon}(x;r^p) - \bar{V}_p^{\epsilon}(x;r^p)
		& \le \bar{V}^{\epsilon}_{p,k}(x;r^p) - \bar{V}_p(x;0)\\
		\nonumber
		& \le \frac{1}{N}\sum\limits_{i=1}^N \langle c^{\transpose}x'_{J_i}+d, y^{(i)}-[\hat{\xi}^{(i)}]_{2k}\rangle
		= \mathcal{O}(r).
	\end{align}
	Finally, the conclusion follows from 
	$\bar{V}_{p,k}^\epsilon(x;r^p)=\bar{v}_{p,k}^\epsilon(x;r^p)$ 
	by Slater's condition for $r>0$, 
	and $\bar{v}_p(x;r^p)\ge\bar{V}_p(x;r^p)$ by~Lemma~\ref{lem:2stage_gap}.
	
\end{proof}

\section{Numerical Experiments}\label{sec:num}
	
	In this section, we experiment our moment relaxations~\eqref{eq:1stage} and~\eqref{eq:2stage} 
	on a box-constrained regression and a two-stage production problems.
	We begin with some common experimental setups.
	
	Since the goal of the experiments is to illustrate the asymptotic consistency 
	of our proposed moment relaxations as the training sample size $N$ grows, 
	we adopt the following procedure for each $N\in\{10,20,\dots,100\}$.
	\begin{description}
		\item[Step~1.] Take $N$ iid samples $\hat{\xi}^{(1)},\dots,\hat{\xi}^{(N)}$ 
		of the uncertain parameters $\xi$, and solve the corresponding ESO or 
		moment relaxations of the $p$-WDRO for candidate Wasserstein radii $r=r(N)$.
		In particular,
		\begin{itemize}
			\item when $r=0$, we solve the ESO~\eqref{eq:ESO} directly;
			\item and for all $r>0$, we solve the corresponding moment relaxation~\eqref{eq:1stage} or~\eqref{eq:2stage}. 
		\end{itemize}
		\item[Step~2.] Record the training objective value and an optimal decision $x$.
		\item[Step~3.] Estimate the out-of-sample test performance 
		(i.e., mean, standard deviation, median, 10\% and 90\% percentiles) of the decision $x$ 
		by redrawing $N^{{\rm test}}=10,000$ iid samples of $\xi$.
	\end{description}
	Here, $r(N)$ is a decreasing function that will be detailed in each of the problem descriptions below.
	To account for the variability induced by random training samples, 
	we conduct $N^{{\rm rep}}$ independent replications of the above experimental procedure and report the mean values $\pm$ standard deviations across replications in the tables of results (Tables~\ref{tab:regression} and~\ref{tab:production}).
	
	Our implementation\footnote{Code repository: \url{https://github.com/shixuan-zhang/MoWDRO.jl}} 
	uses \texttt{Julia} packages \texttt{JuMP} \cite{Lubin2023} and 
	\texttt{SumOfSquares} \cite{weisser2019polynomial,legat2017sos}, 
	where the underlying linear programs (LP) and quadratic programs (QP) are solved by 
	\texttt{Gurobi} (v13.0)~\cite{gurobi} and the underlying semidefinte programs (SDP) are solved by 
	\texttt{Mosek} (v11.0)~\cite{mosek} (both with default settings of feasibility tolerances and optimality thresholds).
	We use the level method as our first-order method to solve 
	our moment relaxations~\eqref{eq:1stage} and~\eqref{eq:2stage}, 
	which is terminated when either the absolute or the relative optimality gap reaches $10^{-4}$.
	Except for the comparison with the Hanasusanto-Kuhn semidefinite relaxation~\cite{hanasusanto2018conic}, 
	we distribute the evaluation of subgradients over 10 worker CPUs, each allocated with 30 GB of RAM, 
	as discussed in Algorithm~\ref{alg:Subgradient}.

	\subsection{Box-constrained regression}
	We first consider the following box-constrained regression problem.
	Suppose we have $n:=n_0-1$ independent variables $\zeta_1,\dots,\zeta_{n}\ge0$ and 
	one dependent variable $\omega$ that satisfy the following quadratic equation
	\begin{equation}
		\omega = \sum_{s=0}^{n}\sum_{t=0}^{s}c_{st}\zeta_s\zeta_t + \epsilon, 
	\end{equation}
	for some unknown but bounded parameters $c_{00},c_{10},\dots,c_{n0},\dots,c_{nn}\in[-1,1]$, 
	where by slight abuse of notation we write $\zeta_0=1$ to simplify the expression, 
	and $\epsilon\sim\operatorname{Normal}(0,\sigma^2)$ is a Gaussian noise. 
	Denote $\xi:=(\zeta_1,\dots,\zeta_{n},\omega)$ as our $n_0$-dimensional vector of uncertain parameters.
	To recover the unknown coefficients $c_{st}$ using samples $\hat{\xi}^{(1)},\dots,\hat{\xi}^{(N)}$ of $\xi$, 
	let $\Delta(s,t):=1+t+s(s+1)/2$ and we define our decision vector as 
	$x=(x_{\Delta(s,t)})_{0\le t\le s\le n}\in[-1,1]^{n_1}$ for $n_1=n_0(n_0+1)/2$.
	Consider the absolute deviation as our cost function
	\begin{equation}\label{eq:Regression}
		\begin{aligned}
			F(x,\xi)&=\Big\lvert\omega-\sum_{s=0}^{n}\sum_{t=0}^{s}x_{\Delta(s,t)}\zeta_s\zeta_t\Big\rvert=\max\{F_1,F_2\},\\
			\text{where}\ F_1(x,\xi)&=-F_2(x,\xi)=\omega-\sum_{s=0}^{n}\sum_{t=0}^{s}x_{\Delta(s,t)}\zeta_s\zeta_t.
		\end{aligned}
	\end{equation}
	In particular, setting $f(x)=0$, the ESO problem~\eqref{eq:ESO} of minimizing the expectation of $F(x,\xi)$ with respect to $\hat{\xi}^{(i)}=(\hat{\zeta}^{(i)}_1,\dots,\hat{\zeta}^{(i)}_n,\hat{\omega}^{(i)})$ for $i\in[N]$ 
	is a box-constrained least absolute deviations problem
	\begin{equation}\label{eq:RegressionESO}
		\min_{x\in[-1,1]^{n_1}}\quad\frac{1}{N}\sum_{i=1}^{N}\Big\lvert\hat{\omega}^{(i)}-\sum_{s=0}^{n}
		\sum_{t=0}^{s}x_{\Delta(s,t)}\hat{\zeta}_s^{(i)}\hat{\zeta}_t^{(i)}\Big\rvert.
	\end{equation}
	While the ESO problem~\eqref{eq:RegressionESO} is a convex program amenable to efficient algorithms, 
	its $p$-WDRO counterpart with $\Xi:=\{\xi=(\zeta,\omega):\zeta\in\bbR^n_{\ge0},\omega\in\bbR\}$
	\begin{equation}\label{eq:RegressionDRO}
		\min_{\substack{x\in[-1,1]^{n_1}\\\lambda\in\bbR_{\ge0}}}\ \frac{1}{N}\sum_{i=1}^{N}\sup_{\xi\in\Xi}
		\Big[\big\lvert\omega-\sum_{s=0}^{n}\sum_{t=0}^{s}x_{\Delta(s,t)}\zeta_s\zeta_t\big\rvert-\lambda\nVert{\xi-\hat{\xi}^{(i)}}^p\Big]
	\end{equation}
	can be much more challenging to solve, due to the lack of concavity of $F(x,\xi)$ in $\xi$.
	Nevertheless, $F(x,\xi)$ satisfies the condition in Proposition~\ref{prop:ConvexVx} 
	so we can apply our moment relaxation~\eqref{eq:1stage} and obtain a convex optimization problem in $x$ and $\lambda$ 
	which can be solved via first-order methods with distributable evaluation of the subgradients (see Algorithm~\ref{alg:Subgradient}).
	
\begin{figure}[h]
		\centering
		\begin{subfigure}{0.49\textwidth}
			\centering
			\includegraphics[width=\textwidth]{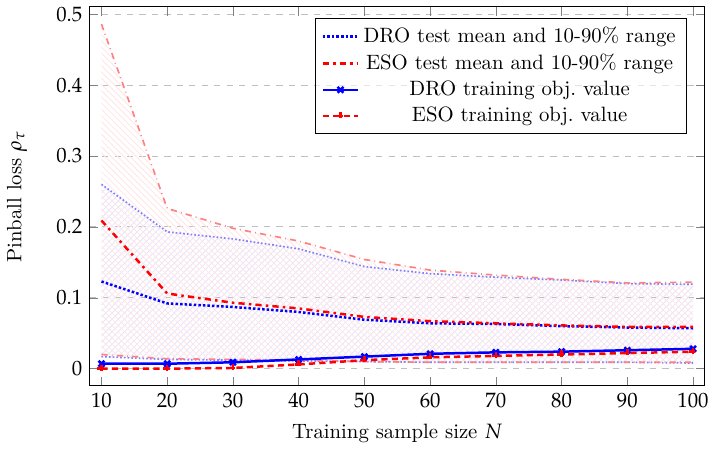}
			\caption{$r_0=0.01$}
			\label{fig:regression1}
		\end{subfigure}
		\hfill
		\begin{subfigure}{0.49\textwidth}
			\centering
			\includegraphics[width=\textwidth]{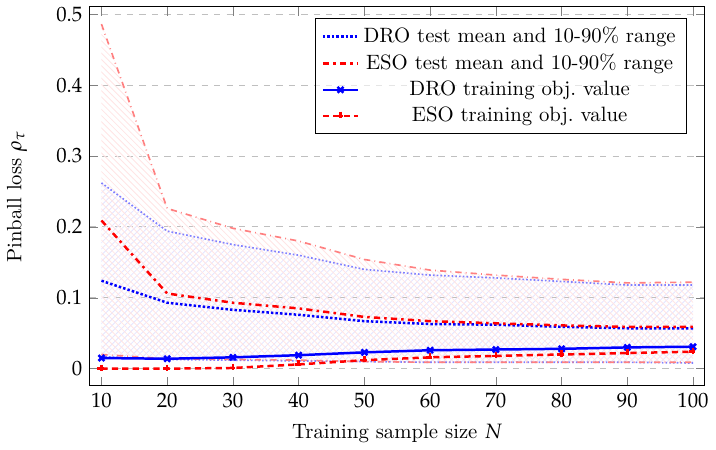}
			\caption{$r_0=0.02$}
			\label{fig:regression2}
		\end{subfigure}\\
		\vspace{2mm}
		\begin{subfigure}{0.49\textwidth}
			\centering
			\includegraphics[width=\textwidth]{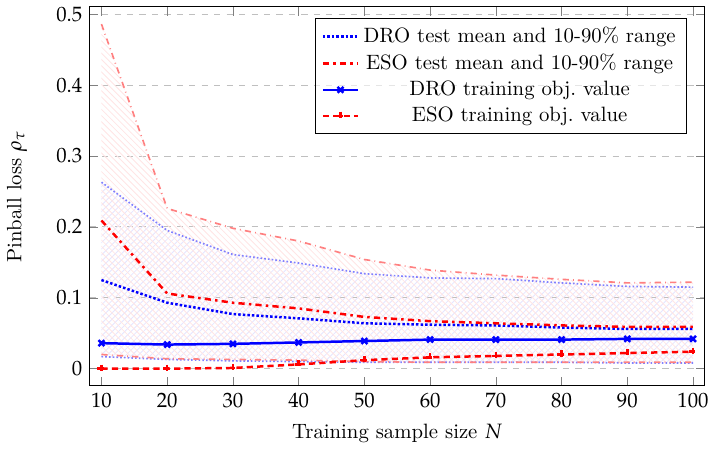}
			\caption{$r_0=0.05$}
			\label{fig:regression3}
		\end{subfigure}
		\hfill
		\begin{subfigure}{0.49\textwidth}
			\centering
			\includegraphics[width=\textwidth]{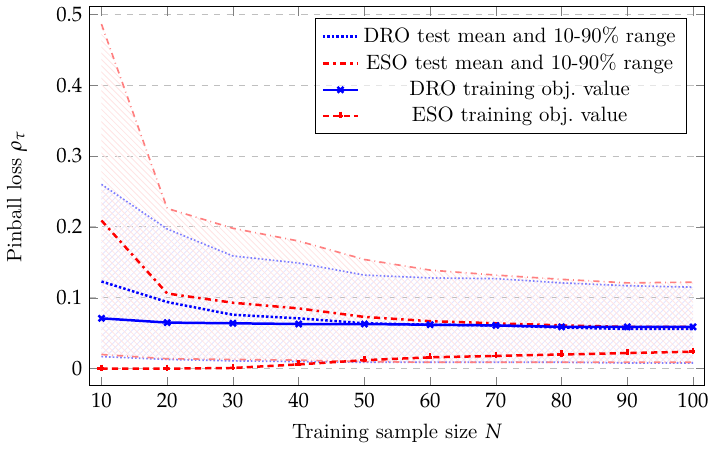}
			\caption{$r_0=0.1$}
			\label{fig:regression4}
		\end{subfigure}
		\caption{Training objective values and test means for box-constrained 
			least absolute deviations regression~\eqref{eq:Regression} with Wasserstein radii $r=r_0(N/10)^{-1/2}$}
		\label{fig:regression}
	\end{figure}

	We choose the number of independent variables to be $n=n_0-1=10$ in our experiments.
	Assume the standard deviation of $\epsilon$ is $\sigma=0.1$ and $(\zeta_1,\dots,\zeta_n)$ 
	is a multivariate folded normal random vector with $0$ mean and randomly generated covariance matrix of Frobenius norm $1$.
	Moreover, since the cost function $F$ is piecewise smooth, we follow~\cite{an2021generalization} 
	and set $r(N)=r_0(N/10)^{-1/2}$ for some initial radius $r_0>0$ while ensuring out-of-sample test performance guarantee with high probability.
	We report the experiment results in Table~\ref{tab:regression} and plot the training objective values and test means against training sample sizes in Figure~\ref{fig:regression}, from which we make the following comments.
	\begin{itemize}
		\item Our moment relaxation~\eqref{eq:1stage} scales well to larger numbers of training samples (e.g., $N=100$).
		Moreover, the gaps among the relaxed DRO training objective value, the ESO training objective value, and their corresponding out-of-sample test means,
		diminish as $N$ grows.
		\item Further, our moment relaxation~\eqref{eq:1stage} prevents overfitting as the original $p$-WDRO, meaning that it significantly reduces the out-of-sample test means and standard deviations when the training sample size is small (e.g., $N=10$).
	\end{itemize}
	Therefore, our moment relaxation~\eqref{eq:1stage} can be employed as a more tractable surrogate than the original WDRO on the constrained regression problems.

	\begin{longtable}{rrrrrr}
		\toprule
		$N$ & $r$ & Time (s) & Training Obj. & Test Mean & Test Std. \\
		\midrule
		\endfirsthead
		\toprule
		$N$ & $r$ & Time (s) & Training Obj. & Test Mean & Test Std. \\
		\midrule
		\endhead
		\midrule
		\multicolumn{3}{r}{\it Continued on next page} \\
		\endfoot
		\endlastfoot
		10 & 0.000 & 2.5 & $0.000 \pm 0.000$ & $0.209 \pm 0.087$ & $0.297 \pm 0.193$ \\
		10 & 0.010 & 146.7 & $0.007 \pm 0.001$ & $0.123 \pm 0.017$ & $0.120 \pm 0.032$ \\
		10 & 0.020 & 55.4 & $0.015 \pm 0.002$ & $0.124 \pm 0.017$ & $0.122 \pm 0.033$ \\
		10 & 0.050 & 68.8 & $0.036 \pm 0.006$ & $0.125 \pm 0.016$ & $0.120 \pm 0.030$ \\
		10 & 0.100 & 158.8 & $0.071 \pm 0.010$ & $0.123 \pm 0.018$ & $0.118 \pm 0.029$ \\
		\midrule
		20 & 0.000 & 7.3 & $0.000 \pm 0.000$ & $0.106 \pm 0.021$ & $0.119 \pm 0.045$ \\
		20 & 0.007 & 96.7 & $0.007 \pm 0.001$ & $0.092 \pm 0.008$ & $0.090 \pm 0.016$ \\
		20 & 0.014 & 52.5 & $0.014 \pm 0.002$ & $0.093 \pm 0.008$ & $0.090 \pm 0.016$ \\
		20 & 0.035 & 65.8 & $0.034 \pm 0.004$ & $0.093 \pm 0.007$ & $0.092 \pm 0.018$ \\
		20 & 0.071 & 72.6 & $0.065 \pm 0.006$ & $0.094 \pm 0.009$ & $0.092 \pm 0.021$ \\
		\midrule
		30 & 0.000 & 8.2 & $0.001 \pm 0.001$ & $0.093 \pm 0.008$ & $0.091 \pm 0.016$ \\
		30 & 0.006 & 60.8 & $0.009 \pm 0.002$ & $0.087 \pm 0.007$ & $0.084 \pm 0.012$ \\
		30 & 0.012 & 56.6 & $0.016 \pm 0.003$ & $0.083 \pm 0.007$ & $0.078 \pm 0.011$ \\
		30 & 0.029 & 56.4 & $0.035 \pm 0.005$ & $0.077 \pm 0.006$ & $0.070 \pm 0.008$ \\
		30 & 0.058 & 82.6 & $0.064 \pm 0.008$ & $0.076 \pm 0.006$ & $0.068 \pm 0.008$ \\
		\midrule
		40 & 0.000 & 12.8 & $0.006 \pm 0.002$ & $0.085 \pm 0.012$ & $0.087 \pm 0.025$ \\
		40 & 0.005 & 152.5 & $0.013 \pm 0.002$ & $0.080 \pm 0.010$ & $0.078 \pm 0.020$ \\
		40 & 0.010 & 75.4 & $0.019 \pm 0.003$ & $0.076 \pm 0.011$ & $0.073 \pm 0.019$ \\
		40 & 0.025 & 97.8 & $0.037 \pm 0.004$ & $0.071 \pm 0.008$ & $0.066 \pm 0.013$ \\
		40 & 0.050 & 106.1 & $0.063 \pm 0.006$ & $0.071 \pm 0.007$ & $0.065 \pm 0.010$ \\
		\midrule
		50 & 0.000 & 20.5 & $0.012 \pm 0.003$ & $0.073 \pm 0.010$ & $0.069 \pm 0.017$ \\
		50 & 0.004 & 161.0 & $0.017 \pm 0.003$ & $0.069 \pm 0.006$ & $0.062 \pm 0.011$ \\
		50 & 0.009 & 85.8 & $0.023 \pm 0.003$ & $0.067 \pm 0.006$ & $0.059 \pm 0.010$ \\
		50 & 0.022 & 120.7 & $0.039 \pm 0.004$ & $0.064 \pm 0.005$ & $0.056 \pm 0.009$ \\
		50 & 0.045 & 144.3 & $0.063 \pm 0.005$ & $0.064 \pm 0.005$ & $0.056 \pm 0.007$ \\
		\midrule
		60 & 0.000 & 26.8 & $0.016 \pm 0.004$ & $0.067 \pm 0.007$ & $0.060 \pm 0.010$ \\
		60 & 0.004 & 192.1 & $0.021 \pm 0.004$ & $0.064 \pm 0.006$ & $0.057 \pm 0.009$ \\
		60 & 0.008 & 93.2 & $0.026 \pm 0.005$ & $0.063 \pm 0.006$ & $0.055 \pm 0.008$ \\
		60 & 0.020 & 146.8 & $0.041 \pm 0.005$ & $0.062 \pm 0.005$ & $0.053 \pm 0.004$ \\
		60 & 0.041 & 122.3 & $0.062 \pm 0.006$ & $0.062 \pm 0.005$ & $0.053 \pm 0.004$ \\
		\midrule
		70 & 0.000 & 33.5 & $0.018 \pm 0.004$ & $0.064 \pm 0.005$ & $0.056 \pm 0.009$ \\
		70 & 0.004 & 151.3 & $0.023 \pm 0.004$ & $0.063 \pm 0.004$ & $0.055 \pm 0.008$ \\
		70 & 0.008 & 86.5 & $0.027 \pm 0.005$ & $0.062 \pm 0.004$ & $0.054 \pm 0.007$ \\
		70 & 0.019 & 112.6 & $0.041 \pm 0.005$ & $0.061 \pm 0.005$ & $0.053 \pm 0.007$ \\
		70 & 0.038 & 123.6 & $0.061 \pm 0.005$ & $0.061 \pm 0.005$ & $0.052 \pm 0.007$ \\
		\midrule
		80 & 0.000 & 40.8 & $0.020 \pm 0.003$ & $0.061 \pm 0.005$ & $0.053 \pm 0.009$ \\
		80 & 0.004 & 131.0 & $0.024 \pm 0.004$ & $0.060 \pm 0.004$ & $0.052 \pm 0.007$ \\
		80 & 0.007 & 99.9 & $0.028 \pm 0.004$ & $0.059 \pm 0.004$ & $0.051 \pm 0.007$ \\
		80 & 0.018 & 130.3 & $0.041 \pm 0.004$ & $0.058 \pm 0.004$ & $0.050 \pm 0.007$ \\
		80 & 0.035 & 138.0 & $0.059 \pm 0.004$ & $0.058 \pm 0.004$ & $0.050 \pm 0.007$ \\
		\midrule
		90 & 0.000 & 65.1 & $0.022 \pm 0.003$ & $0.059 \pm 0.006$ & $0.051 \pm 0.010$ \\
		90 & 0.003 & 154.9 & $0.026 \pm 0.003$ & $0.058 \pm 0.005$ & $0.050 \pm 0.008$ \\
		90 & 0.007 & 108.6 & $0.030 \pm 0.003$ & $0.057 \pm 0.005$ & $0.049 \pm 0.008$ \\
		90 & 0.017 & 131.7 & $0.042 \pm 0.004$ & $0.056 \pm 0.005$ & $0.048 \pm 0.009$ \\
		90 & 0.033 & 151.6 & $0.059 \pm 0.004$ & $0.056 \pm 0.004$ & $0.048 \pm 0.007$ \\
		\midrule
		100 & 0.000 & 57.8 & $0.024 \pm 0.002$ & $0.059 \pm 0.006$ & $0.051 \pm 0.009$ \\
		100 & 0.003 & 278.2 & $0.028 \pm 0.003$ & $0.057 \pm 0.005$ & $0.048 \pm 0.007$ \\
		100 & 0.006 & 162.1 & $0.031 \pm 0.003$ & $0.057 \pm 0.005$ & $0.048 \pm 0.006$ \\
		100 & 0.016 & 195.7 & $0.042 \pm 0.003$ & $0.056 \pm 0.004$ & $0.046 \pm 0.006$ \\
		100 & 0.032 & 209.0 & $0.059 \pm 0.003$ & $0.056 \pm 0.004$ & $0.046 \pm 0.005$ \\
		\bottomrule
		\caption{Experiment results for the box-constrained least absolute deviations regression~\eqref{eq:Regression} over 10 independent replications}
		\label{tab:regression}
	\end{longtable}
	
	To quantitatively see the relaxation gap of our approach~\eqref{eq:1stage}, we solve alternatively the inner supremum of~\eqref{eq:RegressionDRO} as a nonconvex QP using \texttt{PolyJuMP} (provided by the \texttt{SumOfSquares} package).
	In particular, we set the nonconvex option of Gurobi to be \texttt{NonConvex=2} and set the time limit to be 600 seconds.
	The result of a typical run (without replication) is presented in Table~\ref{tab:NonconvexQPComparison}.
	From the table, we see that all of our moment relaxations~\eqref{eq:1stage} can be solved consistently within the time limit, while half of the nonconvex QP runs fail to with at optimality.
	On all the instances on which the nonconvex QP can be solved within the time limit, the difference of the training objectives and test means and standard deviations from our moment relaxation is negligible (no larger than $0.002$).
	These observations further corroborate the high accuracy of the solutions obtained by our moment relaxation for the original $p$-WDRO problems.

\begin{table}[h!]
	\centering
	\scriptsize
	\begin{tabular}{rrrrrrrrrr}
		\toprule
		& & \multicolumn{4}{c}{Moment-WDRO} & \multicolumn{4}{c}{Nonconvex-QP} \\
		\cmidrule(lr){3-6} \cmidrule(lr){7-10}
		& & \multicolumn{2}{c}{Training} & \multicolumn{2}{c}{Test} & \multicolumn{2}{c}{Training} & \multicolumn{2}{c}{Test}\\
		\cmidrule(lr){3-4} \cmidrule(lr){5-6} \cmidrule(lr){7-8} \cmidrule(lr){9-10}
		$N$ & $r$ & Time (s) & Obj. & Mean & Std.
		& Time (s) & Obj. & Mean & Std. \\
		\midrule
		10 & 0.000 & 8.4 & $0.000$ & $0.196$ & $0.168$ & 1.5 & $0.000$ & $0.196$ & $0.168$ \\
		10 & 0.010 & 47.7 & $0.012$ & $0.145$ & $0.165$ & 34.4 & $0.012$ & $0.146$ & $0.165$ \\
		10 & 0.020 & 34.1 & $0.024$ & $0.146$ & $0.164$ & 138.2 & $0.024$ & $0.146$ & $0.164$ \\
		10 & 0.050 & 32.7 & $0.061$ & $0.145$ & $0.163$ & 134.6 & $0.061$ & $0.145$ & $0.163$ \\
		10 & 0.100 & 61.9 & $0.121$ & $0.146$ & $0.165$ & $>600$ & $1.027$ & $0.482$ & $0.567$ \\
		\midrule
		20 & 0.000 & 3.7 & $0.000$ & $0.137$ & $0.111$ & 3.3 & $0.000$ & $0.137$ & $0.111$ \\
		20 & 0.007 & 134.1 & $0.010$ & $0.125$ & $0.104$ & 38.4 & $0.010$ & $0.127$ & $0.108$ \\
		20 & 0.014 & 23.5 & $0.020$ & $0.125$ & $0.105$ & 578.0 & $0.020$ & $0.125$ & $0.104$ \\
		20 & 0.035 & 40.8 & $0.049$ & $0.126$ & $0.110$ & $>600$ & $154.8$ & $4.803$ & $4.222$ \\
		20 & 0.071 & 48.7 & $0.094$ & $0.123$ & $0.130$ & $>600$ & $0.219$ & $0.127$ & $0.125$ \\
		\midrule
		30 & 0.000 & 5.7 & $0.003$ & $0.101$ & $0.084$ & 5.4 & $0.003$ & $0.101$ & $0.084$ \\
		30 & 0.006 & 45.0 & $0.012$ & $0.099$ & $0.080$ & 60.6 & $0.012$ & $0.099$ & $0.080$ \\
		30 & 0.012 & 65.7 & $0.022$ & $0.096$ & $0.075$ & $>600$ & $9.213$ & $4.798$ & $4.217$ \\
		30 & 0.029 & 59.8 & $0.047$ & $0.095$ & $0.077$ & 563.9 & $0.048$ & $0.095$ & $0.078$ \\
		30 & 0.058 & 107.5 & $0.088$ & $0.088$ & $0.075$ & 542.4 & $0.088$ & $0.089$ & $0.076$ \\
		\midrule
		40 & 0.000 & 12.1 & $0.006$ & $0.084$ & $0.068$ & 12.1 & $0.006$ & $0.084$ & $0.068$ \\
		40 & 0.005 & 57.0 & $0.015$ & $0.077$ & $0.063$ & 67.2 & $0.015$ & $0.077$ & $0.063$ \\
		40 & 0.010 & 51.3 & $0.023$ & $0.072$ & $0.058$ & $>600$ & $165.5$ & $4.803$ & $4.222$ \\
		40 & 0.025 & 64.5 & $0.047$ & $0.066$ & $0.052$ & $>600$ & $822.8$ & $4.803$ & $4.222$ \\
		40 & 0.050 & 64.1 & $0.084$ & $0.064$ & $0.050$ & $>600$ & $0.123$ & $0.063$ & $0.048$ \\
		\midrule
		50 & 0.000 & 20.1 & $0.010$ & $0.066$ & $0.051$ & 24.4 & $0.010$ & $0.066$ & $0.051$ \\
		50 & 0.004 & 186.8 & $0.017$ & $0.066$ & $0.052$ & 98.2 & $0.017$ & $0.066$ & $0.052$ \\
		50 & 0.009 & 58.1 & $0.024$ & $0.066$ & $0.052$ & $>600$ & $49.5$ & $4.803$ & $4.222$ \\
		50 & 0.022 & 76.9 & $0.045$ & $0.063$ & $0.051$ & $>600$ & $167.7$ & $4.803$ & $4.222$ \\
		50 & 0.045 & 84.5 & $0.078$ & $0.062$ & $0.049$ & $>600$ & $939.6$ & $4.803$ & $4.222$ \\
		\midrule
		60 & 0.000 & 33.7 & $0.012$ & $0.061$ & $0.050$ & 33.3 & $0.012$ & $0.061$ & $0.050$ \\
		60 & 0.004 & 102.0 & $0.019$ & $0.060$ & $0.049$ & 93.4 & $0.019$ & $0.061$ & $0.049$ \\
		60 & 0.008 & 78.4 & $0.025$ & $0.060$ & $0.048$ & $>600$ & $165.2$ & $4.803$ & $4.222$ \\
		60 & 0.020 & 303.6 & $0.043$ & $0.059$ & $0.048$ & $>600$ & $165.2$ & $4.803$ & $4.222$ \\
		60 & 0.041 & 204.3 & $0.073$ & $0.063$ & $0.050$ & $>600$ & $165.2$ & $4.803$ & $4.222$ \\
		\midrule
		70 & 0.000 & 42.9 & $0.015$ & $0.058$ & $0.046$ & 38.4 & $0.015$ & $0.058$ & $0.046$ \\
		70 & 0.004 & 106.2 & $0.021$ & $0.057$ & $0.045$ & 99.1 & $0.021$ & $0.057$ & $0.045$ \\
		70 & 0.008 & 111.0 & $0.027$ & $0.057$ & $0.045$ & $>600$ & $157.8$ & $4.803$ & $4.222$ \\
		70 & 0.019 & 256.4 & $0.044$ & $0.057$ & $0.045$ & $>600$ & $157.8$ & $4.803$ & $4.222$ \\
		70 & 0.038 & 131.9 & $0.072$ & $0.058$ & $0.045$ & $>600$ & $157.8$ & $4.803$ & $4.222$ \\
		\midrule
		80 & 0.000 & 47.8 & $0.019$ & $0.054$ & $0.043$ & 53.7 & $0.019$ & $0.054$ & $0.043$ \\
		80 & 0.004 & 118.0 & $0.024$ & $0.054$ & $0.043$ & 144.5 & $0.024$ & $0.055$ & $0.043$ \\
		80 & 0.007 & 173.3 & $0.030$ & $0.055$ & $0.043$ & $>600$ & $154.3$ & $4.803$ & $4.222$ \\
		80 & 0.018 & 105.9 & $0.045$ & $0.055$ & $0.043$ & $>600$ & $154.3$ & $4.803$ & $4.222$ \\
		80 & 0.035 & 121.9 & $0.071$ & $0.053$ & $0.041$ & $>600$ & $154.3$ & $4.803$ & $4.222$ \\
		\midrule
		90 & 0.000 & 58.3 & $0.020$ & $0.055$ & $0.044$ & 62.9 & $0.020$ & $0.055$ & $0.044$ \\
		90 & 0.003 & 330.4 & $0.025$ & $0.055$ & $0.043$ & 142.1 & $0.025$ & $0.055$ & $0.043$ \\
		90 & 0.007 & 134.3 & $0.030$ & $0.055$ & $0.043$ & $>600$ & $147.5$ & $4.803$ & $4.222$ \\
		90 & 0.017 & 117.8 & $0.045$ & $0.053$ & $0.041$ & $>600$ & $147.5$ & $4.803$ & $4.222$ \\
		90 & 0.033 & 144.4 & $0.069$ & $0.052$ & $0.040$ & $>600$ & $147.5$ & $4.803$ & $4.222$ \\
		\midrule
		100 & 0.000 & 218.7 & $0.020$ & $0.053$ & $0.042$ & 58.9 & $0.020$ & $0.053$ & $0.042$ \\
		100 & 0.003 & 197.6 & $0.025$ & $0.052$ & $0.041$ & 206.2 & $0.025$ & $0.052$ & $0.041$ \\
		100 & 0.006 & 122.5 & $0.030$ & $0.052$ & $0.041$ & $>600$ & $139.4$ & $4.803$ & $4.222$ \\
		100 & 0.016 & 283.3 & $0.043$ & $0.051$ & $0.040$ & $>600$ & $139.4$ & $4.803$ & $4.222$ \\
		100 & 0.032 & 191.4 & $0.065$ & $0.050$ & $0.039$ & $>600$ & $139.4$ & $4.803$ & $4.222$ \\
		\bottomrule
	\end{tabular}
	\caption{Comparison with nonconvex QP formulation}
	\label{tab:NonconvexQPComparison}
\end{table}

\subsection{Two-stage production problem}
We now turn to a two-stage production problem, which is adapted from~\cite[Chapter 1.3.1]{shapiro2021lectures}.
In such a production problem, the decision-maker needs to determine the inventory of $n_1$ ingredients $x=(x_1,\dots,x_{n_1})\in[0,D]^{n_1}$, each with a unit cost $c^\rmi_t$ for $t\in[n_1]$, in the first stage, where $D>0$ is the capacity limit.
The second-stage cost is described by the value function~\eqref{eq:Recourse} as
\begin{equation}\label{eq:ExperimentProduction}
	\begin{aligned}
		F(x,\xi)=\min_{x'=(x^\rmp,x^\rmq,x^\rmr)}\ & -\sum_{s=1}^{n_\rmp}c^\rmp x^{\rmp}_s-\sum_{t=1}^{n_1}c^\rmq_t(\xi)x^{\rmq}_t+\sum_{t=1}^{n_1}c^\rmr_tx^\rmr_t \\
		\suchthat\quad\quad & x^\rmq_t-x^{\rmr}_t=B_{t}^\rmp(\xi)x_t-\sum_{s=1}^{n_\rmp}A_{ts}^\rmm x^{\rmp}_s, \quad && \forall t\in[n_1],\\
		& 0\le x^\rmp_s \le b_s(\xi), && \forall s\in[n_r],\\
		& x^\rmq_t,x^\rmr_t\ge0, && \forall t\in[n_1].
	\end{aligned}
\end{equation}
Here, $n_\rmp$ is the number of products to be produced; 
$x^\rmp=(x^\rmp_1,\dots,x^\rmp_{n_\rmp})$ is the decision vector on quantities of each of the products to produce, which can be sold at the price of $c^\rmp_1,\dots,c^\rmp_{n_\rmp}$;
$x^\rmq=(x^\rmq_1,\dots,x^\rmq_{n_1})$ is the vector of the quantities of ingredients that are left unused, each of which can be sold at a random salvage price $c^\rmq_t(\xi)$;
$x^\rmr=(x^\rmr_1,\dots,x^\rmr_{n_\rmp})$ is the vector of ingredient purchases with late prices $c^\rmr_1,\dots,c^\rmr_{n_1}$;
$B^\rmp_1(\xi),\dots,B^\rmp_{n_1}(\xi)$ are the random decay percentage of the stored ingredients;
$A^\rmm_{ts}$ is the quantity of ingredient $t$ needed to produce a unit of product $s$, for each $s\in[n_\rmp]$ and $t\in[n_1]$;
and $b_1(\xi),\dots,b_{n_\rmp}(\xi)$ are the random demands for the products.
It is clear that the second-stage problem~\eqref{eq:ExperimentProduction} is always feasible by setting $x^\rmp=0$, $x^\rmr=0$ and $x^\rmq=x$, regardless of the outcomes of the random parameters.
Using linear optimization duality, we reformulate the value function as a maximization to facilitate building our moment relaxation~\eqref{eq:2stage}
\begin{equation}\label{eq:ExperimentProductionDual}
	\begin{aligned}
		F(x,\xi)=\max_{u\in\bbR^{n_1+n_\rmp}}\quad & -(B^\rmp(\xi)x,b(\xi))^\transpose u \\
		\suchthat\quad\quad & \begin{bmatrix}
			I & 0 \\
			-I & 0 \\
			A^\rmm & I \\
			0 & I \\
			0 & -I
		\end{bmatrix} u \ge
		\begin{bmatrix}
			c^\rmq(\xi) \\
			-c^\rmr \\
			c^\rmp \\
			0 \\
			-c^\rmp 
		\end{bmatrix},
	\end{aligned}
\end{equation}
where $A^\rmm=(A^\rmm_{ts})_{t\in[n_1],s\in[n_\rmp]}$ is an $n_1\times n_\rmp$ matrix, 
and $B^\rmp(\xi)$ is the diagonal matrix consisting of $B^\rmp_1(\xi),\dots,B^\rmp_{n_1}(\xi)$.
It is worth noting that the original dual formulation of~\eqref{eq:ExperimentProduction} 
has an unbounded feasibility set for $u$.
We added the last row of the linear constraints 
	without affecting the value $F(x,\xi)$ or the $p$-WDRO problem, 
because whenever $b(\xi)>0$ the last $n_\rmp$ components of an optimal solution $u$ do not exceed $c^\rmp$, similar to the approach taken by~\cite{duque2022distributionally}.

We now describe the settings of the parameters in~\eqref{eq:ExperimentProduction} used for our experiment results below.
We consider the capacity limit $D=5$.
The ingredient purchase prices $c^\rmi_t=2+3(t-1)/(n_1-1)$, and the late purchase prices are $c^\rmr_t=3c^\rmi_t$, for each $t\in[n_1]$.
We set $A^\rmm_{ts}=1/(10\cdot(t-1))$ for each $s\le t-1$, and $A^\rmm_{ts}=9/(10\cdot(t-1))$ for each $s\ge t$.
For the uncertain parameters  $c^\rmq(\xi),B^\rmp(\xi)$, and $b(\xi)$, we define the underlying truth by specifying the probability distributions of each of these parameters.
The demands follow independently lognormal distributions $b_s(\xi)\sim \operatorname{LogNormal}(\log(\bar{b}_s),\sigma)$, where we set $\bar{b}_s=2-(s-1)/(n_\rmp-1)$ for each $s\in[n_\rmp]$, and $\sigma=0.1$.
The decay percentages and salvage prices depend on whether the ingredient is perishable or nonperishable.
Let $T\subseteq[n_1]$ be an index subset of perishable ingredients, which we set to be all of the odd indices.
Then $B^\rmp_t(\xi)\sim\operatorname{Uniform}(0,\bar{B})$ for each $t\in T$, where we $\bar{B}=1$, and $B^\rmp_t(\xi)\equiv1$ is deterministic for each $t\notin T$;
$c^\rmq_t(\xi)\sim\operatorname{Uniform}(0,\bar{c}^\rmq_t)$ for each $t\notin T$, and $c^\rmq_t(\xi)=\bar{c}^\rmq_t$ is deterministic for each $t\in T$, where we set $\bar{c}^\rmq_t=5-3(t-1)/(n_1-1)$ for each $t\in[n_1]$.
By abuse of notation, we may view $\xi$ as the concatenation of an $n_\rmp$-dimensional lognormal random vector, and a $n_1$-dimensional uniformly distributed random vector, so the support set $\Xi=[0,+\infty)^{n_\rmp}\times [0,1]^{n_1}\subset\bbR^{n_0}$ where $n_0=n_\rmp+n_1$.
	Thus, $\deg(B(\xi))=\deg(b(\xi))=1$ and $\deg(G_x)=2$, so $d_g=\max\{2\deg(B(\xi)),2\deg(b(\xi)),\deg(G_x),\deg(g)\}=2$ in~\eqref{eq:d_g}.
	We set $p=d_g=2$ which ensures that the asymptotic consistency is guaranteed by Theorem~\ref{thm:2stage_idc}.

\begin{figure}[h]
		\centering
		\begin{subfigure}{0.49\textwidth}
			\centering
			\includegraphics[width=\textwidth]{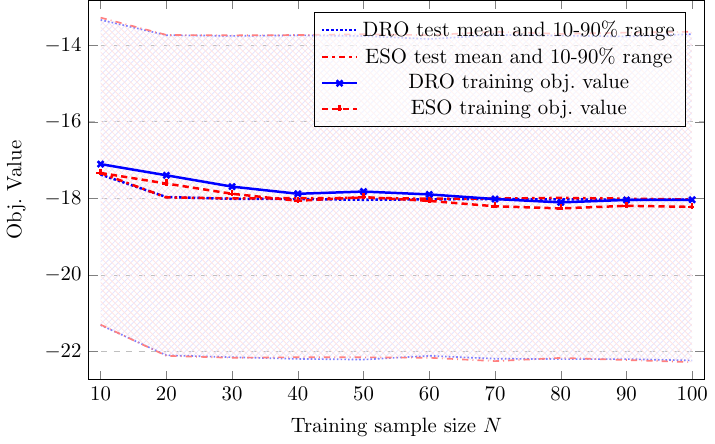}
			\caption{$r_0=0.01$}
			\label{fig:production1}
		\end{subfigure}
		\hfill
		\begin{subfigure}{0.49\textwidth}
			\centering
			\includegraphics[width=\textwidth]{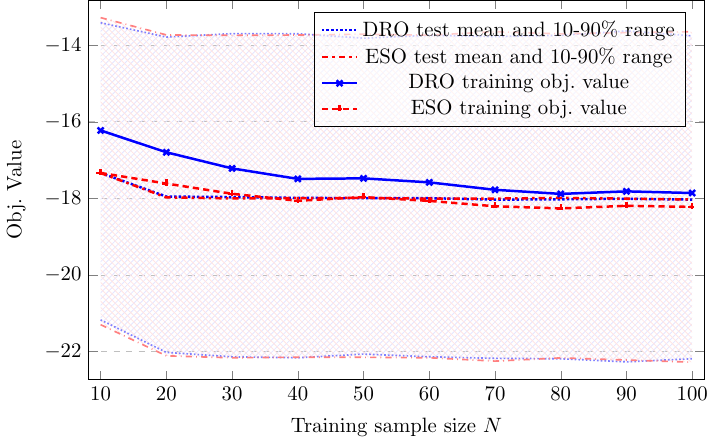}
			\caption{$r_0=0.05$}
			\label{fig:production2}
		\end{subfigure}\\
		\vspace{2mm}
		\begin{subfigure}{0.49\textwidth}
			\centering
			\includegraphics[width=\textwidth]{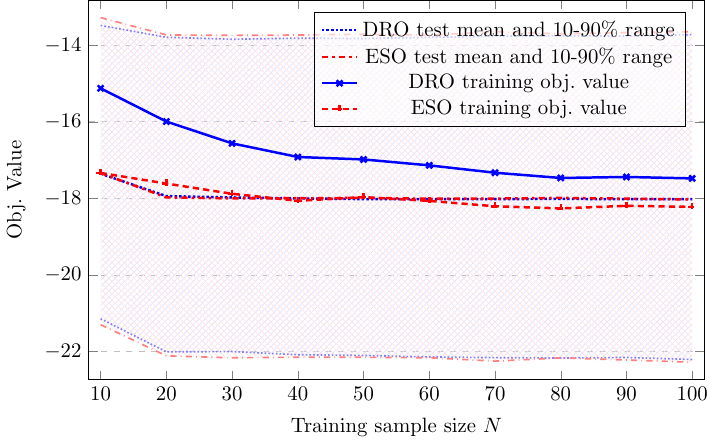}
			\caption{$r_0=0.1$}
			\label{fig:production3}
		\end{subfigure}
		\caption{Training objective values and test means for two-stage production problems~\eqref{eq:ExperimentProduction} with Wasserstein radii $r=r_0(N/10)^{-1/2}$}
		\label{fig:production}
	\end{figure}

	We report the experiment results for $n_1=20$ and $n_\rmp=10$.
	The Wasserstein radius function $r(N)$ is set heuristically to $r(N)=r_0(N/10)^{-1/2}$, despite the large dimension $n_0=n_\rmp+n_1=30$ of the random vector $\xi$.
	Such square-root decay scheme has been widely used for $p=2$~\cite{gao2023finite,blanchet2022confidence}, and we adopt it mainly for the purpose of easy illustration of our asymptotic consistency as $N$ grows (cf.~the slow decrease of the optimality gap in Table~\ref{tab:HKComparison}).
	We summarize the results over $N^{{\mathrm{rep}}}=3$ replications in Table~\ref{tab:production} and plot the training objective values and test means against training sample sizes in Figure~\ref{fig:production}, from which we have observed that the similar decreasing behavior of the gaps between relaxed DRO, ESO training objective values and their corresponding out-of-sample test means as we did in the constrained regression problems.
	Although being less significant, our moment relaxation~\eqref{eq:2stage} could also reduce the test mean and standard deviation on the two-stage production problem, when the training sample size is small $N=10$.
	
	\begin{table}[h]
		\scriptsize
		\centering
		\begin{tabular}{rrrrrr}
			\toprule
			$N$ & $r$ & Time (s) & Training Obj. & Test Mean & Test Std. \\
			\midrule
			10 & 0.000 & 8.4 & $-17.334 \pm 1.488$ & $-17.337 \pm 0.815$ & $3.091 \pm 0.481$ \\
			10 & 0.010 & 402.4 & $-17.103 \pm 1.480$ & $-17.376 \pm 0.754$ & $3.068 \pm 0.445$ \\
			10 & 0.050 & 279.9 & $-16.223 \pm 1.455$ & $-17.333 \pm 0.836$ & $2.990 \pm 0.453$ \\
			10 & 0.100 & 244.5 & $-15.123 \pm 1.410$ & $-17.353 \pm 0.820$ & $2.947 \pm 0.385$ \\
			\midrule
			20 & 0.000 & 10.7 & $-17.611 \pm 1.146$ & $-17.973 \pm 0.073$ & $3.225 \pm 0.220$ \\
			20 & 0.007 & 627.0 & $-17.395 \pm 1.119$ & $-17.964 \pm 0.037$ & $3.221 \pm 0.242$ \\
			20 & 0.035 & 424.0 & $-16.794 \pm 1.112$ & $-17.950 \pm 0.090$ & $3.174 \pm 0.233$ \\
			20 & 0.071 & 560.1 & $-15.991 \pm 1.087$ & $-17.939 \pm 0.091$ & $3.167 \pm 0.228$ \\
			\midrule
			30 & 0.000 & 13.8 & $-17.878 \pm 1.013$ & $-18.000 \pm 0.027$ & $3.251 \pm 0.177$ \\
			30 & 0.006 & 933.8 & $-17.691 \pm 1.003$ & $-18.009 \pm 0.014$ & $3.235 \pm 0.171$ \\
			30 & 0.029 & 487.8 & $-17.215 \pm 0.987$ & $-17.964 \pm 0.021$ & $3.248 \pm 0.252$ \\
			30 & 0.058 & 659.6 & $-16.561 \pm 1.000$ & $-17.974 \pm 0.056$ & $3.141 \pm 0.213$ \\
			\midrule
			40 & 0.000 & 19.8 & $-18.061 \pm 0.837$ & $-17.993 \pm 0.008$ & $3.239 \pm 0.113$ \\
			40 & 0.005 & 1035.2 & $-17.878 \pm 0.841$ & $-18.016 \pm 0.015$ & $3.257 \pm 0.133$ \\
			40 & 0.025 & 723.2 & $-17.489 \pm 0.837$ & $-17.985 \pm 0.015$ & $3.255 \pm 0.161$ \\
			40 & 0.050 & 766.6 & $-16.917 \pm 0.828$ & $-17.999 \pm 0.027$ & $3.186 \pm 0.143$ \\
			\midrule
			50 & 0.000 & 23.1 & $-17.965 \pm 0.816$ & $-17.985 \pm 0.013$ & $3.259 \pm 0.164$ \\
			50 & 0.004 & 1253.9 & $-17.820 \pm 0.806$ & $-18.035 \pm 0.010$ & $3.254 \pm 0.147$ \\
			50 & 0.022 & 972.2 & $-17.475 \pm 0.815$ & $-17.992 \pm 0.032$ & $3.176 \pm 0.096$ \\
			50 & 0.045 & 851.3 & $-16.981 \pm 0.806$ & $-18.017 \pm 0.035$ & $3.189 \pm 0.124$ \\
			\midrule
			60 & 0.000 & 27.7 & $-18.064 \pm 0.939$ & $-18.008 \pm 0.030$ & $3.246 \pm 0.170$ \\
			60 & 0.004 & 1512.1 & $-17.895 \pm 0.912$ & $-18.022 \pm 0.048$ & $3.191 \pm 0.165$ \\
			60 & 0.020 & 1033.7 & $-17.580 \pm 0.907$ & $-17.992 \pm 0.012$ & $3.240 \pm 0.185$ \\
			60 & 0.041 & 1057.2 & $-17.136 \pm 0.903$ & $-18.019 \pm 0.040$ & $3.215 \pm 0.155$ \\
			\midrule
			70 & 0.000 & 31.6 & $-18.204 \pm 1.000$ & $-18.004 \pm 0.038$ & $3.308 \pm 0.181$ \\
			70 & 0.004 & 2661.6 & $-18.015 \pm 0.981$ & $-18.019 \pm 0.023$ & $3.255 \pm 0.114$ \\
			70 & 0.019 & 1308.5 & $-17.774 \pm 0.955$ & $-18.031 \pm 0.007$ & $3.236 \pm 0.130$ \\
			70 & 0.038 & 1008.0 & $-17.328 \pm 0.973$ & $-18.017 \pm 0.025$ & $3.237 \pm 0.126$ \\
			\midrule
			80 & 0.000 & 39.8 & $-18.260 \pm 0.876$ & $-17.985 \pm 0.026$ & $3.261 \pm 0.171$ \\
			80 & 0.004 & 3656.7 & $-18.103 \pm 0.861$ & $-18.022 \pm 0.025$ & $3.254 \pm 0.111$ \\
			80 & 0.018 & 1801.9 & $-17.881 \pm 0.863$ & $-18.022 \pm 0.006$ & $3.249 \pm 0.162$ \\
			80 & 0.035 & 1086.4 & $-17.465 \pm 0.870$ & $-18.012 \pm 0.030$ & $3.242 \pm 0.145$ \\
			\midrule
			90 & 0.000 & 64.7 & $-18.192 \pm 0.817$ & $-18.005 \pm 0.010$ & $3.299 \pm 0.132$ \\
			90 & 0.003 & 3923.4 & $-18.038 \pm 0.805$ & $-18.036 \pm 0.020$ & $3.251 \pm 0.134$ \\
			90 & 0.017 & 2550.9 & $-17.815 \pm 0.812$ & $-18.008 \pm 0.013$ & $3.324 \pm 0.130$ \\
			90 & 0.033 & 1187.0 & $-17.439 \pm 0.816$ & $-18.018 \pm 0.037$ & $3.240 \pm 0.180$ \\
			\midrule
			100 & 0.000 & 47.5 & $-18.220 \pm 0.731$ & $-18.028 \pm 0.018$ & $3.333 \pm 0.108$ \\
			100 & 0.003 & 2450.4 & $-18.034 \pm 0.709$ & $-18.027 \pm 0.005$ & $3.285 \pm 0.093$ \\
			100 & 0.016 & 2028.4 & $-17.858 \pm 0.710$ & $-18.030 \pm 0.021$ & $3.252 \pm 0.133$ \\
			100 & 0.032 & 1246.5 & $-17.476 \pm 0.713$ & $-18.019 \pm 0.003$ & $3.276 \pm 0.096$ \\
			\bottomrule
		\end{tabular}
		\caption{Experiment results for the two-stage production problem~\eqref{eq:ExperimentProduction} over 3 independent replications}
		\label{tab:production}
	\end{table}
	
\begin{figure}
		\centering
		\includegraphics[width=0.6\linewidth]{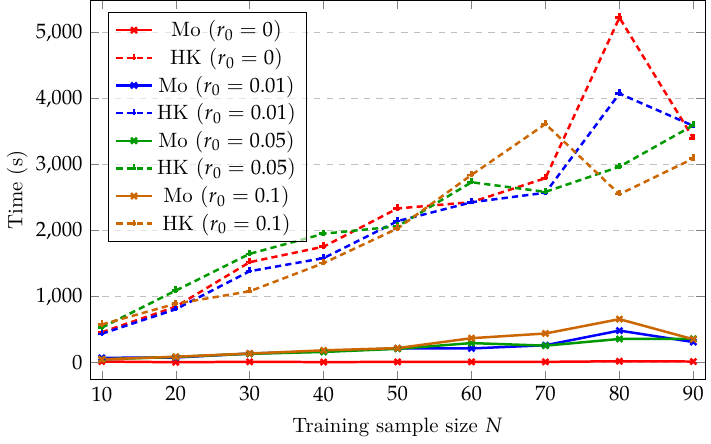}
		\caption{Computational time comparison against Hanasusanto-Kuhn semidefinite relaxation (Mo for our moment relaxation~\eqref{eq:2stage} and HK for the HK relaxation)}
		\label{fig:HKComparison}
	\end{figure}

	To further demonstrate the scalability of our moment relaxation for two-stage $2$-WDRO, we compare it against the copositive reformulation-based semidefinite relaxation proposed in~\cite{hanasusanto2018conic}, shortened as the HK relaxation below.
	Since the HK relaxation is coupled across $N$ training samples, we now use a single CPU with 100 GB of RAM for fairer comparison.
	Moreover, we reduce the size of the two-stage production problem to $n_1=10$ and $n_\rmp=5$, for which the HK relaxation can be solved with up to $N=90$ training samples before we run out of memory.
	Besides, to ensure that the HK relaxation provides a valid overestimation, we set the parameter $\delta$ in~\cite[Theorem 5]{hanasusanto2018conic} to 0 and report the best feasible objective value returned by \texttt{Mosek} whenever \texttt{slow progress} is encountered.
	To save computational resources, we do not replicate the experiment runs.
	The Wasserstein radius is set to be $r(N)=r_0(N/10)^{-1/15}$ in this case.
	The results are summarized in Table~\ref{tab:HKComparison} and the computational times are plotted in Figure~\ref{fig:HKComparison} for easy comparison, from which we have the following remarks.
	\begin{itemize}
		\item In terms of computational times, our moment relaxation~\eqref{eq:2stage} scales consistently better than the HK relaxation as $N$ grows, for all choices of radii $r_0$.
		Since we disabled the distributed subgradient evaluation, we believe that the better scalability is mainly due to the independence of the SDP sizes in our moment relaxations from the training sample size $N$.
		\item For larger Wasserstein radii $r_0=0.05$ or $0.1$, the HK relaxation can lead to a tighter upper bound on the training optimal value than our moment relaxation.
		However, for the smaller radius $r_0=0.01$, there is no significant difference between the two models.
		The two models also yield decisions with similar out-of-sample test performance.
	\end{itemize}
	
	\begin{table}[htb!]
		\centering
		\scriptsize
		\begin{tabular}{rrrrrrrrrr}
			\toprule
			& & \multicolumn{4}{c}{Moment-WDRO} & \multicolumn{4}{c}{Hanasusanto-Kuhn} \\
			\cmidrule(lr){3-6} \cmidrule(lr){7-10}
			& & \multicolumn{2}{c}{Training} & \multicolumn{2}{c}{Test} & \multicolumn{2}{c}{Training} & \multicolumn{2}{c}{Test}\\
			\cmidrule(lr){3-4} \cmidrule(lr){5-6} \cmidrule(lr){7-8} \cmidrule(lr){9-10}
			$N$ & $r$ & Time (s) & Obj. & Mean & Std.
			& Time (s) & Obj. & Mean & Std. \\
			\midrule
			10 & 0.000 & 13.4 & $-10.24$ & $-9.85$ & $2.39$ & 457.0 & $-10.13$ & $-9.88$ & $2.32$ \\
			10 & 0.010 & 69.5 & $-10.07$ & $-9.84$ & $2.29$ & 434.2 & $-10.12$ & $-9.88$ & $2.31$ \\
			10 & 0.050 & 46.1 & $-9.40$ & $-9.87$ & $2.27$ & 527.8 & $-9.75$ & $-9.91$ & $2.19$ \\
			10 & 0.100 & 46.9 & $-8.55$ & $-9.91$ & $2.25$ & 574.9 & $-9.26$ & $-9.87$ & $2.13$ \\
			\midrule
			20 & 0.000 & 3.1 & $-9.80$ & $-9.85$ & $2.27$ & 844.1 & $-9.69$ & $-9.91$ & $2.25$ \\
			20 & 0.010 & 74.8 & $-9.63$ & $-9.88$ & $2.21$ & 806.7 & $-9.63$ & $-9.91$ & $2.25$ \\
			20 & 0.048 & 83.7 & $-9.02$ & $-9.91$ & $2.21$ & 1088.9 & $-9.35$ & $-9.90$ & $2.23$ \\
			20 & 0.095 & 85.9 & $-8.22$ & $-9.87$ & $2.12$ & 890.4 & $-8.88$ & $-9.86$ & $2.07$ \\
			\midrule
			30 & 0.000 & 9.4 & $-9.81$ & $-9.91$ & $2.33$ & 1520.7 & $-9.63$ & $-9.92$ & $2.26$ \\
			30 & 0.009 & 134.2 & $-9.67$ & $-9.92$ & $2.22$ & 1382.8 & $-9.63$ & $-9.92$ & $2.26$ \\
			30 & 0.046 & 130.6 & $-9.05$ & $-9.89$ & $2.28$ & 1646.9 & $-9.37$ & $-9.91$ & $2.23$ \\
			30 & 0.093 & 135.1 & $-8.28$ & $-9.93$ & $2.18$ & 1078.7 & $-8.89$ & $-9.90$ & $2.13$ \\
			\midrule
			40 & 0.000 & 5.2 & $-9.92$ & $-9.89$ & $2.24$ & 1755.6 & $-9.73$ & $-9.93$ & $2.26$ \\
			40 & 0.009 & 168.4 & $-9.78$ & $-9.89$ & $2.29$ & 1581.7 & $-9.73$ & $-9.93$ & $2.26$ \\
			40 & 0.046 & 158.9 & $-9.17$ & $-9.89$ & $2.24$ & 1951.4 & $-9.49$ & $-9.93$ & $2.20$ \\
			40 & 0.091 & 183.3 & $-8.44$ & $-9.89$ & $2.15$ & 1506.5 & $-9.04$ & $-9.91$ & $2.13$ \\
			\midrule
			50 & 0.000 & 7.8 & $-9.66$ & $-9.91$ & $2.27$ & 2335.0 & $-9.44$ & $-9.95$ & $2.21$ \\
			50 & 0.009 & 212.5 & $-9.52$ & $-9.94$ & $2.28$ & 2142.2 & $-9.44$ & $-9.95$ & $2.21$ \\
			50 & 0.045 & 206.5 & $-8.94$ & $-9.93$ & $2.21$ & 2064.0 & $-9.23$ & $-9.94$ & $2.18$ \\
			50 & 0.090 & 217.3 & $-8.19$ & $-9.86$ & $2.05$ & 2030.9 & $-8.83$ & $-9.86$ & $1.96$ \\
			\midrule
			60 & 0.000 & 8.4 & $-9.56$ & $-9.92$ & $2.26$ & 2430.5 & $-9.33$ & $-9.95$ & $2.19$ \\
			60 & 0.009 & 212.0 & $-9.42$ & $-9.92$ & $2.17$ & 2428.5 & $-9.32$ & $-9.94$ & $2.18$ \\
			60 & 0.044 & 293.2 & $-8.86$ & $-9.95$ & $2.24$ & 2732.2 & $-9.08$ & $-9.92$ & $2.11$ \\
			60 & 0.089 & 367.5 & $-8.13$ & $-9.88$ & $2.06$ & 2843.3 & $-8.74$ & $-9.84$ & $1.93$ \\
			\midrule
			70 & 0.000 & 8.5 & $-9.50$ & $-9.93$ & $2.25$ & 2794.1 & $-9.25$ & $-9.94$ & $2.18$ \\
			70 & 0.009 & 260.8 & $-9.36$ & $-9.92$ & $2.15$ & 2573.0 & $-9.26$ & $-9.95$ & $2.19$ \\
			70 & 0.044 & 254.1 & $-8.79$ & $-9.90$ & $2.13$ & 2586.2 & $-9.05$ & $-9.89$ & $2.02$ \\
			70 & 0.088 & 438.7 & $-8.09$ & $-9.85$ & $1.99$ & 3615.2 & $-8.66$ & $-9.84$ & $1.92$ \\
			\midrule
			80 & 0.000 & 17.6 & $-9.56$ & $-9.90$ & $2.22$ & 5228.1 & $-9.35$ & $-9.94$ & $2.19$ \\
			80 & 0.009 & 482.7 & $-9.42$ & $-9.93$ & $2.24$ & 4076.3 & $-9.29$ & $-9.94$ & $2.18$ \\
			80 & 0.044 & 355.7 & $-8.85$ & $-9.89$ & $2.16$ & 2967.6 & $-9.08$ & $-9.92$ & $2.10$ \\
			80 & 0.087 & 656.1 & $-8.16$ & $-9.92$ & $2.12$ & 2550.2 & $-8.72$ & $-9.83$ & $1.92$ \\
			\midrule
			90 & 0.000 & 13.6 & $-9.60$ & $-9.92$ & $2.28$ & 3397.0 & $-9.35$ & $-9.95$ & $2.20$ \\
			90 & 0.009 & 308.9 & $-9.46$ & $-9.91$ & $2.15$ & 3587.2 & $-9.32$ & $-9.94$ & $2.19$ \\
			90 & 0.043 & 359.5 & $-8.91$ & $-9.92$ & $2.15$ & 3595.8 & $-9.12$ & $-9.93$ & $2.15$ \\
			90 & 0.086 & 350.7 & $-8.20$ & $-9.84$ & $1.96$ & 3092.3 & $-8.76$ & $-9.84$ & $1.94$ \\
			\bottomrule
		\end{tabular}
		\caption{Comparison with Hanasusanto-Kuhn semidefinite relaxation}
		\label{tab:HKComparison}
	\end{table}

\section{Concluding Remarks}
\label{sec:con}

In this work, we have investigated moment relaxations for $p$-WDRO problems.
The main contribution of our analysis is to show that such moment relaxations can preserve the asymptotic consistency property of $p$-WDRO model with polynomially sized semidefinite optimization formulations that admit direct implementation using standard SDP solvers.
To the best of our knowledge, this is the first time that convergence of moment relaxations is achieved by increasing the data size $N\to\infty$, instead of raising the relaxation order $k\to\infty$.
We believe that this has brought exciting perspectives on moment relaxations in the general context of data-driven global optimization.

The following questions can be of interest for future studies.
First, can we further improve the scalability of the proposed moment relaxations while keeping the asymptotic consistency guarantee?
Possible methods include sparse moment relaxations~\cite{wang2021tssos,ahmadi2019dsos,blekherman2022sparse}, and low-rank (Burer-Monteiro) methods~\cite{burer2003nonlinear,boumal2016non,cifuentes2021burer,legat2023low,blekherman2024spurious}.
Also we observe that the strengthened relaxation proposed in~Section~\ref{sec:2stage-strengthened} may sometimes run into numerical instability, and thus methods such as facial reduction~\cite{pataki2013strong,permenter2018partial} could be useful for more efficient implementation.
Second, can we extend the proposed moment relaxations to general (non-polynomial) cost functions?
The polynomial approximation technique in~\cite{zhong2024towards} may be considered for such extension.
	It would also be very interesting to extend the relaxation framework to include the case $p=1$, 
	which is very commonly used choice in the WDRO literature.

\vspace{0.1in}
\noindent{\bf Acknowledgements}
	Portions of this research were conducted with the advanced
	computing resources provided by Texas A\&M High Performance Research Computing.
	The authors used the large language models Claude Opus 4.7 and GPT-5.5 
	to assist with the implementation of the numerical experiments, 
	and independently verified the correctness of the implementation.

\printbibliography
\end{document}